
\magnification=1100
\overfullrule0pt

\input amssym.def
\input prepictex
\input pictex
\input postpictex



\def\starsub#1#2#3#4
{{\beginpicture
\put{$\scriptstyle{#1}$} at 5 13
\put{$\scriptstyle{\downarrow}$} at 5 6
\put{$\scriptstyle{#2}$} at 5 0
\put{$\scriptstyle{#3}^{{}^\swarrow}$} at -3 -6
\put{${\phantom{i}}^{{}_\searrow}\scriptstyle{#4}$} at 13 -6
\endpicture}}

\def\trisub#1#2#3{
{\beginpicture
\put{$\scriptstyle{#1}$} at 0 2
\put{$\scriptstyle{#2}$} at 7 2
\put{$\scriptstyle{#3}$} at 3 -5
\endpicture}
}
\def\vsub#1#2{
{\beginpicture
\put{$\scriptstyle{#1}$} at 3 2
\put{$\scriptstyle{#2}$} at 3 -5
\endpicture}
}

\def\dtrisub#1#2#3#4#5{
{\beginpicture
\put{$\scriptstyle{#1}$} at 0 2
\put{$\scriptstyle{#2}$} at 7 2
\put{$\scriptstyle{#3}$} at 0 -5
\put{$\scriptstyle{#4}$} at 7 -5
\put{$\scriptstyle{#5}$} at 3 -12
\endpicture}
}

\def\ddtrisub#1#2#3#4#5#6#7{
{\beginpicture
\put{$\scriptstyle{#1}$} at 0 2
\put{$\scriptstyle{#2}$} at 7 2
\put{$\scriptstyle{#3}$} at 0 -5
\put{$\scriptstyle{#4}$} at 7 -5
\put{$\scriptstyle{#5}$} at 0 -12
\put{$\scriptstyle{#6}$} at 7 -12
\put{$\scriptstyle{#7}$} at 3 -19
\endpicture}
}

\def\CC{{\Bbb C}}
\def\FF{{\Bbb F}}
\def\KK{{\Bbb K}}

\def\ZZ{{\Bbb Z}}

\def\Card{{\rm Card}}
\def\End{{\rm End}}
\def\Hom{{\rm Hom}}
\def\Id{{\rm Id}}
\def\im{{\rm im}}

\def\Ind{{\rm Ind}}
\def\Res{{\rm Res}}

\def\Tr{{\rm Tr}}
\def\tr{{\rm tr}}

\def\Rad{\hbox{Rad}}
\def\Card{\hbox{Card}}

\def\pkp{{p_{{k + {1 \over 2}}}}}

\def\Akp{{A_{k + {1 \over 2}}}}
\def\Akm{{A_{k - {1 \over 2}}}}

\def\deg{\hbox{deg}}
\def\End{\hbox{End}}
\def\Hom{\hbox{Hom}}

\def\Ind{\hbox{Ind}}

\def\Res{\hbox{Res}}

\def\Tr{\hbox{Tr}}

\newsymbol\ltimes 226E 
\newsymbol\rtimes 226F 


\def\qed{\hbox{\hskip 1pt\vrule width4pt height 6pt depth1.5pt \hskip 1pt}}

\def\sqr#1#2{{\vcenter{\vbox{\hrule height.#2pt
\hbox{\vrule width.#2pt height#1pt \kern#1pt
\vrule width.2pt}
\hrule height.2pt}}}}
\def\square{\mathchoice\sqr54\sqr54\sqr{3.5}3\sqr{2.5}3}


\font\smallcaps=cmcsc10
\font\titlefont=cmr10 scaled \magstep2

\font\sectionfont=cmbx10
\font\tinyrm=cmr10 at 8pt


\newcount\sectno
\newcount\subsectno
\newcount\resultno

\def\section #1. #2\par{
\sectno=#1
\resultno=0
\bigskip
\centerline{\sectionfont #1.  #2} }

\def\subsection #1\par{\bigskip\noindent{\it  #1} \medbreak}


\def\prop{ \global\advance\resultno by 1
\bigskip\noindent{\bf Proposition \the\sectno.\the\resultno. }\sl}
\def\lemma{ \global\advance\resultno by 1
\bigskip\noindent{\bf Lemma \the\sectno.\the\resultno. }
\sl}
\def\remark{ \global\advance\resultno by 1
\bigskip\noindent{\bf Remark \the\sectno.\the\resultno. }}
\def\example{ \global\advance\resultno by 1
\bigskip\noindent{\bf Example \the\sectno.\the\resultno. }}
\def\cor{ \global\advance\resultno by 1
\bigskip\noindent{\bf Corollary \the\sectno.\the\resultno. }\sl}
\def\thm{ \global\advance\resultno by 1
\bigskip\noindent{\bf Theorem \the\sectno.\the\resultno. }\sl}
\def\defn{ \global\advance\resultno by 1
\bigskip\noindent{\it Definition \the\sectno.\the\resultno. }\slrm}
\def\endthm{\rm\bigskip}

\def\endprop{\rm\bigskip}

\def\pf{\rm\bigskip\noindent{\it Proof. }}
\def\endpf{\qed\hfil\bigskip}
\def\pfend{\qed\hfil\bigskip}


\def\formula{\global\advance\resultno by 1
\eqno{(\the\sectno.\the\resultno)}}
\def\formulano{\global\advance\resultno by 1 (\the\sectno.\the\resultno)}
\def\tableno{\global\advance\resultno by 1
\the\sectno.\the\resultno. }
\def\lformula{\global\advance\resultno by 1
\leqno(\the\sectno.\the\resultno)}


\def\mapright#1{\smash{\mathop
        {\longrightarrow}\limits^{#1}}}

\def\mapleftright#1{\smash{
   \mathop{\longleftrightarrow}\limits^{#1}}}


\nopagenumbers
\def\runningtitle{\smallcaps partition algebras }
\headline={\ifnum\pageno>1\eoheadline\else\firstheadline\fi}
\def\names{\smallcaps t.\ halverson \enspace
and \enspace a.\ ram}
\def\firstheadline{}
\def\eoheadline{\ifodd\pageno\oddheadline\else\evenheadline\fi}
\def\oddheadline{\tenrm\hfil\runningtitle\hfil\folio}
\def\evenheadline{\tenrm\folio\hfil{\names}\hfil}

\def\monthname {\ifcase\month\or January\or February\or March\or April\or
May\or June\or
July\or August\or September\or October\or November\or December\fi}

\newcount\mins  \newcount\hours  \hours=\time \mins=\time
\def\now{\divide\hours by60 \multiply\hours by60 \advance\mins by-\hours
     \divide\hours by60         
     \ifnum\hours>12 \advance\hours by-12
       \number\hours:\ifnum\mins<10 0\fi\number\mins\ P.M.\else
       \number\hours:\ifnum\mins<10 0\fi\number\mins\ A.M.\fi}

\vphantom{$ $}  
\vskip.75truein

\centerline{\titlefont Partition Algebras}
\bigskip\bigskip

\centerline{
\vbox{
\centerline{Tom Halverson${}^\dagger$}
\centerline{Mathematics and Computer Science}
\centerline{Macalester College}
\centerline{Saint Paul, MN 55105}
\centerline{\tt halverson@macalester.edu}
}
\hskip-3.75truein
\vbox{
\centerline{ Arun Ram${}^\ast$}
\centerline{Department of Mathematics}
\centerline{University of Wisconsin--Madison}
\centerline{Madison, WI 53706}
\centerline{{\tt ram@math.wisc.edu}}
}
}

\footnote{}{\tinyrm ${}^\dagger$
Research supported in part
by National Science Foundation Grant DMS-0100975.}
\footnote{}{\tinyrm ${}^\ast$
Research supported in part by the
National Science Foundation (DMS-0097977) and the National Security
Agency (MDA904-01-1-0032).}

\bigskip


\bigskip\bigskip

\section 0. Introduction

\bigskip\noindent
A centerpiece of representation theory is the Schur-Weyl duality,
which says that,
\medskip\itemitem{(a)} the general linear group $GL_n(\CC)$ and
the symmetric group $S_k$ both act on tensor space
$$
V^{\otimes k} = \underbrace{V\otimes\cdots\otimes V}_{k{\rm \  factors}},
\qquad\hbox{ with } \quad \dim(V) = n,
$$
\smallskip\itemitem{(b)}  these two actions commute and
\smallskip\itemitem{(c)} each action generates the full centralizer of the
other, so that
\smallskip\itemitem{(d)} as a $(GL_n(\CC),S_k)$-bimodule, the tensor space 
has
a multiplicity free decomposition,
$$V^{\otimes k}\cong \bigoplus_\lambda L_{GL_n}(\lambda)\otimes S_k^\lambda,
\formula$$
where the $L_{GL_n}(\lambda)$ are irreducible $GL_n(\CC)$-modules
and the $S_k^\lambda$ are irreducible $S_k$-modules.
\medskip\noindent
The decomposition in (0.1) essentially makes the study
of the representations of $GL_n(\CC)$ and the study of representations
of the symmetric group $S_k$ two sides of the same coin.

The group $GL_n(\CC)$ has interesting subgroups,
$$\matrix{
&&GL_n(\CC)&\supseteq &O_n(\CC)&\supseteq &S_n&\supseteq &S_{n-1}, \cr
\hbox{and corresponding centralizer algebras,}\hfill  & \cr
&&\CC S_k &\subseteq &\CC B_k(n) &\subseteq &\CC A_k(n)&\subseteq
&\CC A_{k+{1\over2}}(n), \cr}$$
which are combinatorially defined in terms of the
``multiplication of diagrams'' (see Section 1) and
which play exactly analogous
``Schur-Weyl duality'' roles with their corresponding subgroup
of $GL_n(\CC)$.  The Brauer algebras $\CC B_k(n)$ were introduced
in 1937 by R. Brauer [Bra].  The partition algebras
$\CC A_k(n)$ arose in the early 1990s in the work of
P.\ Martin [Ma1-4] and later, independently, in the work of V.\ Jones [Jo]. 
Martin
and Jones discovered the partiton algebra as a generalization of the 
Temperley-Lieb
algebra and the Potts model in statistical mechanics.  The partition 
algebras
$\CC A_{k+{1\over2}}(n)$ apear in [Ma4] and [MR], and their existence and 
importance
was pointed out to us by C.\ Grood [Gr].  In this paper we follow the method
of [Ma4] and show that if the algebras $\CC A_{k+{1\over2}}(n)$ are given 
the same stature as the algebras $A_k(n)$,
then well-known methods from the theory
of the ``basic construction'' (see Section 4) allow for easy analysis of
the whole tower of algebras
$$
\CC  A_0(n) \subseteq \CC  A_{1\over 2}(n) \subseteq \CC  A_1(n)  \subseteq
\CC A_{1 {1 \over 2}}(n)  \subseteq \cdots,
$$
all at once.

\medskip\noindent
Let $\ell\in {1\over2}\ZZ_{\ge 0}$.  In this paper we prove:
\smallskip\noindent
\item{(a)} A presentation by generators and relations
for the algebras $\CC A_\ell(n)$.
\smallskip\noindent
\item{(b)} $\CC A_\ell(n)$ has
$$\hbox{an ideal}\qquad \CC I_\ell(n),
\qquad\hbox{with}\qquad
{\CC A_\ell(n)\over \CC I_\ell(n)}\cong \CC S_\ell,
$$
such that $\CC I_\ell(n)$ is isomorphic to a ``basic construction''
(see Section 4).
Thus the structure of the ideal $\CC I_{\ell}(n)$ can be
analyzed with the general theory of the basic construction
and the structure of the quotient $\CC A_\ell(n)/(\CC I_\ell(n))$ follows
from the general theory of the representations of the symmetric group.
\medskip\noindent
\item{(c)} The algebras
$\CC A_\ell(n)$ are in ``Schur-Weyl duality'' with the symmetric groups
$S_n$ and $S_{n-1}$ on $V^{\otimes k}$.
\smallskip\noindent
\item{(d)} The general theory of the basic construction
provides a construction of ``Specht modules'' for the
partition algebras, i.e. integral lattices
in the (generically) irreducible $\CC A_\ell(n)$-modules.
\smallskip\noindent
\item{(e)} Except for a few special cases,
the algebras $\CC A_{\ell}(n)$ are semisimple if and only if
$\ell\le (n+1)/2$.
\smallskip\noindent
\item{(f)} There are ``Murphy elements'' $M_i$ for the partition
algebras that play exactly analogous roles to the classical
Murphy elements for the group algebra of the symmetric group.
In particular, the $M_i$ commute with each other in $\CC A_\ell(n)$, and 
when
$\CC A_\ell(n)$ is semisimple each irreducible $\CC A_\ell(n)$-module has a
unique, up to constants, basis of simultaneous eigenvectors for the $M_i$.

\medskip
The primary new results in this paper are (a) and (f).
There has been work towards a presentation theorem for the
partition monoid by Fitzgerald and Leech [FL], and it is possible that
by now they have proved a similar presentation theorem.
The statement in (b) has appeared implicitly and explicitly
throughout the literature on the partition algebra, depending
on what one considers as the definition of a ``basic construction''.
The treatment of this connection between the partition algebras
and the basic construction is explained very nicely
and thoroughly in [Ma4].  We consider this connection 
an important part of the understanding
of the structure of the partition algebras.
The Schur-Weyl duality for the partition algebras $\CC A_k(n)$
appears in [Ma1], [Ma4], and [MR] and was one of the motivations
for the introduction of these algebras in [Jo].   The Schur-Weyl
duality for  $\CC A_{k+{1\over2}}(n)$ appears in [Ma4] and [MW].
Most of the previous literature (for example [Ma3], [MW1-2],  [DW])
 on the partition algebras has studied the structure of the partition 
algebras using the ``Specht'' modules of (d).  Our point here is that
their existence follows from the general theory of the
basic construction.  This is a special case of the fact
that quasi-hereditary algebras are iterated sequences of basic 
constructions, as proved by Dlab and Ringel [DR].  
The statements about the semsimplicity of $\CC A_{\ell}(n)$ have mostly,
if not completely, appeared in the work of Martin and Saleur [Ma3], [MS]. 
The Murphy elements for the partition algebras are new. Their form was 
conjectured by Owens [Ow], who proved that  the sum of the first $k$ of 
them is a central element in $\CC A_k(n)$.  Here we prove all of Owens' 
theorems and conjectures (by a different technique than he was using).
We have not taken the next natural step and
provided formulas for the action of the generators of the
partition algebra in the ``seminormal'' representations.
We hope that someone will do this in the near future.

The ``basic construction'' is a fundamental tool in the study of
algebras such as the partition algebra.  Of course, like any fundamental
construct, it appears in the literature and is rediscovered
over and over in various forms.  For example, one finds this
construction in Bourbaki [Bou1, Ch.\ 2, \S 4.2 Remark 1], in [Bro1-2],
in [GHJ, Ch. 2], and in the wonderful paper of
Dlab and Ringel [DR] where it is explained
that this construction is also the algebraic construct
that ``controls'' the theory of quasi-hereditary algebras, recollement
and highest weight categories [CPS] and some aspects of
the theory of perverse sheaves [MiV].  

Though this paper contains new results 
in the study of partition algebras we have made a
distinct effort to present this material in a ``survey'' style
so that it may be accessible to nonexperts and to newcomers
to the field.  For this reason we have included, in 
Sections 4 and 5, expositions, from scratch, of
\smallskip\noindent
\item{(a)} the theory of the
basic construction (see also [GHJ, Ch.\ 2]), and
\item{(b)} the theory of semisimple algebras, in particular,
Maschke's theorem, the Artin-Wedderburn theorem and the
Tits deformation theorem (see also [CR, \S 3B and \S 68]).
\smallskip\noindent
Here the reader will find statements of the main theorems which are
in exactly the correct form for our applications (generally difficult
to find in the literature), and short slick proofs of all the results on
the basic construction and on 
semisimple algebras that we need for the study of the partition algebras.

There are two sets of results on partition algebras that
we have not had the space to treat in this paper:
\item{(a)} The ``Frobenius formula," ``Murnaghan-Nakayama'' rule, and
 orthogonality rule for the irreducible characters
given by Halverson [Ha] and  Farina-Halverson [FH], and
\item{(b)} The cellularity of the partition algebras proved
by Xi [Xi] (see also Doran and Wales [DW]).

\smallskip
The techniques in this paper apply, in exactly the same
fashion, to the study of other diagram algebras;
in particular, the planar partition algebras $\CC P_k(n)$, the
Temperley-Lieb algebras $\CC T_k(n)$, and the
Brauer algebras $\CC B_k(n)$.  It was our original intent
to include in this paper results (mostly known)
for these algebras analogous to those which we have proved
for the algebras $\CC A_\ell(n)$, but the restrictions
of time and space have prevented this.
While perusing this paper, the reader should keep in mind  that the
techniques we have used do apply to these other algebras.

\bigskip

\section 1. The Partition Monoid

\bigskip
For $k \in \ZZ_{> 0}$, let
$$
\eqalign{
A_k & =\big\{ \hbox{set partitions of
$\{1,2,\ldots,k,1',2',\ldots, k'\}$}\big\},
\qquad\hbox{and}\cr
\Akp & =
\left\{ d \in A_{k+1} ~ \big| ~
\hbox{$(k+1)$ and $(k+1)'$ are in the same block}
\right\}.
\cr}
\formula
$$
The {\it propagating number} of $d \in  A_k$ is
$$
pn(d) = \left(
\matrix{
\hbox{the number of blocks in $d$ that contain both an element}\cr
\hbox{of $\{1, 2, \ldots, k\}$ and an element of
$\{1',2',\ldots,k'\}$}\hfill \cr}\right).
\formula
$$
For convenience, represent a set partition $d\in A_k$
by a graph with $k$ vertices in the top row,
labeled $1, \ldots,k$ left to right, and $k$ vertices
in the bottom row, labeled $1', \ldots, k'$ left to right,
with vertex $i$ and vertex $j$ connected by a path if $i$
and $j$ are in the same block of the set partition $d$.
For example,
$$
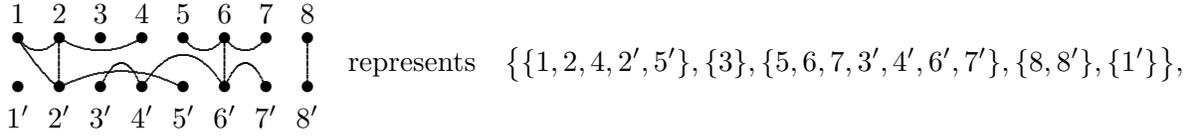
{\beginpicture
\setcoordinatesystem units <0.5cm,0.2cm> 
\setplotarea x from 0 to 7, y from 0 to 3    
\linethickness=0.5pt                        
\put{1} at 0 3.5
\put{2} at 1 3.5
\put{3} at 2 3.5
\put{4} at 3 3.5
\put{5} at 4 3.5
\put{6} at 5 3.5
\put{7} at 6 3.5
\put{8} at 7 3.5
\put{$1'$} at  0 -2.75
\put{$2'$} at 1 -2.75
\put{$3'$} at 2 -2.75
\put{$4'$} at 3 -2.75
\put{$5'$} at 4 -2.75
\put{$6'$} at 5 -2.75
\put{$7'$} at 6 -2.75
\put{$8'$} at 7 -2.75
\put{$\bullet$} at 0 -1 \put{$\bullet$} at 0 2
\put{$\bullet$} at 1 -1 \put{$\bullet$} at 1 2
\put{$\bullet$} at 2 -1 \put{$\bullet$} at 2 2
\put{$\bullet$} at 3 -1 \put{$\bullet$} at 3 2
\put{$\bullet$} at 4 -1 \put{$\bullet$} at 4 2
\put{$\bullet$} at 5 -1 \put{$\bullet$} at 5 2
\put{$\bullet$} at 6 -1 \put{$\bullet$} at 6 2
\put{$\bullet$} at 7 -1 \put{$\bullet$} at 7 2
\plot 1 -1 1 2  /
\plot 5 -1 5 2  /
\plot 7 -1 7 2  /
\setquadratic
\plot 0  2  .5 1.25 1 2 /
\plot 4  2 4.5 1.25 5 2 /
\plot 5  2 5.5 1.25 6 2 /
\plot 0  2  .5 .25 1 -1 /
\plot 1 -1 2.5 0 4 -1 /
\plot 1  2  2  1.25 3 2 /
\plot 2 -1 2.5 0.5 3 -1 /
\plot 3 -1 4 1 5 -1 /
\plot 5 -1 5.5 .5 6 -1 /
\endpicture}
\quad\hbox{represents}\quad
\big\{ \{ 1, 2, 4, 2',5'\},\{3\},
\{5, 6, 7, 3', 4', 6', 7'\},\{8, 8'\}, \{1'\}\big\},
$$
and has propagating number 3.
The graph representing $d$ is not unique.

Define the composition $d_1\circ d_2$ of partition diagrams
$d_1,d_2\in A_k$ to be the
set partition $d_1 \circ d_2\in A_k$
obtained by placing $d_1$ above $d_2$ and identifying the bottom
dots of $d_1$ with the top dots of $d_2$,
removing any connected components that live entirely in the middle row.
For example,
$$\hbox{if}\qquad
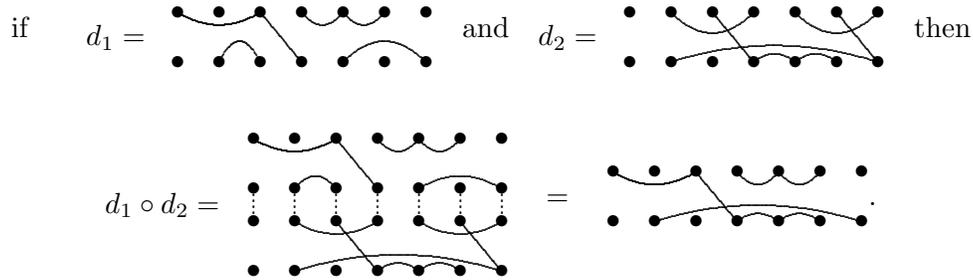
{\beginpicture
\setcoordinatesystem units <0.5cm,0.2cm> 
\setplotarea x from 0 to 7, y from 0 to 3    
\linethickness=0.5pt
\put{$d_1 = $} at -.5 0
\put{$\bullet$} at 1 -1.5 \put{$\bullet$} at 1 1.5
\put{$\bullet$} at 2 -1.5 \put{$\bullet$} at 2 1.5
\put{$\bullet$} at 3 -1.5 \put{$\bullet$} at 3 1.5
\put{$\bullet$} at 4 -1.5 \put{$\bullet$} at 4 1.5
\put{$\bullet$} at 5 -1.5 \put{$\bullet$} at 5 1.5
\put{$\bullet$} at 6 -1.5 \put{$\bullet$} at 6 1.5
\put{$\bullet$} at 7 -1.5 \put{$\bullet$} at 7 1.5
\setquadratic
\plot 1 1.5 2 0.75 3 1.5 /
\plot 1 1.5 2 0.75 3 1.5 /
\plot 4 1.5 4.5 0.75 5 1.5 /
\plot 5 1.5 5.5 0.75 6 1.5 /
\plot 2 -1.5 2.5 -0.25 3 -1.5 /
\plot 5 -1.5 6 -0.25 7 -1.5 /
\setlinear
\plot 3 1.5 4 -1.5 /
\endpicture}
\quad\hbox{and}\quad
{\beginpicture
\setcoordinatesystem units <0.5cm,0.2cm> 
\setplotarea x from 0 to 7, y from 0 to 3    
\linethickness=0.5pt
\put{$d_2 = $} at -.5 0
\put{$\bullet$} at 1 1.5 \put{$\bullet$} at 1 -1.5
\put{$\bullet$} at 2 1.5 \put{$\bullet$} at 2 -1.5
\put{$\bullet$} at 3 1.5 \put{$\bullet$} at 3 -1.5
\put{$\bullet$} at 4 1.5 \put{$\bullet$} at 4 -1.5
\put{$\bullet$} at 5 1.5 \put{$\bullet$} at 5 -1.5
\put{$\bullet$} at 6 1.5 \put{$\bullet$} at 6 -1.5
\put{$\bullet$} at 7 1.5 \put{$\bullet$} at 7 -1.5

\setquadratic
\plot 2 1.5 3 0.25 4 1.5 /
\plot 5 1.5 6 0.25 7 1.5 /
\plot 2 -1.5 4.5 -0.5 7 -1.5 /
\plot 4 -1.5 4.5 -1 5 -1.5 /
\plot 5 -1.5 5.5 -1 6 -1.5 /
\setlinear
\plot 3 1.5 4 -1.5 /
\plot 6 1.5 7 -1.5 /
\endpicture}
\quad\hbox{then}
$$
\medskip
$$
{\beginpicture
\setcoordinatesystem units <0.5cm,0.2cm> 
\setplotarea x from 0 to 7, y from 0 to 3    
\linethickness=0.5pt
\put{$d_1\circ d_2 = $\qquad} at -.5 0
\put{$\bullet$} at 1 1 \put{$\bullet$} at 1 4
\put{$\bullet$} at 2 1 \put{$\bullet$} at 2 4
\put{$\bullet$} at 3 1 \put{$\bullet$} at 3 4
\put{$\bullet$} at 4 1 \put{$\bullet$} at 4 4
\put{$\bullet$} at 5 1 \put{$\bullet$} at 5 4
\put{$\bullet$} at 6 1 \put{$\bullet$} at 6 4
\put{$\bullet$} at 7 1 \put{$\bullet$} at 7 4

\put{$\bullet$} at 1 -1 \put{$\bullet$} at 1 -4
\put{$\bullet$} at 2 -1 \put{$\bullet$} at 2 -4
\put{$\bullet$} at 3 -1 \put{$\bullet$} at 3 -4
\put{$\bullet$} at 4 -1 \put{$\bullet$} at 4 -4
\put{$\bullet$} at 5 -1 \put{$\bullet$} at 5 -4
\put{$\bullet$} at 6 -1 \put{$\bullet$} at 6 -4
\put{$\bullet$} at 7 -1 \put{$\bullet$} at 7 -4

\setquadratic
\plot 1 4 2 3.25 3 4 /
\plot 1 4 2 3.25 3 4 /
\plot 4 4 4.5 3.25 5 4 /
\plot 5 4 5.5 3.25 6 4 /
\plot 2 1 2.5 1.75 3 1 /
\plot 5 1 6 1.75 7 1 /
\plot 2 -1 3 -1.75 4 -1 /
\plot 5 -1 6 -1.75 7 -1 /
\plot 2 -4 4.5 -3 7 -4 /
\plot 4 -4 4.5 -3.5 5 -4 /
\plot 5 -4 5.5 -3.5 6 -4 /
\setlinear
\plot 3 4 4 1 /
\plot 3 -1 4 -4 /
\plot 6 -1 7 -4 /
\setdashpattern <.02cm,.05cm>
\plot 1 1 1 -1 /
\plot 2 1 2 -1 /
\plot 3 1 3 -1 /
\plot 4 1 4 -1 /
\plot 5 1 5 -1 /
\plot 6 1 6 -1 /
\plot 7 1 7 -1 /
\endpicture}
\quad = \quad {\beginpicture
\setcoordinatesystem units <0.5cm,0.2cm> 
\setplotarea x from 1 to 7, y from 0 to 3    
\linethickness=0.5pt                        
\put{$\bullet$} at 1 -1 \put{$\bullet$} at 1 2
\put{$\bullet$} at 2 -1 \put{$\bullet$} at 2 2
\put{$\bullet$} at 3 -1 \put{$\bullet$} at 3 2
\put{$\bullet$} at 4 -1 \put{$\bullet$} at 4 2
\put{$\bullet$} at 5 -1 \put{$\bullet$} at 5 2
\put{$\bullet$} at 6 -1 \put{$\bullet$} at 6 2
\put{$\bullet$} at 7 -1 \put{$\bullet$} at 7 2
\plot 3 2 4 -1 /
\setquadratic
\plot 1 2 2 1.25 3 2 /
\plot 1 2 2 1.25 3 2 /
\plot 4 2 4.5 1.25 5 2 /
\plot 5 2 5.5 1.25 6 2 /
\plot 2 -1 4.5 0 7 -1 /
\plot 4 -1 4.5 -.5 5 -1 /
\plot 5 -1 5.5 -.5 6 -1 /
\endpicture}.
$$
Diagram multiplication makes  $A_k$ into an
associative monoid with identity,
$1 =
{\beginpicture
\setcoordinatesystem units <0.4cm,0.12cm> 
\setplotarea x from 1 to 4, y from -1 to 2   
\linethickness=0.5pt                        
\put{$\bullet$} at 1 -1 \put{$\bullet$} at 1 2 \plot 1 -1 1 2 /
\put{$\bullet$} at 2 -1 \put{$\bullet$} at 2 2 \plot 2 -1 2 2 /
\put{$\cdots$} at 3 .5
\put{$\bullet$} at 4 -1  \put{$\bullet$} at 4 2
\plot 4 -1 4 2 /
\endpicture}.
$
The propagating number satisfies
$$
pn(d_1 \circ d_2) \le \min(pn(d_1),pn(d_2)). \formula
$$

A set partition is {\it planar} [Jo] if it can be represented
as a graph without edge crossings inside of the rectangle formed
by its vertices.  For each $k \in {1 \over 2} \ZZ_{>0}$,
the following are
submonoids of the partition monoid $A_k$:
$$\matrix{
S_k = \{ d \in A_k \ | \ pn(d) = k \},\quad
I_k = \{ d \in A_k \ | \ pn(d) < k \}, \quad
P_k = \{ d \in A_k \ | \ \hbox{$d$ is planar} \}, \cr
\cr
B_k = \{ d\in A_k \ |\ \hbox{all blocks of $d$ have size 2}\},
\quad\hbox{and}\quad
T_k = P_k \cap B_k. \cr }
\formula
$$
Examples are
$$
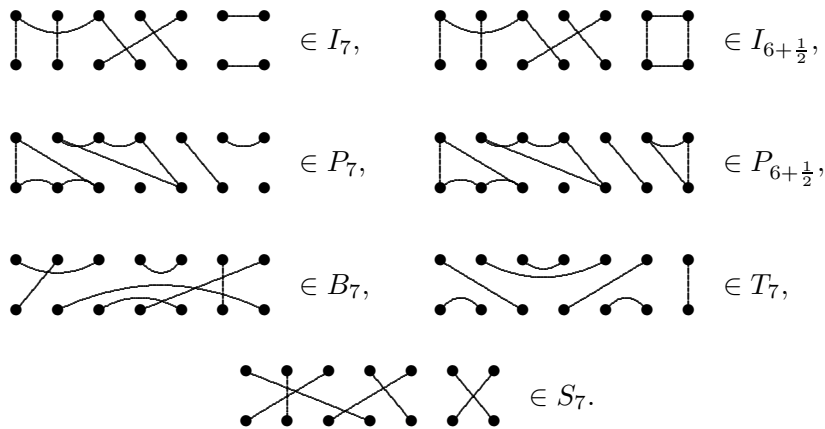
\matrix{
{\beginpicture
\setcoordinatesystem units <0.5cm,0.2cm> 
\setplotarea x from 0 to 7, y from 0 to 3    
\linethickness=0.5pt                        
\put{$\bullet$} at 1 -1 \put{$\bullet$} at 1 2
\put{$\bullet$} at 2 -1 \put{$\bullet$} at 2 2
\put{$\bullet$} at 3 -1 \put{$\bullet$} at 3 2
\put{$\bullet$} at 4 -1 \put{$\bullet$} at 4 2
\put{$\bullet$} at 5 -1 \put{$\bullet$} at 5 2
\put{$\bullet$} at 6 -1 \put{$\bullet$} at 6 2
\put{$\bullet$} at 7 -1 \put{$\bullet$} at 7 2
\plot 3 2 4 -1 /
\plot 2 2 2 -1 /
\plot 1 2 1 -1 /
\plot 4 2 5 -1 /
\plot 5 2 3 -1 /
\plot 6 -1 7 -1 /
\plot 7 2 6 2 /
\setquadratic
\plot 1 2 2 1 3 2 /
\endpicture}
\ \ \in  I_7, \hfill
&{\beginpicture
\setcoordinatesystem units <0.5cm,0.2cm> 
\setplotarea x from 0 to 7, y from 0 to 3    
\linethickness=0.5pt                        
\put{$\bullet$} at 1 -1 \put{$\bullet$} at 1 2
\put{$\bullet$} at 2 -1 \put{$\bullet$} at 2 2
\put{$\bullet$} at 3 -1 \put{$\bullet$} at 3 2
\put{$\bullet$} at 4 -1 \put{$\bullet$} at 4 2
\put{$\bullet$} at 5 -1 \put{$\bullet$} at 5 2
\put{$\bullet$} at 6 -1 \put{$\bullet$} at 6 2
\put{$\bullet$} at 7 -1 \put{$\bullet$} at 7 2
\plot 3 2 4 -1 /
\plot 2 2 2 -1 /
\plot 1 2 1 -1 /
\plot 4 2 5 -1 /
\plot 5 2 3 -1 /
\plot 6 -1 7 -1 /
\plot 7 2 6 2 /
\plot 7 2 7 -1 /
\plot 6 2 6 -1 /
\setquadratic
\plot 1 2 2 1 3 2 /
\endpicture}
\ \ \in  I_{6 + {1\over2}}, \hfill \cr
\cr
\cr
{\beginpicture
\setcoordinatesystem units <0.5cm,0.2cm> 
\setplotarea x from 0 to 7, y from 0 to 3    
\linethickness=0.5pt                        
\put{$\bullet$} at 1 -1 \put{$\bullet$} at 1 2
\put{$\bullet$} at 2 -1 \put{$\bullet$} at 2 2
\put{$\bullet$} at 3 -1 \put{$\bullet$} at 3 2
\put{$\bullet$} at 4 -1 \put{$\bullet$} at 4 2
\put{$\bullet$} at 5 -1 \put{$\bullet$} at 5 2
\put{$\bullet$} at 6 -1 \put{$\bullet$} at 6 2
\put{$\bullet$} at 7 -1 \put{$\bullet$} at 7 2
\plot 1 2 3 -1 /
\plot 5 2 6 -1 /
\plot 2 2 5 -1 /
\plot 4 2 5 -1 /
\plot 1 2 1 -1 /
\setquadratic
\plot 2 2 2.5 1.5 3 2 /
\plot 3 2 3.5 1.5 4 2 /
\plot 6 2 6.5 1.5 7 2 /
\plot 1 -1 1.5 -.5 2 -1 /
\plot 2 -1 2.5 -.5 3 -1 /
\endpicture}
\ \ \in P_{7}, \hfill
&{\beginpicture
\setcoordinatesystem units <0.5cm,0.2cm> 
\setplotarea x from 0 to 7, y from 0 to 3    
\linethickness=0.5pt                        
\put{$\bullet$} at 1 -1 \put{$\bullet$} at 1 2
\put{$\bullet$} at 2 -1 \put{$\bullet$} at 2 2
\put{$\bullet$} at 3 -1 \put{$\bullet$} at 3 2
\put{$\bullet$} at 4 -1 \put{$\bullet$} at 4 2
\put{$\bullet$} at 5 -1 \put{$\bullet$} at 5 2
\put{$\bullet$} at 6 -1 \put{$\bullet$} at 6 2
\put{$\bullet$} at 7 -1 \put{$\bullet$} at 7 2
\plot 1 2 3 -1 /
\plot 7 2 7 -1 /
\plot 6 2 7 -1 /
\plot 5 2 6 -1 /
\plot 2 2 5 -1 /
\plot 4 2 5 -1 /
\plot 1 2 1 -1 /
\setquadratic
\plot 2 2 2.5 1.5 3 2 /
\plot 3 2 3.5 1.5 4 2 /
\plot 6 2 6.5 1.5 7 2 /
\plot 1 -1 1.5 -.5 2 -1 /
\plot 2 -1 2.5 -.5 3 -1 /
\endpicture}
\ \ \in P_{6+{1 \over 2}}, \hfill \cr
\cr\cr
{\beginpicture
\setcoordinatesystem units <0.5cm,0.2cm> 
\setplotarea x from 0 to 7, y from 0 to 3    
\linethickness=0.5pt                        
\put{$\bullet$} at 1 -1 \put{$\bullet$} at 1 2
\put{$\bullet$} at 2 -1 \put{$\bullet$} at 2 2
\put{$\bullet$} at 3 -1 \put{$\bullet$} at 3 2
\put{$\bullet$} at 4 -1 \put{$\bullet$} at 4 2
\put{$\bullet$} at 5 -1 \put{$\bullet$} at 5 2
\put{$\bullet$} at 6 -1 \put{$\bullet$} at 6 2
\put{$\bullet$} at 7 -1 \put{$\bullet$} at 7 2
\plot 7 2 4 -1 /
\plot 6 2 6 -1 /
\plot 2 2 1 -1 /
\setquadratic
\plot 1 2 2 1.25 3 2 /
\plot 4 2 4.5 1.25 5 2 /
\plot 2 -1 4.5 0.5 7 -1 /
\plot 3 -1 4 -.25 5 -1 /
\endpicture}
\ \ \in  B_7, \hfill
&{\beginpicture
\setcoordinatesystem units <0.5cm,0.2cm> 
\setplotarea x from 0 to 7, y from 0 to 3    
\linethickness=0.5pt                        
\put{$\bullet$} at 1 -1 \put{$\bullet$} at 1 2
\put{$\bullet$} at 2 -1 \put{$\bullet$} at 2 2
\put{$\bullet$} at 3 -1 \put{$\bullet$} at 3 2
\put{$\bullet$} at 4 -1 \put{$\bullet$} at 4 2
\put{$\bullet$} at 5 -1 \put{$\bullet$} at 5 2
\put{$\bullet$} at 6 -1 \put{$\bullet$} at 6 2
\put{$\bullet$} at 7 -1 \put{$\bullet$} at 7 2
\plot 1 2 3 -1 /
\plot 6 2 4 -1 /
\plot 7 2 7 -1 /
\setquadratic
\plot 3 2 3.5 1.5 4 2 /
\plot 2 2 3.5 1.0 5 2 /
\plot 1 -1 1.5 -.25 2 -1 /
\plot 5 -1 5.5 -.25 6 -1 /
\endpicture}
\ \ \in T_7, \hfill \cr
}$$
$$
{\beginpicture
\setcoordinatesystem units <0.5cm,0.2cm> 
\setplotarea x from 0 to 7, y from 0 to 3    
\linethickness=0.5pt                        
\put{$\bullet$} at 1 -1 \put{$\bullet$} at 1 2
\put{$\bullet$} at 2 -1 \put{$\bullet$} at 2 2
\put{$\bullet$} at 3 -1 \put{$\bullet$} at 3 2
\put{$\bullet$} at 4 -1 \put{$\bullet$} at 4 2
\put{$\bullet$} at 5 -1 \put{$\bullet$} at 5 2
\put{$\bullet$} at 6 -1 \put{$\bullet$} at 6 2
\put{$\bullet$} at 7 -1 \put{$\bullet$} at 7 2
\plot 1 2 4 -1 /
\plot 2 2 2 -1 /
\plot 3 2 1 -1 /
\plot 4 2 5 -1 /
\plot 5 2 3 -1 /
\plot 6 2 7 -1 /
\plot 7 2 6 -1 /
\endpicture}
\ \ \in  S_7.
$$

For $k \in {1 \over 2} \ZZ_{>0}$, there is an isomorphism of monoids
$$P_k\ \  \mapleftright{{\rm 1-1}}\ \  T_{2k},\formula
$$
which is best illustrated by examples.
For $k=7$ we
have
$$
\matrix{
{\beginpicture
\setcoordinatesystem units <0.5cm,0.35cm> 
\setplotarea x from 1 to 7, y from -1 to 2    
\linethickness=0.5pt                        
\put{$\bullet$} at 1 -1 \put{$\bullet$} at 1 2
\put{$\bullet$} at 2 -1 \put{$\bullet$} at 2 2
\put{$\bullet$} at 3 -1 \put{$\bullet$} at 3 2
\put{$\bullet$} at 4 -1 \put{$\bullet$} at 4 2
\put{$\bullet$} at 5 -1 \put{$\bullet$} at 5 2
\put{$\bullet$} at 6 -1 \put{$\bullet$} at 6 2
\put{$\bullet$} at 7 -1 \put{$\bullet$} at 7 2
\plot 1 2 3 -1 /
\plot 5 2 6 -1 /
\plot 2 2 5 -1 /
\plot 4 2 5 -1 /
\plot 1 2 1 -1 /
\setquadratic
\plot 2 2 2.5 1.65 3 2 /
\plot 3 2 3.5 1.65 4 2 /
\plot 6 2 6.5 1.5 7 2 /
\plot 1 -1 1.5 -.5 2 -1 /
\plot 2 -1 2.5 -.5 3 -1 /
\endpicture}
&\quad\leftrightarrow\quad&
{\beginpicture
\setcoordinatesystem units <.6cm,0.35cm> 
\setplotarea x from 1 to 7, y from -1 to 2   
\linethickness=0.5pt                        
\put{$\bullet$} at 1 -1 \put{$\bullet$} at 1 2
\put{$\circ$} at .75 -1 \put{$\circ$} at .75 2
\put{$\circ$} at 1.25 -1 \put{$\circ$} at 1.25 2
\put{$\bullet$} at 2 -1 \put{$\bullet$} at 2 2
\put{$\circ$} at 1.75 -1 \put{$\circ$} at 1.75 2
\put{$\circ$} at 2.25 -1 \put{$\circ$} at 2.25 2
\put{$\bullet$} at 3 -1 \put{$\bullet$} at 3 2
\put{$\circ$} at 2.75 -1 \put{$\circ$} at 2.75 2
\put{$\circ$} at 3.25 -1 \put{$\circ$} at 3.25 2
\put{$\bullet$} at 4 -1 \put{$\bullet$} at 4 2
\put{$\circ$} at 3.75 -1 \put{$\circ$} at 3.75 2
\put{$\circ$} at 4.25 -1 \put{$\circ$} at 4.25 2
\put{$\bullet$} at 5 -1 \put{$\bullet$} at 5 2
\put{$\circ$} at 4.75 -1 \put{$\circ$} at 4.75 2
\put{$\circ$} at 5.25 -1 \put{$\circ$} at 5.25 2
\put{$\bullet$} at 6 -1 \put{$\bullet$} at 6 2
\put{$\circ$} at 5.75 -1 \put{$\circ$} at 5.75 2
\put{$\circ$} at 6.25 -1 \put{$\circ$} at 6.25 2
\put{$\bullet$} at 7 -1 \put{$\bullet$} at 7 2
\put{$\circ$} at 6.75 -1 \put{$\circ$} at 6.75 2
\put{$\circ$} at 7.25 -1 \put{$\circ$} at 7.25 2
\plot 1 2  3 -1 /
\plot 5 2 6 -1 /
\plot 2 2 5 -1 /
\plot 4 2 5 -1 /
\plot 1 2 1 -1 /
\setquadratic
\plot 2 2 2.5 1.65 3 2 /
\plot 3 2 3.5 1.65 4 2 /
\plot 6 2 6.5 1.5 7 2 /
\plot 1 -1 1.5 -.5 2 -1 /
\plot 2 -1 2.5 -.5 3 -1 /
\setdashes <.04cm>
\plot 2.25 2 2.5 1.75 2.75 2 /
\plot 3.25 2 3.5 1.75 3.75 2 /
\plot 6.25 2 6.5 1.75 6.75 2 /
\plot 5.75 2 6.5 1.25 7.25 2 /
\plot 6.75 -1 7 -.65 7.25 -1 /
\plot 3.75 -1 4 -.65 4.25 -1 /
\plot 1.25 -1 1.5 -.65 1.75 -1 /
\plot 2.25 -1 2.5 -.65 2.75 -1 /
\setlinear
\plot .75 2 .75 -1 /
\plot 1.25 2 3.25 -1 /
\plot 1.75 2 4.75 -1 /
\plot 4.25 2 5.25 -1 /
\plot 4.75 2 5.75 -1 /
\plot 5.25 2 6.25 -1 /
\endpicture}
&\quad\leftrightarrow\quad&
{\beginpicture
\setcoordinatesystem units <.5cm,0.35cm> 
\setplotarea x from 1 to 7, y from -1 to 2   
\linethickness=0.5pt                        
\put{$\circ$} at .75 -1 \put{$\circ$} at .75 2
\put{$\circ$} at 1.25 -1 \put{$\circ$} at 1.25 2
\put{$\circ$} at 1.75 -1 \put{$\circ$} at 1.75 2
\put{$\circ$} at 2.25 -1 \put{$\circ$} at 2.25 2
\put{$\circ$} at 2.75 -1 \put{$\circ$} at 2.75 2
\put{$\circ$} at 3.25 -1 \put{$\circ$} at 3.25 2
\put{$\circ$} at 3.75 -1 \put{$\circ$} at 3.75 2
\put{$\circ$} at 4.25 -1 \put{$\circ$} at 4.25 2
\put{$\circ$} at 4.75 -1 \put{$\circ$} at 4.75 2
\put{$\circ$} at 5.25 -1 \put{$\circ$} at 5.25 2
\put{$\circ$} at 5.75 -1 \put{$\circ$} at 5.75 2
\put{$\circ$} at 6.25 -1 \put{$\circ$} at 6.25 2
\put{$\circ$} at 6.75 -1 \put{$\circ$} at 6.75 2
\put{$\circ$} at 7.25 -1 \put{$\circ$} at 7.25 2
\setquadratic
\plot 2.25 2 2.5 1.75 2.75 2 /
\plot 3.25 2 3.5 1.75 3.75 2 /
\plot 6.25 2 6.5 1.75 6.75 2 /
\plot 5.75 2 6.5 1.25 7.25 2 /
\plot 6.75 -1 7 -.65 7.25 -1 /
\plot 3.75 -1 4 -.65 4.25 -1 /
\plot 1.25 -1 1.5 -.65 1.75 -1 /
\plot 2.25 -1 2.5 -.65 2.75 -1 /
\setlinear
\plot .75 2 .75 -1 /
\plot 1.25 2 3.25 -1 /
\plot 1.75 2 4.75 -1 /
\plot 4.25 2 5.25 -1 /
\plot 4.75 2 5.75 -1 /
\plot 5.25 2 6.25 -1 /
\endpicture}, \cr
}
$$
and for $k = 6 + {1 \over 2}$ we have
$$
\matrix{
{\beginpicture
\setcoordinatesystem units <0.5cm,0.35cm> 
\setplotarea x from 1 to 7, y from -1 to 2    
\linethickness=0.5pt                        
\put{$\bullet$} at 1 -1 \put{$\bullet$} at 1 2
\put{$\bullet$} at 2 -1 \put{$\bullet$} at 2 2
\put{$\bullet$} at 3 -1 \put{$\bullet$} at 3 2
\put{$\bullet$} at 4 -1 \put{$\bullet$} at 4 2
\put{$\bullet$} at 5 -1 \put{$\bullet$} at 5 2
\put{$\bullet$} at 6 -1 \put{$\bullet$} at 6 2
\put{$\bullet$} at 7 -1 \put{$\bullet$} at 7 2
\plot 1 2 3 -1 /
\plot 5 2 6 -1 /
\plot 2 2 5 -1 /
\plot 4 2 5 -1 /
\plot 1 2 1 -1 /
\plot 7 2 7 -1 /
\setquadratic
\plot 2 2 2.5 1.65 3 2 /
\plot 3 2 3.5 1.65 4 2 /
\plot 6 2 6.5 1.5 7 2 /
\plot 1 -1 1.5 -.5 2 -1 /
\plot 2 -1 2.5 -.5 3 -1 /
\endpicture}
&\quad\leftrightarrow\quad&
{\beginpicture
\setcoordinatesystem units <.6cm,0.35cm> 
\setplotarea x from 1 to 7, y from -1 to 2   
\linethickness=0.5pt                        
\put{$\bullet$} at 1 -1 \put{$\bullet$} at 1 2
\put{$\circ$} at .75 -1 \put{$\circ$} at .75 2
\put{$\circ$} at 1.25 -1 \put{$\circ$} at 1.25 2
\put{$\bullet$} at 2 -1 \put{$\bullet$} at 2 2
\put{$\circ$} at 1.75 -1 \put{$\circ$} at 1.75 2
\put{$\circ$} at 2.25 -1 \put{$\circ$} at 2.25 2
\put{$\bullet$} at 3 -1 \put{$\bullet$} at 3 2
\put{$\circ$} at 2.75 -1 \put{$\circ$} at 2.75 2
\put{$\circ$} at 3.25 -1 \put{$\circ$} at 3.25 2
\put{$\bullet$} at 4 -1 \put{$\bullet$} at 4 2
\put{$\circ$} at 3.75 -1 \put{$\circ$} at 3.75 2
\put{$\circ$} at 4.25 -1 \put{$\circ$} at 4.25 2
\put{$\bullet$} at 5 -1 \put{$\bullet$} at 5 2
\put{$\circ$} at 4.75 -1 \put{$\circ$} at 4.75 2
\put{$\circ$} at 5.25 -1 \put{$\circ$} at 5.25 2
\put{$\bullet$} at 6 -1 \put{$\bullet$} at 6 2
\put{$\circ$} at 5.75 -1 \put{$\circ$} at 5.75 2
\put{$\circ$} at 6.25 -1 \put{$\circ$} at 6.25 2
\put{$\bullet$} at 7 -1 \put{$\bullet$} at 7 2
\put{$\circ$} at 6.75 -1 \put{$\circ$} at 6.75 2
\plot 1 2  3 -1 /
\plot 5 2 6 -1 /
\plot 2 2 5 -1 /
\plot 4 2 5 -1 /
\plot 1 2 1 -1 /
\plot 7 2 7 -1 /
\setquadratic
\plot 2 2 2.5 1.65 3 2 /
\plot 3 2 3.5 1.65 4 2 /
\plot 6 2 6.5 1.5 7 2 /
\plot 1 -1 1.5 -.5 2 -1 /
\plot 2 -1 2.5 -.5 3 -1 /
\setdashes <.04cm>
\plot 2.25 2 2.5 1.75 2.75 2 /
\plot 3.25 2 3.5 1.75 3.75 2 /
\plot 6.25 2 6.5 1.75 6.75 2 /
\plot 3.75 -1 4 -.65 4.25 -1 /
\plot 1.25 -1 1.5 -.65 1.75 -1 /
\plot 2.25 -1 2.5 -.65 2.75 -1 /
\setlinear
\plot .75 2 .75 -1 /
\plot 1.25 2 3.25 -1 /
\plot 1.75 2 4.75 -1 /
\plot 4.25 2 5.25 -1 /
\plot 4.75 2 5.75 -1 /
\plot 5.25 2 6.25 -1 /
\plot 5.75 2 6.75 -1 /
\endpicture}
&\quad\leftrightarrow\quad&
{\beginpicture
\setcoordinatesystem units <.5cm,0.35cm> 
\setplotarea x from 1 to 7, y from -1 to 2   
\linethickness=0.5pt                        
\put{$\circ$} at .75 -1 \put{$\circ$} at .75 2
\put{$\circ$} at 1.25 -1 \put{$\circ$} at 1.25 2
\put{$\circ$} at 1.75 -1 \put{$\circ$} at 1.75 2
\put{$\circ$} at 2.25 -1 \put{$\circ$} at 2.25 2
\put{$\circ$} at 2.75 -1 \put{$\circ$} at 2.75 2
\put{$\circ$} at 3.25 -1 \put{$\circ$} at 3.25 2
\put{$\circ$} at 3.75 -1 \put{$\circ$} at 3.75 2
\put{$\circ$} at 4.25 -1 \put{$\circ$} at 4.25 2
\put{$\circ$} at 4.75 -1 \put{$\circ$} at 4.75 2
\put{$\circ$} at 5.25 -1 \put{$\circ$} at 5.25 2
\put{$\circ$} at 5.75 -1 \put{$\circ$} at 5.75 2
\put{$\circ$} at 6.25 -1 \put{$\circ$} at 6.25 2
\put{$\circ$} at 6.75 -1 \put{$\circ$} at 6.75 2
\setquadratic
\plot 2.25 2 2.5 1.75 2.75 2 /
\plot 3.25 2 3.5 1.75 3.75 2 /
\plot 6.25 2 6.5 1.75 6.75 2 /
\plot 3.75 -1 4 -.65 4.25 -1 /
\plot 1.25 -1 1.5 -.65 1.75 -1 /
\plot 2.25 -1 2.5 -.65 2.75 -1 /
\setlinear
\plot .75 2 .75 -1 /
\plot 1.25 2 3.25 -1 /
\plot 1.75 2 4.75 -1 /
\plot 4.25 2 5.25 -1 /
\plot 4.75 2 5.75 -1 /
\plot 5.25 2 6.25 -1 /
\plot 5.75 2 6.75 -1 /
\endpicture}. \cr
}
$$

Let $k \in \ZZ_{>0}$.
By permuting the vertices in the top row and in  the bottom row
each $d\in A_k$ can be written as a product
$d=\sigma_1 t\sigma_2$, with $\sigma_1,\sigma_2\in S_k$ and $t\in P_k$,
and so
$$A_k = S_{k}\, P_k\, S_{k}.\qquad\hbox{For example,}\qquad
\matrix{
&{\beginpicture
\setcoordinatesystem units <0.4cm,0.18cm>         
\setplotarea x from 0 to 7, y from 0 to 3    
\linethickness=0.5pt                         
\put{$\bullet$} at 0 0 \put{$\bullet$} at 0 3
\put{$\bullet$} at 1 0 \put{$\bullet$} at 1 3
\put{$\bullet$} at 2 0 \put{$\bullet$} at 2 3
\put{$\bullet$} at 3 0 \put{$\bullet$} at 3 3
\put{$\bullet$} at 4 0 \put{$\bullet$} at 4 3
\put{$\bullet$} at 5 0 \put{$\bullet$} at 5 3
\put{$\bullet$} at 6 0 \put{$\bullet$} at 6 3
\put{$\bullet$} at 7 0 \put{$\bullet$} at 7 3
\plot 0 3 3 0 /
\plot 1 3 4 0 /
\plot 2 3 0 0 /
\plot 3 3 1 0 /
\plot 4 3 5 0 /
\plot 5 3 6 0 /
\plot 6 3 7 0 /
\plot 7 3 2 0 /
\endpicture}
\cr
{\beginpicture
\setcoordinatesystem units <0.4cm,0.18cm>         
\setplotarea x from 0 to 8.5, y from 0 to 3    
\linethickness=0.5pt                          
\put{$\bullet$} at 0 0 \put{$\bullet$} at 0 3
\put{$\bullet$} at 1 0 \put{$\bullet$} at 1 3
\put{$\bullet$} at 2 0 \put{$\bullet$} at 2 3
\put{$\bullet$} at 3 0 \put{$\bullet$} at 3 3
\put{$\bullet$} at 4 0 \put{$\bullet$} at 4 3
\put{$\bullet$} at 5 0 \put{$\bullet$} at 5 3
\put{$\bullet$} at 6 0 \put{$\bullet$} at 6 3
\put{$\bullet$} at 7 0 \put{$\bullet$} at 7 3
\put{$=$} at 8 1.5
\plot 1 0 2 0 /
\plot 0 3 1 3 /
\plot 4 3 5 3 /
\plot 1 3 1 0  /
\plot 3 0 4 3  /
\plot 5 3 6 3  /
\plot 6 0 5 0  /
\plot 5 3 5 0  /
\setquadratic
\plot 3 3 5 2 7 3 /
\plot 4 0 5.5 1 7 0  /
\plot 4 0 3 1 2 0  /
\endpicture}
&{\beginpicture
\setcoordinatesystem units <0.4cm,0.18cm>         
\setplotarea x from 0 to 7, y from 0 to 3   
\linethickness=0.5pt                        
\put{$\bullet$} at 0 0 \put{$\bullet$} at 0 3
\put{$\bullet$} at 1 0 \put{$\bullet$} at 1 3
\put{$\bullet$} at 2 0 \put{$\bullet$} at 2 3
\put{$\bullet$} at 3 0 \put{$\bullet$} at 3 3
\put{$\bullet$} at 4 0 \put{$\bullet$} at 4 3
\put{$\bullet$} at 5 0 \put{$\bullet$} at 5 3
\put{$\bullet$} at 6 0 \put{$\bullet$} at 6 3
\put{$\bullet$} at 7 0 \put{$\bullet$} at 7 3
\plot 1 3 2 3 /
\plot 3 3 4 3 /
\plot 1 0 3 3 /
\plot 1 0 4 0 /
\plot 4 3 4 0 /
\plot 5 3 7 3 /
\plot 5 0 5 3 /
\plot 7 0 7 3 /
\plot 5 0 7 0 /
\endpicture}\cr
&{\beginpicture
\setcoordinatesystem units <0.4cm,0.18cm>         
\setplotarea x from 0 to 7, y from 0 to 3    
\linethickness=0.5pt                         
\put{$\bullet$} at 0 0 \put{$\bullet$} at 0 3
\put{$\bullet$} at 1 0 \put{$\bullet$} at 1 3
\put{$\bullet$} at 2 0 \put{$\bullet$} at 2 3
\put{$\bullet$} at 3 0 \put{$\bullet$} at 3 3
\put{$\bullet$} at 4 0 \put{$\bullet$} at 4 3
\put{$\bullet$} at 5 0 \put{$\bullet$} at 5 3
\put{$\bullet$} at 6 0 \put{$\bullet$} at 6 3
\put{$\bullet$} at 7 0 \put{$\bullet$} at 7 3
\plot 0 3 0 0 /
\plot 1 3 1 0 /
\plot 2 3 2 0 /
\plot 3 3 4 0 /
\plot 4 3 7 0 /
\plot 5 3 3 0 /
\plot 6 3 5 0 /
\plot 7 3 6 0 /
\endpicture}  \cr
}.\formula$$
For $\ell\in \ZZ_{>0}$, define
$$\matrix{
\hbox{the {\it Bell number,}} \hfill
&B(\ell) =
\hbox{(the number of set partitions of $\{1, 2, \ldots, \ell\}$)}, \cr
\cr
\hbox{the {\it Catalan number,}}\hfill
&\displaystyle{C(\ell) = {1 \over \ell + 1} {2 \ell \choose \ell}
= {2 \ell \choose \ell} - {2 \ell \choose \ell+ 1},} \hfill \cr
\cr
(2 \ell)!! = (2\ell-1)\cdot (2 \ell -3)&\hskip-.1in \cdots 5 \cdot 3 \cdot
1,
\qquad\hbox{and} \qquad
\ell ! = \ell \cdot (\ell-1) \cdots  2 \cdot 1, \hfill\cr
}
\formula
$$
with generating functions (see [Sta, 1.24f, and 6.2]),
$$
\matrix{
\displaystyle{
\sum_{\ell \ge 0} B(\ell) {z^\ell \over \ell !}  = \exp(e^z - 1),}\hfill
&\qquad
&\displaystyle{
\sum_{\ell \ge 0} C(\ell-1) z^\ell
= {1 - \sqrt{1 - 4z}\over 2z},} \hfill \cr
\cr
\displaystyle{\sum_{\ell \ge 0} (2 (\ell-1))!! {z^\ell \over \ell !}
= {1 - \sqrt{1 - 2 z}\over z},} \hfill
&&\displaystyle{\sum_{\ell \ge 0} \ell ! {z^\ell \over \ell !}
= {1 \over 1 - z}.} \hfill \cr
}
\formula
$$
Then
$$\matrix{
\hbox{for $k \in {1 \over 2} \ZZ_{>0}$,}
\qquad
&\Card(A_k) = B(2 k)\hfill &\quad\hbox{and}\quad
&\Card(P_k) = \Card(T_{2k}) = C(2k)\,,\hfill  \cr
\cr
\hbox{for $k \in \ZZ_{>0}$,} \hfill
&\Card(B_k) = (2k)!!, \hfill
&\quad\hbox{and}\quad
&\Card(S_k) = k!.\hfill \cr}
\formula
$$

\eject
\subsection{Presentation of the Partition Monoid}

In this section, for convenience, we will write
$$d_1 d_2 = d_1 \circ d_2, \qquad\hbox{for $d_1, d_2 \in A_k$}
$$
Let $k \in \ZZ_{>0}$. For $1 \le i \le k-1$ and $1 \le j \le k$, define
$$
\matrix{
p_{i+{1\over2}} =
{\beginpicture
\setcoordinatesystem units <0.5cm,0.2cm> 
\setplotarea x from 1 to 5.5, y from -1 to 4   
\linethickness=0.5pt                        
\put{$\bullet$} at 1 -1 \put{$\bullet$} at 1 2 \plot 1 -1 1 2 /
\put{$\cdots$} at 2 .5
\put{$\bullet$} at 3 -1  \put{$\bullet$} at 3 2 \plot 3 -1 3 2 /
\put{$\bullet$} at 4 -1 \put{$\bullet$} at 4 2
\put{$\bullet$} at 5 -1 \put{$\bullet$} at 5 2
\put{$\scriptstyle{i}$}[b] at 4 3.5
\put{$\scriptstyle{i+1}$}[b] at 5 3.45
\put{$\bullet$} at 6 -1 \put{$\bullet$} at 6 2  \plot 6 -1 6 2 /
\put{$\cdots$} at 7 .5
\put{$\bullet$} at 8 -1 \put{$\bullet$} at 8 2  \plot 8 -1 8 2 /
\plot 4 2 5 2 /
\plot 4 -1 5 -1 /
\plot 4 -1 4 2 /
\plot 5 -1 5 2 /
\endpicture},
&\qquad\qquad&
p_j =
{\beginpicture
\setcoordinatesystem units <0.5cm,0.2cm> 
\setplotarea x from 2 to 5.5, y from -1 to 4   
\linethickness=0.5pt                        
\put{$\bullet$} at 2 -1 \put{$\bullet$} at 2 2 \plot 2 -1 2 2 /
\put{$\cdots$} at 3 .5
\put{$\bullet$} at 4 -1  \put{$\bullet$} at 4 2 \plot 4 -1 4 2 /
\put{$\bullet$} at 5 -1 \put{$\bullet$} at 5 2
\put{$\scriptstyle{j}$}[b] at 5 3.5
\put{$\bullet$} at 6 -1 \put{$\bullet$} at 6 2  \plot 6 -1 6 2 /
\put{$\cdots$} at 7 .5
\put{$\bullet$} at 8 -1 \put{$\bullet$} at 8 2  \plot 8 -1 8 2 /
\endpicture},\hfill
\cr\cr\cr
\hfill e_{i} =
{\beginpicture
\setcoordinatesystem units <0.5cm,0.2cm> 
\setplotarea x from 2 to 5.5, y from -1 to 4   
\linethickness=0.5pt                        
\put{$\bullet$} at 2 -1 \put{$\bullet$} at 2 2 \plot 2 -1 2 2 /
\put{$\cdots$} at 3 .5
\put{$\bullet$} at 4 -1  \put{$\bullet$} at 4 2 \plot 4 -1 4 2 /
\put{$\bullet$} at 5 -1 \put{$\bullet$} at 5 2
\put{$\bullet$} at 6 -1 \put{$\bullet$} at 6 2
\put{$\scriptstyle{i}$}[b] at 5 3.45
\put{$\scriptstyle{i+1}$}[b] at 6 3.5
\put{$\bullet$} at 7 -1  \put{$\bullet$} at 7 2 \plot 7 -1 7 2 /
\put{$\cdots$} at 8 .5
\put{$\bullet$} at 9 -1  \put{$\bullet$} at 9 2 \plot 9 -1 9 2 /
\setquadratic
\plot 5 2 5.5 1.5 6 2 /
\plot 5 -1 5.5 -.5 6 -1 /
\endpicture} ,
& &
s_{i} =
{\beginpicture
\setcoordinatesystem units <0.5cm,0.2cm> 
\setplotarea x from 1 to 5.5, y from -1 to 4   
\linethickness=0.5pt                        
\put{$\bullet$} at 1 -1 \put{$\bullet$} at 1 2 \plot 1 -1 1 2 /
\put{$\cdots$} at 2 .5
\put{$\bullet$} at 3 -1  \put{$\bullet$} at 3 2 \plot 3 -1 3 2 /
\put{$\bullet$} at 4 -1 \put{$\bullet$} at 4 2
\put{$\bullet$} at 5 -1 \put{$\bullet$} at 5 2
\put{$\scriptstyle{i}$}[b] at 4 3.5
\put{$\scriptstyle{i+1}$}[b] at 5 3.45
\put{$\bullet$} at 6 -1 \put{$\bullet$} at 6 2  \plot 6 -1 6 2 /
\put{$\cdots$} at 7 .5
\put{$\bullet$} at 8 -1 \put{$\bullet$} at 8 2  \plot 8 -1 8 2 /
\plot 4 2 5 -1 /
\plot 5 2 4 -1 /
\endpicture}.
} \formula
$$
Note that
$e_{i} = p_{i+{1\over2}} p_i p_{i+1}  p_{i+{1\over2}}.$

\thm
\item{(a)} The monoid $T_k$ is presented by generators $e_1, \ldots,
e_{k-1}$
and relations
$$
e_i^2 = e_i,\qquad
e_i e_{i\pm1} e_i = e_i,
\qquad\hbox{and}\qquad
 e_ie_j=e_je_i,\quad\hbox{for $|i-j|>1$.}
$$
\item{(b)} The monoid $P_k$ is presented by generators $p_{1 \over
2},p_1,p_{3 \over 2},
\ldots, p_{k}$ and relations
$$
p_i^2 = p_i,\qquad
p_i p_{i\pm {1\over 2}} p_i = p_i,
\qquad\hbox{and}\qquad
p_i p_j=p_jp_i,\quad\hbox{for $|i-j|>1/2$.}
$$
\item{(c)} The group $S_k$ is presented by generators $s_1, \ldots, s_{k-1}$
and relations
$$
s_i^2 = 1,\qquad
s_i s_{i+1} s_i = s_{i+1} s_i  s_{i+1},
\qquad\hbox{and}\qquad
 s_i s_j=s_j s_i,\quad\hbox{for $|i-j|>1$.}
$$
\item{(d)} The monoid $A_k$ is presented by generators $s_1, \ldots,
s_{k-1}$ and
$p_{1 \over 2},p_1,p_{3 \over 2}, \ldots, p_{k}$
and relations in (b) and (c) and
$$
\matrix{
s_i p_i p_{i+1} = p_ip_{i+1}s_i = p_ip_{i+1},\qquad
s_i p_{i+{1\over2}}=p_{i+{1\over2}}s_i=p_{i+{1\over2}}, \qquad
s_i p_{i} s_i = p_{i+1},\qquad\qquad\cr
s_is_{i+1}p_{i+{1\over2}}s_{i+1}s_i = p_{i+{3\over2}},
\qquad\hbox{and}\qquad
s_i p_j = p_j s_i, \quad
\hbox{ for $j\ne i-{1\over2}, i, i+{1\over2}, i+1, i+{3\over2}$}.
\cr}
$$

\pf  Parts (a) and (c) are standard.  See [GHJ, Prop.\ 2.8.1] and
[Bou2, Ch.\ IV \S 1.3, Ex.\ 2], respectively.
Part (b) is a consequence of (a) and the monoid isomorphism in (1.5).

\smallskip\noindent
(d)
The right way to think of this is to realize that $A_k$ is defined as
a presentation by the generators $d\in A_k$ and the relations  which
specify the composition of diagrams.  To prove the presentation in
the statement of the theorem we need to establish that the generators
and relations in each of these two presentations can be derived from
each other.  Thus it is  sufficient to show that
\smallskip\noindent
\itemitem{(1)} The generators in (1.10) satisfy the relations in
Theorem (1.11).
\smallskip\noindent
\itemitem{(2)} Every set partition $d\in A_k$ can be written as a
product of the generators in (1.10).
\smallskip\noindent
\itemitem{(3)} Any product $d_1\circ d_2$ can be computed using the
relations in Theorem (1.11).

\smallskip\noindent
(1) is established by a direct check using the definition of the
multiplication of diagrams.
(2) follows from (b) and (c) and the fact (1.6) that $A_k = S_k P_k S_k$.
The bulk of the work is in proving (3).
\smallskip\noindent
{\it Step 1.}\
First note that the relations in (a--d) imply the following relations:
\smallskip\noindent
\itemitem{$\matrix{\hbox{\sl (e1)}\cr\cr}$}
$\eqalign{
p_{i+{1\over2}}s_{i-1}p_{i+{1\over2}}
&=p_{i+{1\over2}}s_is_{i-1}p_{i+{1\over2}}
=p_{i+{1\over2}}s_is_{i-1}p_{i+{1\over2}}s_{i-1}s_is_is_{i-1} \cr
&=p_{i-{1\over2}}p_{i+{1\over2}}s_is_{i-1}
=p_{i-{1\over2}}p_{i+{1\over2}}s_{i-1}
=p_{i+{1\over2}}p_{i-{1\over2}}s_{i-1}
=p_{i+{1\over2}}p_{i-{1\over2}}. \cr
}$
\smallskip
\itemitem{{\sl (e2)}}
$p_is_ip_i = s_is_ip_is_ip_i = s_ip_{i+1}p_i = p_{i+1}p_i$.
\smallskip
\itemitem{{\sl (f1)}}
$p_ip_{i+{1\over2}}p_{i+1}
=p_ip_{i+{1\over2}}s_ip_{i+1}
=p_ip_{i+{1\over2}}p_is_i
=p_is_i$.
\smallskip
\itemitem{{\sl (f2)}}
$p_{i+1}p_{i+{1\over2}}p_i = p_{i+1}s_ip_{i+{1\over2}}p_i
=s_ip_ip_{i+{1\over2}}p_i
=s_ip_i$.

\medskip\noindent
{\it Step 2.}\ Analyze how elements of $P_k$ can be efficiently
expressed in terms of the generators.
\smallskip
Let $t\in P_k$. The blocks of $t$ partition $\{1, \ldots, k\}$ into
{\it top blocks\/} and partition $\{1', \ldots, k'\}$ into
{\it bottom blocks\/}.
In $t$, some top blocks are connected to bottom blocks by an
edge, but no top block is connected to two bottom blocks,
for then by transitivity the two bottom blocks are actually a
single block.  Draw the diagram of $t$, such that if a top block
connects to a bottom block, then it connects with a single edge
joining the leftmost vertices in each block.
The element $t \in P_k$ can be decomposed in {\it block form\/} as
$$t =
(p_{i_1+{1\over2}} \cdots p_{i_r+{1\over2}})
(p_{j_1} \cdots p_{j_s})
\tau
(p_{\ell_1} \cdots p_{\ell_m})
(p_{r_1+{1\over2}} \cdots p_{r_n+{1\over2}}) \formula
$$
with $\tau\in S_k$, $i_1 < i_2 < \cdots < i_r$, $j_1 < j_2 < \cdots < j_s$,
$\ell_1 < \ell_2 < \cdots < \ell_m$, and $r_1 < r_2 < \cdots < r_n$.
The left product of $p_i$s corresponds to the top blocks of $t$,
the right product of $p_i$s corresponds to the bottom blocks of $t$
and the permutation $\tau$ corresponds to the propagation pattern
of the edges connecting top blocks of $t$ to bottom blocks of $t$.
For example,
$$ \eqalign{
t  =
{\beginpicture
\setcoordinatesystem units <0.3cm,0.18cm>         
\setplotarea x from 1 to 8, y from -1 to 2    
\linethickness=0.5pt                         
\put{$\bullet$} at 1 -1 \put{$\bullet$} at 1 2
\put{$\bullet$} at 2 -1 \put{$\bullet$} at 2 2
\put{$\bullet$} at 3 -1 \put{$\bullet$} at 3 2
\put{$\bullet$} at 4 -1 \put{$\bullet$} at 4 2
\put{$\bullet$} at 5 -1 \put{$\bullet$} at 5 2
\put{$\bullet$} at 6 -1 \put{$\bullet$} at 6 2
\put{$\bullet$} at 7 -1 \put{$\bullet$} at 7 2
\put{$\bullet$} at 8 -1 \put{$\bullet$} at 8 2
\plot 1 -1 1 2 /
\plot 2 2 5 -1 /
\plot 5 2 6 -1 /
\plot 8 2 8 -1 /
\plot 2 2 4 2 /
\plot 1 -1 3 -1 /
\plot 6 2 7 2 /
\endpicture}
& =
\matrix{
{\beginpicture
\setcoordinatesystem units <0.3cm,0.18cm>         
\setplotarea x from 1 to 8, y from -1 to 2    
\linethickness=0.5pt                         
\put{$\bullet$} at 1 -1 \put{$\bullet$} at 1 2
\put{$\bullet$} at 2 -1 \put{$\bullet$} at 2 2
\put{$\bullet$} at 3 -1 \put{$\bullet$} at 3 2
\put{$\bullet$} at 4 -1 \put{$\bullet$} at 4 2
\put{$\bullet$} at 5 -1 \put{$\bullet$} at 5 2
\put{$\bullet$} at 6 -1 \put{$\bullet$} at 6 2
\put{$\bullet$} at 7 -1 \put{$\bullet$} at 7 2
\put{$\bullet$} at 8 -1 \put{$\bullet$} at 8 2
\plot 1 -1 1 2 /
\plot 2 -1 2 2 4 2 /
\plot 5 -1 5 2 /
\plot 6 2 7 2  /
\plot 8 -1 8 2 /
\endpicture} \cr
{\beginpicture
\setcoordinatesystem units <0.3cm,0.18cm>         
\setplotarea x from 1 to 8, y from -1 to 2    
\linethickness=0.5pt                         
\put{$\bullet$} at 1 -1 \put{$\bullet$} at 1 2
\put{$\bullet$} at 2 -1 \put{$\bullet$} at 2 2
\put{$\bullet$} at 3 -1 \put{$\bullet$} at 3 2
\put{$\bullet$} at 4 -1 \put{$\bullet$} at 4 2
\put{$\bullet$} at 5 -1 \put{$\bullet$} at 5 2
\put{$\bullet$} at 6 -1 \put{$\bullet$} at 6 2
\put{$\bullet$} at 7 -1 \put{$\bullet$} at 7 2
\put{$\bullet$} at 8 -1 \put{$\bullet$} at 8 2
\plot 1 -1 1 2 /
\plot 2 2 5 -1 /
\plot 5 2 6 -1 /
\plot 8 2 8 -1 /
\endpicture} \cr
{\beginpicture
\setcoordinatesystem units <0.3cm,0.18cm>         
\setplotarea x from 1 to 8, y from -1 to 2    
\linethickness=0.5pt                         
\put{$\bullet$} at 1 -1 \put{$\bullet$} at 1 2
\put{$\bullet$} at 2 -1 \put{$\bullet$} at 2 2
\put{$\bullet$} at 3 -1 \put{$\bullet$} at 3 2
\put{$\bullet$} at 4 -1 \put{$\bullet$} at 4 2
\put{$\bullet$} at 5 -1 \put{$\bullet$} at 5 2
\put{$\bullet$} at 6 -1 \put{$\bullet$} at 6 2
\put{$\bullet$} at 7 -1 \put{$\bullet$} at 7 2
\put{$\bullet$} at 8 -1 \put{$\bullet$} at 8 2
\plot 1 2 1 -1 3 -1 /
\plot 5 2 5 -1 /
\plot 6 2 6 -1 /
\plot 8 2 8 -1 /
\endpicture}
} =
\matrix{
{\beginpicture
\setcoordinatesystem units <0.3cm,0.18cm>         
\setplotarea x from 1 to 8, y from -1 to 2    
\linethickness=0.5pt                         
\put{$\bullet$} at 1 -1 \put{$\bullet$} at 1 2
\put{$\bullet$} at 2 -1 \put{$\bullet$} at 2 2
\put{$\bullet$} at 3 -1 \put{$\bullet$} at 3 2
\put{$\bullet$} at 4 -1 \put{$\bullet$} at 4 2
\put{$\bullet$} at 5 -1 \put{$\bullet$} at 5 2
\put{$\bullet$} at 6 -1 \put{$\bullet$} at 6 2
\put{$\bullet$} at 7 -1 \put{$\bullet$} at 7 2
\put{$\bullet$} at 8 -1 \put{$\bullet$} at 8 2
\plot 1 -1 1 2 /
\plot 2 -1 2 2 4 2 /
\plot 5 -1 5 2 /
\plot 6 2 7 2  /
\plot 8 -1 8 2 /
\endpicture} \cr
{\beginpicture
\setcoordinatesystem units <0.3cm,0.18cm>         
\setplotarea x from 1 to 8, y from -1 to 2    
\linethickness=0.5pt                         
\put{$\bullet$} at 1 -1 \put{$\bullet$} at 1 2
\put{$\bullet$} at 2 -1 \put{$\bullet$} at 2 2
\put{$\bullet$} at 3 -1 \put{$\bullet$} at 3 2
\put{$\bullet$} at 4 -1 \put{$\bullet$} at 4 2
\put{$\bullet$} at 5 -1 \put{$\bullet$} at 5 2
\put{$\bullet$} at 6 -1 \put{$\bullet$} at 6 2
\put{$\bullet$} at 7 -1 \put{$\bullet$} at 7 2
\put{$\bullet$} at 8 -1 \put{$\bullet$} at 8 2
\plot 1 -1 1 2 /
\plot 2 2 5 -1 /
\plot 5 2 6 -1 /
\plot 8 2 8 -1 /
\setdashes <.04cm>
\plot 3 2 2 -1 /
\plot 4 2 3 -1 /
\plot 6 2 4 -1 /
\plot 7 2 7 -1 /
\endpicture} \cr
{\beginpicture
\setcoordinatesystem units <0.3cm,0.18cm>         
\setplotarea x from 1 to 8, y from -1 to 2    
\linethickness=0.5pt                         
\put{$\bullet$} at 1 -1 \put{$\bullet$} at 1 2
\put{$\bullet$} at 2 -1 \put{$\bullet$} at 2 2
\put{$\bullet$} at 3 -1 \put{$\bullet$} at 3 2
\put{$\bullet$} at 4 -1 \put{$\bullet$} at 4 2
\put{$\bullet$} at 5 -1 \put{$\bullet$} at 5 2
\put{$\bullet$} at 6 -1 \put{$\bullet$} at 6 2
\put{$\bullet$} at 7 -1 \put{$\bullet$} at 7 2
\put{$\bullet$} at 8 -1 \put{$\bullet$} at 8 2
\plot 1 2 1 -1 3 -1 /
\plot 5 2 5 -1 /
\plot 6 2 6 -1 /
\plot 8 2 8 -1 /
\endpicture}
} =
\matrix{
{\beginpicture
\setcoordinatesystem units <0.3cm,0.18cm>         
\setplotarea x from 1 to 8, y from -1 to 2    
\linethickness=0.5pt                         
\put{$\bullet$} at 1 -1 \put{$\bullet$} at 1 2
\put{$\bullet$} at 2 -1 \put{$\bullet$} at 2 2
\put{$\bullet$} at 3 -1 \put{$\bullet$} at 3 2
\put{$\bullet$} at 4 -1 \put{$\bullet$} at 4 2
\put{$\bullet$} at 5 -1 \put{$\bullet$} at 5 2
\put{$\bullet$} at 6 -1 \put{$\bullet$} at 6 2
\put{$\bullet$} at 7 -1 \put{$\bullet$} at 7 2
\put{$\bullet$} at 8 -1 \put{$\bullet$} at 8 2
\plot 1 -1 1 2 /
\plot 2 -1 2 2 4 2 4 -1 2 -1 /
\plot 5 -1 5 2 /
\plot 6 2 7 2  7 -1 6 -1 6 2 /
\plot 8 -1 8 2 /
\endpicture} \cr
{\beginpicture
\setcoordinatesystem units <0.3cm,0.18cm>         
\setplotarea x from 1 to 8, y from -1 to 2    
\linethickness=0.5pt                         
\put{$\bullet$} at 1 -1 \put{$\bullet$} at 1 2
\put{$\bullet$} at 2 -1 \put{$\bullet$} at 2 2
\put{$\bullet$} at 3 -1 \put{$\bullet$} at 3 2
\put{$\bullet$} at 4 -1 \put{$\bullet$} at 4 2
\put{$\bullet$} at 5 -1 \put{$\bullet$} at 5 2
\put{$\bullet$} at 6 -1 \put{$\bullet$} at 6 2
\put{$\bullet$} at 7 -1 \put{$\bullet$} at 7 2
\put{$\bullet$} at 8 -1 \put{$\bullet$} at 8 2
\plot 1 -1 1 2 /
\plot 2 -1 2 2 /
\plot 5 -1 5 2 /
\plot 8 -1 8 2 /
\endpicture} \cr
{\beginpicture
\setcoordinatesystem units <0.3cm,0.18cm>         
\setplotarea x from 1 to 8, y from -1 to 2    
\linethickness=0.5pt                         
\put{$\bullet$} at 1 -1 \put{$\bullet$} at 1 2
\put{$\bullet$} at 2 -1 \put{$\bullet$} at 2 2
\put{$\bullet$} at 3 -1 \put{$\bullet$} at 3 2
\put{$\bullet$} at 4 -1 \put{$\bullet$} at 4 2
\put{$\bullet$} at 5 -1 \put{$\bullet$} at 5 2
\put{$\bullet$} at 6 -1 \put{$\bullet$} at 6 2
\put{$\bullet$} at 7 -1 \put{$\bullet$} at 7 2
\put{$\bullet$} at 8 -1 \put{$\bullet$} at 8 2
\plot 1 -1 1 2 /
\plot 2 2 5 -1 /
\plot 5 2 6 -1 /
\plot 8 2 8 -1 /
\setdashes <.04cm>
\plot 3 2 2 -1 /
\plot 4 2 3 -1 /
\plot 6 2 4 -1 /
\plot 7 2 7 -1 /
\endpicture} \cr
{\beginpicture
\setcoordinatesystem units <0.3cm,0.18cm>         
\setplotarea x from 1 to 8, y from -1 to 2    
\linethickness=0.5pt                         
\put{$\bullet$} at 1 -1 \put{$\bullet$} at 1 2
\put{$\bullet$} at 2 -1 \put{$\bullet$} at 2 2
\put{$\bullet$} at 3 -1 \put{$\bullet$} at 3 2
\put{$\bullet$} at 4 -1 \put{$\bullet$} at 4 2
\put{$\bullet$} at 5 -1 \put{$\bullet$} at 5 2
\put{$\bullet$} at 6 -1 \put{$\bullet$} at 6 2
\put{$\bullet$} at 7 -1 \put{$\bullet$} at 7 2
\put{$\bullet$} at 8 -1 \put{$\bullet$} at 8 2
\plot 1 2 1 -1  /
\plot 5 2 5 -1 /
\plot 6 2 6 -1 /
\plot 8 2 8 -1 /
\endpicture}\cr
{\beginpicture
\setcoordinatesystem units <0.3cm,0.18cm>         
\setplotarea x from 1 to 8, y from -1 to 2    
\linethickness=0.5pt                         
\put{$\bullet$} at 1 -1 \put{$\bullet$} at 1 2
\put{$\bullet$} at 2 -1 \put{$\bullet$} at 2 2
\put{$\bullet$} at 3 -1 \put{$\bullet$} at 3 2
\put{$\bullet$} at 4 -1 \put{$\bullet$} at 4 2
\put{$\bullet$} at 5 -1 \put{$\bullet$} at 5 2
\put{$\bullet$} at 6 -1 \put{$\bullet$} at 6 2
\put{$\bullet$} at 7 -1 \put{$\bullet$} at 7 2
\put{$\bullet$} at 8 -1 \put{$\bullet$} at 8 2
\plot 1 2 1 -1 3 -1 3 2 1 2 /
\plot 4 2 4 -1 /
\plot 5 2 5 -1 /
\plot 6 2 6 -1 /
\plot 7 2 7 -1 /
\plot 8 2 8 -1 /
\endpicture}
} \cr
& = (p_{2{1\over2}}p_{3{1\over2}}p_{6{1\over2}})
(p_3p_4p_6p_7)\tau
(p_2p_3p_4p_7)(p_{1{1\over2}}p_{2{1\over2}}),\cr
& = (p_{2{1\over2}}p_{3{1\over2}}p_{6{1\over2}})
(p_3p_4p_6p_7) s_2s_3s_5s_4
(p_2p_3p_4p_7)(p_{1{1\over2}}p_{2{1\over2}}),
}$$
The dashed edges of $\tau$ are ``non-propagating" edges, and
they may be chosen so that they do not cross each other.
The propagating edges of $\tau$ do not cross, since $t$ is planar.

Using the relations {\sl (f1)} and {\sl (f2)}, the
non-propagating edges of $\tau$ can be ``removed'',
leaving a planar diagram which is
written in terms of the generators $p_i$ and $p_{i+{1\over2}}$. In
our example, this process will replace $\tau$ by
$p_{2{1\over2}}p_2p_{3{1\over2}}p_3 p_{5{1\over2}}p_5 p_{4{1\over2}}p_4,$
so that
$$
t =
{\beginpicture
\setcoordinatesystem units <0.3cm,0.18cm>         
\setplotarea x from 1 to 8, y from -1 to 2    
\linethickness=0.5pt                         
\put{$\bullet$} at 1 -1 \put{$\bullet$} at 1 2
\put{$\bullet$} at 2 -1 \put{$\bullet$} at 2 2
\put{$\bullet$} at 3 -1 \put{$\bullet$} at 3 2
\put{$\bullet$} at 4 -1 \put{$\bullet$} at 4 2
\put{$\bullet$} at 5 -1 \put{$\bullet$} at 5 2
\put{$\bullet$} at 6 -1 \put{$\bullet$} at 6 2
\put{$\bullet$} at 7 -1 \put{$\bullet$} at 7 2
\put{$\bullet$} at 8 -1 \put{$\bullet$} at 8 2
\plot 1 -1 1 2 /
\plot 2 2 5 -1 /
\plot 5 2 6 -1 /
\plot 8 2 8 -1 /
\plot 2 2 4 2 /
\plot 1 -1 3 -1 /
\plot 6 2 7 2 /
\endpicture}
 =
\matrix{
{\beginpicture
\setcoordinatesystem units <0.3cm,0.18cm>         
\setplotarea x from 1 to 8, y from -1 to 2    
\linethickness=0.5pt                         
\put{$\bullet$} at 1 -1 \put{$\bullet$} at 1 2
\put{$\bullet$} at 2 -1 \put{$\bullet$} at 2 2
\put{$\bullet$} at 3 -1 \put{$\bullet$} at 3 2
\put{$\bullet$} at 4 -1 \put{$\bullet$} at 4 2
\put{$\bullet$} at 5 -1 \put{$\bullet$} at 5 2
\put{$\bullet$} at 6 -1 \put{$\bullet$} at 6 2
\put{$\bullet$} at 7 -1 \put{$\bullet$} at 7 2
\put{$\bullet$} at 8 -1 \put{$\bullet$} at 8 2
\plot 1 -1 1 2 /
\plot 2 -1 2 2 4 2 /
\plot 5 -1 5 2 /
\plot 6 2 7 2  /
\plot 8 -1 8 2 /
\endpicture} \cr
{\beginpicture
\setcoordinatesystem units <0.3cm,0.15cm>         
\setplotarea x from 1 to 8, y from -1 to 2    
\linethickness=0.5pt                         
\put{$\bullet$} at 1 -1 \put{$\bullet$} at 1 2
\put{$\bullet$} at 2 -1 \put{$\bullet$} at 2 2
\put{$\bullet$} at 3 -1 \put{$\bullet$} at 3 2
\put{$\bullet$} at 4 -1 \put{$\bullet$} at 4 2
\put{$\bullet$} at 5 -1 \put{$\bullet$} at 5 2
\put{$\bullet$} at 6 -1 \put{$\bullet$} at 6 2
\put{$\bullet$} at 7 -1 \put{$\bullet$} at 7 2
\put{$\bullet$} at 8 -1 \put{$\bullet$} at 8 2
\plot 1 -1 1 2 /
\plot 2 2 4 2 4 -1 /
\plot 5 2 6 2 6 -1 /
\plot 7 2 7 -1 /
\plot 8 2 8 -1 /
\endpicture} \cr
{\beginpicture
\setcoordinatesystem units <0.3cm,0.15cm>         
\setplotarea x from 1 to 8, y from -1 to 2    
\linethickness=0.5pt                         
\put{$\bullet$} at 1 -1 \put{$\bullet$} at 1 2
\put{$\bullet$} at 2 -1 \put{$\bullet$} at 2 2
\put{$\bullet$} at 3 -1 \put{$\bullet$} at 3 2
\put{$\bullet$} at 4 -1 \put{$\bullet$} at 4 2
\put{$\bullet$} at 5 -1 \put{$\bullet$} at 5 2
\put{$\bullet$} at 6 -1 \put{$\bullet$} at 6 2
\put{$\bullet$} at 7 -1 \put{$\bullet$} at 7 2
\put{$\bullet$} at 8 -1 \put{$\bullet$} at 8 2
\plot 1 -1 1 2 /
\plot 2 -1 2 2 /
\plot 3 -1 3 2 /
\plot 4 2 5 2 5 -1 /
\plot 6 -1 6 2 /
\plot 7 -1 7 2 /
\plot 8 2 8 -1 /
\endpicture} \cr
{\beginpicture
\setcoordinatesystem units <0.3cm,0.18cm>         
\setplotarea x from 1 to 8, y from -1 to 2    
\linethickness=0.5pt                         
\put{$\bullet$} at 1 -1 \put{$\bullet$} at 1 2
\put{$\bullet$} at 2 -1 \put{$\bullet$} at 2 2
\put{$\bullet$} at 3 -1 \put{$\bullet$} at 3 2
\put{$\bullet$} at 4 -1 \put{$\bullet$} at 4 2
\put{$\bullet$} at 5 -1 \put{$\bullet$} at 5 2
\put{$\bullet$} at 6 -1 \put{$\bullet$} at 6 2
\put{$\bullet$} at 7 -1 \put{$\bullet$} at 7 2
\put{$\bullet$} at 8 -1 \put{$\bullet$} at 8 2
\plot 1 2 1 -1 3 -1 /
\plot 5 2 5 -1 /
\plot 6 2 6 -1 /
\plot 8 2 8 -1 /
\endpicture}
} = \matrix{
(p_{2{1\over2}}p_{3{1\over2}}p_{6{1\over2}}) (p_3p_4p_6p_7) \cr
\cr
\qquad\cdot\ p_{2{1\over2}}p_2p_{3{1\over2}}p_3 p_{5{1\over2}}p_5
p_{4{1\over2}}p_4 \cr
\cr
\qquad\qquad\cdot\
(p_2p_3p_4p_7)(p_{1{1\over2}}p_{2{1\over2}}).  \cr}
$$

\smallskip\noindent
{\it Step 3}: If $t\in P_k$ and
$\sigma_1\in S_k$ which permutes the
the top blocks of the planar diagram $t$,
then there is a permutation $\sigma_2$ of the bottom
blocks of $t$ such that $\sigma_1t\sigma_2$ is planar.
Furthermore, this can be accomplished using the relations.
For example, suppose
$$
t =
{\beginpicture
\setcoordinatesystem units <0.3cm,0.18cm>         
\setplotarea x from 1 to 7, y from -1 to 2    
\linethickness=0.5pt                         
\put{$\bullet$} at 1 -1 \put{$\bullet$} at 1 2
\put{$\bullet$} at 2 -1 \put{$\bullet$} at 2 2
\put{$\bullet$} at 3 -1 \put{$\bullet$} at 3 2
\put{$\bullet$} at 4 -1 \put{$\bullet$} at 4 2
\put{$\bullet$} at 5 -1 \put{$\bullet$} at 5 2
\put{$\bullet$} at 6 -1 \put{$\bullet$} at 6 2
\put{$\bullet$} at 7 -1 \put{$\bullet$} at 7 2
\put{$\overbrace{\hphantom{XX}}^{T_1}$} at 2 4
\put{$\overbrace{\hphantom{i}}^{T_2}$} at 6.5 4
\put{$\underbrace{\hphantom{XXX}}_{B_1}$} at 2.5 -3
\put{$\underbrace{\hphantom{XX}}_{B_2}$} at 6 -3
\plot 1 2 3 2 /
\plot 1 2 1 -1 /
\plot 6 2 7 2 /
\plot 6 2 5 -1 /
\plot 1 -1 4 -1 /
\plot 5 -1 7 -1 /
\endpicture} =
\underbrace{(p_{1{1\over 2}} p_{2{1\over2}})(p_2p_3)}_{T_1}
\underbrace{(p_{6{1 \over 2}})(p_7)}_{T_2}
p_4p_5
s_5
\underbrace{(p_2 p_3 p_4)
(p_{1{1 \over 2}} p_{2{1 \over 2}}p_{3{1 \over 2}})}_{B_1}
\underbrace{(p_6 p_7)(p_{5{1 \over 2}} p_{6{1 \over 2}})}_{B_2}
$$
is a planar diagram with top blocks $T_1$ and $T_2$ connected
respectively to bottom blocks $B_1$ and $B_2$ and
$$\sigma_1=
{\beginpicture
\setcoordinatesystem units <0.3cm,0.2cm>         
\setplotarea x from 1 to 7, y from -1 to 2    
\linethickness=0.5pt                         
\put{$\bullet$} at 1 -1 \put{$\bullet$} at 1 2
\put{$\bullet$} at 2 -1 \put{$\bullet$} at 2 2
\put{$\bullet$} at 3 -1 \put{$\bullet$} at 3 2
\put{$\bullet$} at 4 -1 \put{$\bullet$} at 4 2
\put{$\bullet$} at 5 -1 \put{$\bullet$} at 5 2
\put{$\bullet$} at 6 -1 \put{$\bullet$} at 6 2
\put{$\bullet$} at 7 -1 \put{$\bullet$} at 7 2
\plot 1 2 6 -1 /
\plot 2 2 7 -1 /
\plot 1 -1 5 2 /
\plot 2 -1 6 2 /
\plot 3 -1 7 2 /
\plot 4 -1 3 2 /
\plot 5 -1 4 2 /
\endpicture}
\qquad\hbox{so that}\qquad
\sigma_1t = \matrix{
{\beginpicture
\setcoordinatesystem units <0.3cm,0.18cm>         
\setplotarea x from 1 to 7, y from -1 to 2    
\linethickness=0.5pt                         
\put{$\bullet$} at 1 -1 \put{$\bullet$} at 1 2
\put{$\bullet$} at 2 -1 \put{$\bullet$} at 2 2
\put{$\bullet$} at 3 -1 \put{$\bullet$} at 3 2
\put{$\bullet$} at 4 -1 \put{$\bullet$} at 4 2
\put{$\bullet$} at 5 -1 \put{$\bullet$} at 5 2
\put{$\bullet$} at 6 -1 \put{$\bullet$} at 6 2
\put{$\bullet$} at 7 -1 \put{$\bullet$} at 7 2
\plot 1 2 6 -1 /
\plot 2 2 7 -1 /
\plot 1 -1 5 2 /
\plot 2 -1 6 2 /
\plot 3 -1 7 2 /
\plot 4 -1 3 2 /
\plot 5 -1 4 2 /
\endpicture}\cr
{\beginpicture
\setcoordinatesystem units <0.3cm,0.18cm>         
\setplotarea x from 1 to 7, y from -1 to 2    
\linethickness=0.5pt                         
\put{$\bullet$} at 1 -1 \put{$\bullet$} at 1 2
\put{$\bullet$} at 2 -1 \put{$\bullet$} at 2 2
\put{$\bullet$} at 3 -1 \put{$\bullet$} at 3 2
\put{$\bullet$} at 4 -1 \put{$\bullet$} at 4 2
\put{$\bullet$} at 5 -1 \put{$\bullet$} at 5 2
\put{$\bullet$} at 6 -1 \put{$\bullet$} at 6 2
\put{$\bullet$} at 7 -1 \put{$\bullet$} at 7 2
\plot 1 2 3 2 /
\plot 1 2 1 -1 /
\plot 6 2 7 2 /
\plot 6 2 5 -1 /
\plot 1 -1 4 -1 /
\plot 5 -1 7 -1 /
\endpicture}}
=
{\beginpicture
\setcoordinatesystem units <0.3cm,0.18cm>         
\setplotarea x from 1 to 7, y from -1 to 2    
\linethickness=0.5pt                         
\put{$\bullet$} at 1 -1 \put{$\bullet$} at 1 2
\put{$\bullet$} at 2 -1 \put{$\bullet$} at 2 2
\put{$\bullet$} at 3 -1 \put{$\bullet$} at 3 2
\put{$\bullet$} at 4 -1 \put{$\bullet$} at 4 2
\put{$\bullet$} at 5 -1 \put{$\bullet$} at 5 2
\put{$\bullet$} at 6 -1 \put{$\bullet$} at 6 2
\put{$\bullet$} at 7 -1 \put{$\bullet$} at 7 2
\plot 1 2 2 2 /
\plot 1 2 5 -1 /
\plot 5 2 7 2 /
\plot 5 2 1 -1 /
\plot 1 -1 4 -1 /
\plot 5 -1 7 -1 /
\put{$\overbrace{\hphantom{i}}^{T_2'}$} at 1.5 4
\put{$\overbrace{\hphantom{i}}^{T_1'}$} at 6 4
\put{$\underbrace{\hphantom{XXX}}_{B_1}$} at 2.5 -3
\put{$\underbrace{\hphantom{XX}}_{B_2}$} at 6 -3
\endpicture}
= t',$$
then transposition of $B_1$ and $B_2$ can be accomplished
with the permutation
$$
\sigma_2=
{\beginpicture
\setcoordinatesystem units <0.3cm,0.18cm>         
\setplotarea x from 1 to 7, y from -1 to 2    
\linethickness=0.5pt                         
\put{$\bullet$} at 1 -1 \put{$\bullet$} at 1 2
\put{$\bullet$} at 2 -1 \put{$\bullet$} at 2 2
\put{$\bullet$} at 3 -1 \put{$\bullet$} at 3 2
\put{$\bullet$} at 4 -1 \put{$\bullet$} at 4 2
\put{$\bullet$} at 5 -1 \put{$\bullet$} at 5 2
\put{$\bullet$} at 6 -1 \put{$\bullet$} at 6 2
\put{$\bullet$} at 7 -1 \put{$\bullet$} at 7 2
\plot 5 2 1 -1  /
\plot 6 2 2 -1  /
\plot 7 2 3 -1  /
\plot 1 2 4 -1  /
\plot 2 2 5 -1  /
\plot 3 2 6 -1  /
\plot 4 2 7 -1  /
\endpicture}
\qquad\hbox{so that}\qquad
\sigma_1t\sigma_2=t'\sigma_2=
\matrix{
{\beginpicture
\setcoordinatesystem units <0.3cm,0.18cm>         
\setplotarea x from 1 to 7, y from -1 to 2    
\linethickness=0.5pt                         
\put{$\bullet$} at 1 -1 \put{$\bullet$} at 1 2
\put{$\bullet$} at 2 -1 \put{$\bullet$} at 2 2
\put{$\bullet$} at 3 -1 \put{$\bullet$} at 3 2
\put{$\bullet$} at 4 -1 \put{$\bullet$} at 4 2
\put{$\bullet$} at 5 -1 \put{$\bullet$} at 5 2
\put{$\bullet$} at 6 -1 \put{$\bullet$} at 6 2
\put{$\bullet$} at 7 -1 \put{$\bullet$} at 7 2
\plot 1 2 2 2 /
\plot 1 2 5 -1 /
\plot 5 2 7 2 /
\plot 5 2 1 -1 /
\plot 1 -1 4 -1 /
\plot 5 -1 7 -1 /
\endpicture}  \cr
{\beginpicture
\setcoordinatesystem units <0.3cm,0.18cm>         
\setplotarea x from 1 to 7, y from -1 to 2    
\linethickness=0.5pt                         
\put{$\bullet$} at 1 -1 \put{$\bullet$} at 1 2
\put{$\bullet$} at 2 -1 \put{$\bullet$} at 2 2
\put{$\bullet$} at 3 -1 \put{$\bullet$} at 3 2
\put{$\bullet$} at 4 -1 \put{$\bullet$} at 4 2
\put{$\bullet$} at 5 -1 \put{$\bullet$} at 5 2
\put{$\bullet$} at 6 -1 \put{$\bullet$} at 6 2
\put{$\bullet$} at 7 -1 \put{$\bullet$} at 7 2
\plot 5 2 1 -1  /
\plot 6 2 2 -1  /
\plot 7 2 3 -1  /
\plot 1 2 4 -1  /
\plot 2 2 5 -1  /
\plot 3 2 6 -1  /
\plot 4 2 7 -1  /
\endpicture}}
=
{\beginpicture
\setcoordinatesystem units <0.3cm,0.18cm>         
\setplotarea x from 1 to 7, y from -1 to 2    
\linethickness=0.5pt                         
\put{$\bullet$} at 1 -1 \put{$\bullet$} at 1 2
\put{$\bullet$} at 2 -1 \put{$\bullet$} at 2 2
\put{$\bullet$} at 3 -1 \put{$\bullet$} at 3 2
\put{$\bullet$} at 4 -1 \put{$\bullet$} at 4 2
\put{$\bullet$} at 5 -1 \put{$\bullet$} at 5 2
\put{$\bullet$} at 6 -1 \put{$\bullet$} at 6 2
\put{$\bullet$} at 7 -1 \put{$\bullet$} at 7 2
\plot 1 2 2 2 /
\plot 1 2 1 -1 /
\plot 5 2 7 2 /
\plot 5 2 4 -1 /
\plot 1 -1 3 -1 /
\plot 4 -1 7 -1 /
\put{$\overbrace{\hphantom{i}}^{T_2'}$} at 1.5 4
\put{$\overbrace{\hphantom{}}^{T_1'}$} at 6 4
\put{$\underbrace{\hphantom{XX}}_{B_2'}$} at 2 -3
\put{$\underbrace{\hphantom{XXX}}_{B_1'}$} at 5.5 -3
\endpicture}
$$
is planar.  It is possible to accomplish these products using the
relations from the statement of the theorem.
In our example, with
$\sigma_1 = s_2s_1s_3s_2s_4s_3s_5s_2s_4s_6s_1s_3s_5s_2s_4s_3$
and  with
$\sigma_2 = s_4 s_5 s_6 s_3 s_4 s_5  s_2 s_3s_4  s_1 s_2 s_3$,
$$
\sigma_1 T_1 T_2\, p_4p_5\, s_5\, B_1 B_2 \sigma_2
=
(\sigma_1 T_1 T_2\,p_4p_5\,\sigma_1^{-1})
(\sigma_1 s_5\sigma_2)
(\sigma_2^{-1} B_1 B_2 \sigma_2)
=
T_2' T_1'\,p_3p_4\, s_4\,  B_2' B_1',
$$
where
$T_2' T_1' = (p_{1{1\over2}}p_2)(p_{5{1\over2}}p_{6{1\over2}}p_6p_7)$
and $B_2'B_1'=(p_2p_3p_{1{1\over2}}p_{2{1\over2}})(p_5p_6p_7p_{4{1\over2}}
p_{5{1\over2}}p_{6{1\over2}})$.

\medskip\noindent
{\it Step 4}:  Let $t,b\in P_k$ and let $\pi\in S_k$. Then $t \pi b
= t x \sigma$ where  $x\in P_k$ and
$\sigma\in S_k$, and this transformation
can be acccomplished using the relations in (b), (c) and (d).

\smallskip
Suppose $T$ is a block of bottom dots of $t$ containing more
than one dot and which is connected, by edges of $\pi$, to two
top blocks $B_1$ and $B_2$ of $b$.
Using Step 3 find permutations
$\gamma_1,\gamma_2\in S_k$ and $\sigma_1,\sigma_2\in S_k$
such that
$$t' = \gamma_1 t \gamma_2
\qquad\hbox{and}\qquad
b' = \sigma_1 b \sigma_2
$$
are planar diagrams with $T$ as the leftmost bottom
block of $t'$ and $B_1$ and $B_2$ as the two leftmost top blocks
of $b'$.
Then
$$\eqalign{
t\pi b
= \gamma_1^{-1}t' \gamma_2^{-1}\pi \sigma_1^{-1} b' \sigma_2^{-1}
= \gamma_1^{-1}t' (\gamma_2^{-1}\pi \sigma_1^{-1}) b' \sigma_2^{-1}
= \gamma_1^{-1}t' (\gamma_2^{-1}\pi \sigma_1^{-1}) b'' \sigma_2^{-1}
= t\pi \sigma_1^{-1} b'' \sigma_2^{-1},
}$$
where $b''$ is a planar diagram with fewer top blocks than $b$ has.  This
is best seen from the following picture, where $t\pi b$ equals
$$
\matrix{
{\beginpicture
\setcoordinatesystem units <0.3cm,0.2cm>         
\setplotarea x from 1 to 10, y from -1 to 2    
\linethickness=0.5pt
\put{$\scriptstyle{T}$} at 4 .7
\setquadratic
\plot 3 -1 3.25 0 3.5 -1 /
\plot 3.5 -1 3.75 0 4 -1 /
\plot 4 -1 4.25 0 4.5 -1 /
\plot 4.5 -1 4.75 0 5 -1 /
\setlinear
\setdashpattern <.02cm,.05cm>
\put{$t$} at 0 .5
\plot 1 -1 1 2 10 2 10 -1 1 -1 /
\endpicture}  \cr
{\beginpicture
\setcoordinatesystem units <0.3cm,0.2cm>         
\setplotarea x from 1 to 10, y from -1 to 2    
\linethickness=0.5pt                         
\put{$\pi$} at .1 .5
\plot 4.5 2 2.5 -1 /
\plot 5 2 6.5 -1 /
\setdashpattern <.02cm,.05cm>
\plot 1 -1 1 2 10 2 10 -1 1 -1 /
\endpicture} \cr
{\beginpicture
\setcoordinatesystem units <0.3cm,0.2cm>         
\setplotarea x from 1 to 10, y from -1 to 2    
\linethickness=0.5pt                         
\put{$b$} at 0 .5
\put{$\scriptstyle{B_1}$} at 2.1 .2
\put{$\scriptstyle{B_2}$} at 7.5 .2
\setquadratic
\plot 1.5 2 1.75 1 2 2 /
\plot 2 2 2.25 1 2.5 2 /
\plot 6.5 2 6.75 1 7 2 /
\plot 7 2 7.25 1 7.5 2 /
\plot 7.5 2 7.75 1 8 2 /
\setlinear
\setdashpattern <.02cm,.05cm>
\plot 1 -1 1 2 10 2 10 -1 1 -1 /
\endpicture}}
=
\matrix{
{\beginpicture
\setcoordinatesystem units <0.3cm,0.2cm>         
\setplotarea x from 1 to 10, y from -1 to 2    
\linethickness=0.5pt
\setquadratic
\setlinear
\setdashpattern <.02cm,.05cm>
\put{$\scriptstyle{\gamma_1^{-1}}$} at 0.3 .5
\plot 1.1 -1 1.1 2 9.9 2 9.9 -1 9.9 -1 /
\endpicture}  \cr
{\beginpicture
\setcoordinatesystem units <0.3cm,0.2cm>         
\setplotarea x from 1 to 10, y from -1 to 2    
\linethickness=0.5pt
\put{$\scriptstyle{T}$} at 2.2 .7
\setquadratic
\plot 1.2 -1 1.45 0 1.7 -1 /
\plot 1.7 -1 1.95 0 2.2 -1 /
\plot 2.2 -1 2.45 0 2.7 -1 /
\plot 2.7 -1 2.95 0 3.2 -1 /
\setlinear
\setdashpattern <.02cm,.05cm>
\put{$t'$} at 0 .5
\plot 1.2 -1 1.2 2 9.96 2 9.96 -1 1.2 -1 /
\endpicture}  \cr
{\beginpicture
\setcoordinatesystem units <0.3cm,0.2cm>         
\setplotarea x from 1 to 10, y from -1 to 2    
\linethickness=0.5pt
\setquadratic
\setlinear
\plot 1.2 2 2.8 -1 /
\plot 1.7 2 3.3 -1 /
\plot 2.2 2 3.8 -1 /
\plot 2.7 2 4.2 -1 /
\plot 3.2 2 4.8 -1 /
\setdashpattern <.02cm,.05cm>
\put{$\scriptstyle{\gamma_2^{-1}}$} at 0.3 .5
\plot 1.1 -1 1.1 2 9.9 2 9.9 -1 1.1 -1 /
\endpicture}  \cr
{\beginpicture
\setcoordinatesystem units <0.3cm,0.2cm>         
\setplotarea x from 1 to 10, y from -1 to 2    
\linethickness=0.5pt                         
\put{$\pi$} at .1 .5
\plot 4.5 2 3.5 -1 /
\plot 5 2 6.5 -1 /
\setdashpattern <.02cm,.05cm>
\plot 1.2 -1 1.2 2 10.03 2 10.03 -1 1.2 -1 /
\endpicture} \cr
{\beginpicture
\setcoordinatesystem units <0.3cm,0.2cm>         
\setplotarea x from 1 to 10, y from -1 to 2    
\linethickness=0.5pt
\setquadratic
\setlinear
\plot 1.2 -1 2.3 2 /
\plot 1.7 -1 2.8 2 /
\plot 2.2 -1 3.3 2 /
\plot 2.7 -1 6.3 2 /
\plot 3.2 -1   6.8 2 /
\plot 3.7 -1 7.3 2 /
\plot 4.2 -1 7.8 2 /
\setdashpattern <.02cm,.05cm>
\put{$\scriptstyle{\sigma_1^{-1}}$} at 0.3 .5
\plot 1.1 -1 1.1 2 9.89 2 9.89 -1 1.1 -1 /
\endpicture}  \cr
{\beginpicture
\setcoordinatesystem units <0.3cm,0.2cm>         
\setplotarea x from 1 to 10, y from -1 to 2    
\linethickness=0.5pt                         
\put{$b'$} at 0 .5
\put{$\scriptstyle{B_1}$} at 1.8 .2
\put{$\scriptstyle{B_2}$} at 3.5 .2
\setquadratic
\plot 1.2 2 1.45 1 1.7 2 /
\plot 1.7 2 1.95 1 2.2 2 /
\plot 2.7 2 2.95 1 3.2 2 /
\plot 3.2 2 3.45 1 3.7 2 /
\plot 3.7 2 3.95 1 4.2 2 /
\setlinear
\setdashpattern <.02cm,.05cm>
\plot 1.2 -1 1.2 2 9.95 2 9.95 -1 1.2 -1 /
\endpicture}\cr
{\beginpicture
\setcoordinatesystem units <0.3cm,0.2cm>         
\setplotarea x from 1 to 10, y from -1 to 2    
\linethickness=0.5pt
\setquadratic
\setlinear
\setdashpattern <.02cm,.05cm>
\put{$\scriptstyle{\sigma_2^{-1}}$} at 0.3 .5
\plot 1.1 -1 1.1 2 9.89 2 9.89 -1 1.1 -1 /
\endpicture}  \cr}
=
\matrix{
{\beginpicture
\setcoordinatesystem units <0.3cm,0.2cm>         
\setplotarea x from 1 to 10, y from -1 to 2    
\linethickness=0.5pt
\setquadratic
\setlinear
\setdashpattern <.02cm,.05cm>
\put{$\scriptstyle{\gamma_1^{-1}}$} at 0.3 .5
\plot 1.1 -1 1.1 2 9.9 2 9.9 -1 1.1 -1 /
\endpicture}  \cr
{\beginpicture
\setcoordinatesystem units <0.3cm,0.2cm>         
\setplotarea x from 1 to 10, y from -1 to 2    
\linethickness=0.5pt
\put{$\scriptstyle{T}$} at 2.2 .7
\setquadratic
\plot 1.2 -1 1.45 0 1.7 -1 /
\plot 1.7 -1 1.95 0 2.2 -1 /
\plot 2.2 -1 2.45 0 2.7 -1 /
\plot 2.7 -1 2.95 0 3.2 -1 /
\setlinear
\setdashpattern <.02cm,.05cm>
\put{$t'$} at 0 .5
\plot 1.2 -1 1.2 2 9.96 2 9.96 -1 1.2 -1 /
\endpicture}  \cr
{\beginpicture
\setcoordinatesystem units <0.3cm,0.2cm>         
\setplotarea x from -.5 to 10, y from -1 to 2    
\linethickness=0.5pt
\setlinear
\plot 2.7 2 2.2 -1 /
\plot 3.2 2 2.7 -1 /
\put{$\scriptstyle{\gamma_2^{-1} \pi \sigma_1^{-1}}$} at 6 .5
\setdashpattern <.02cm,.05cm>
\plot 1.1 -1 1.1 2 9.9 2 9.9 -1 1.1 -1 /
\endpicture}  \cr
{\beginpicture
\setcoordinatesystem units <0.3cm,0.2cm>         
\setplotarea x from 1 to 10, y from -1 to 2    
\linethickness=0.5pt                         
\put{$b'$} at 0.2 .5
\setquadratic
\plot 1.2 2 1.45 1 1.7 2 /
\plot 1.7 2 1.95 1 2.2 2 /
\plot 2.7 2 2.95 1 3.2 2 /
\plot 3.2 2 3.45 1 3.7 2 /
\plot 3.7 2 3.95 1 4.2 2 /
\setlinear
\setdashpattern <.02cm,.05cm>
\plot 1.2 -1 1.2 2 10.05 2 10.05 -1 1.2 -1 /
\endpicture}\cr
{\beginpicture
\setcoordinatesystem units <0.3cm,0.2cm>         
\setplotarea x from 1 to 10, y from -1 to 2    
\linethickness=0.5pt
\setquadratic
\setlinear
\setdashpattern <.02cm,.05cm>
\put{$\scriptstyle{\sigma_2^{-1}}$} at 0.3 .5
\plot 1.1 -1 1.1 2 9.89 2 9.89 -1 1.1 -1 /
\endpicture}  \cr}
=
\matrix{
{\beginpicture
\setcoordinatesystem units <0.3cm,0.2cm>         
\setplotarea x from 1 to 10, y from -1 to 2    
\linethickness=0.5pt
\setquadratic
\setlinear
\setdashpattern <.02cm,.05cm>
\put{$\scriptstyle{\gamma_1^{-1}}$} at 0.3 .5
\plot 1.1 -1 1.1 2 9.9 2 9.9 -1 1.1 -1 /
\endpicture}  \cr
{\beginpicture
\setcoordinatesystem units <0.3cm,0.2cm>         
\setplotarea x from 1 to 10, y from -1 to 2    
\linethickness=0.5pt
\put{$\scriptstyle{T}$} at 2.2 .7
\setquadratic
\plot 1.2 -1 1.45 0 1.7 -1 /
\plot 1.7 -1 1.95 0 2.2 -1 /
\plot 2.2 -1 2.45 0 2.7 -1 /
\plot 2.7 -1 2.95 0 3.2 -1 /
\setlinear
\setdashpattern <.02cm,.05cm>
\put{$t'$} at 0 .5
\plot 1.2 -1 1.2 2 9.96 2 9.96 -1 1.2 -1 /
\endpicture}  \cr
{\beginpicture
\setcoordinatesystem units <0.3cm,0.2cm>         
\setplotarea x from -.5 to 10, y from -1 to 2    
\linethickness=0.5pt
\setlinear
\plot 2.7 2 2.2 -1 /
\plot 3.2 2 2.7 -1 /
\put{$\scriptstyle{\gamma_2^{-1} \pi \sigma_1^{-1}}$} at 6 .5
\setdashpattern <.02cm,.05cm>
\plot 1.1 -1 1.1 2 9.9 2 9.9 -1 1.1 -1 /
\endpicture}  \cr
{\beginpicture
\setcoordinatesystem units <0.3cm,0.2cm>         
\setplotarea x from 1 to 10, y from -1 to 2    
\linethickness=0.5pt                         
\put{$b''$} at 0.2 .5
\setquadratic
\plot 1.2 2 1.45 1 1.7 2 /
\plot 1.7 2 1.95 1 2.2 2 /
\plot 2.2 2 2.45 1 2.7 2 /
\plot 2.7 2 2.95 1 3.2 2 /
\plot 3.2 2 3.45 1 3.7 2 /
\plot 3.7 2 3.95 1 4.2 2 /
\setlinear
\setdashpattern <.02cm,.05cm>
\plot 1.2 -1 1.2 2 10 2 10 -1 1.2 -1 /
\endpicture}\cr
{\beginpicture
\setcoordinatesystem units <0.3cm,0.2cm>         
\setplotarea x from 1 to 10, y from -1 to 2    
\linethickness=0.5pt
\setquadratic
\setlinear
\setdashpattern <.02cm,.05cm>
\put{$\scriptstyle{\sigma_2^{-1}}$} at 0.3 .5
\plot 1.1 -1 1.1 2 9.9 2 9.9 -1 1.1 -1 /
\endpicture}  \cr}
$$
and the last equality uses the relations $p_{i+{1\over2}}
=p_{i+{1\over2}}^2$ and fourth relation in (d) (multiple times).
Then $t\pi b = \gamma_1^{-1}t'\gamma_2^{-1}\pi\sigma_1^{-1} b''\sigma_2^{-1}
=t\pi' b'' \sigma_2^{-1}$, where $\pi' = \pi\sigma_1^{-1}$.

By iteration of this process it is sufficient to assume that
in proving Step 4 we are analyzing $t\pi b$ where each bottom
block of $t$ connects
to a single top block of $b$.  Then, since $\pi$ is a permutation,
the bottom blocks of $t$ must have the same sizes as the top
blocks of $b$ and  $\pi$ is the permutation that permutes
the bottom blocks of $t$ to the top blocks of $b$.  Thus, by Step 1,
there is $\sigma\in S_k$ such that $x=\pi b \sigma^{-1}$ is planar
and
$$t\pi b = t (\pi b \sigma^{-1})\sigma= tx\sigma.$$

\smallskip\noindent
{\it Completion of the proof}: If $d_1, d_2\in A_k$ then use
the decomposition $A_k = S_kP_kS_k$ (from (1.6)) to write $d_1$ and $d_2$ in
the form
$$d_1 = \pi_1 t \pi_2
\qquad\hbox{and}\qquad d_2= \sigma_1 b\sigma_2,
\qquad\hbox{with }
t,b\in P_k, \pi_1,\pi_2, \sigma_1,\sigma_2\in S_k,
$$
and use (b) and (c) to write these products in terms of the generators.
Let $\pi = \pi_2\sigma_1$.  Then Step 4 tell us that the relations
give $\sigma\in S_k$ and $x\in P_k$ such that
$$
d_1 d_2 = \pi_1 t \pi_2 \sigma_1 b \sigma_2
= \pi_1 t \pi b \sigma_2
= \pi_1 t x \sigma \sigma_2,
$$
Using Step 2 and that $A_k = S_kP_kS_k$,
this product can be identified with the product
diagram $d_1d_2$.  Thus, the relations are sufficient to compose
any two elements of $A_k$.
\pfend

\section 2. Partition Algebras

\bigskip
For $k \in {1 \over 2} \ZZ_{> 0}$ and $n \in \CC$, the {\it partition
algebra} $\CC A_k(n)$ is the associative algebra over $\CC$ with basis
$A_k$,
$$
\CC A_k(n) = \CC\hbox{span-}\{ d\in A_k\},
\qquad\hbox{and multiplication defined by}\qquad
d_1d_2 = n^\ell (d_1 \circ d_2),$$
where, for $d_1,d_2\in A_k$,
$d_1 \circ d_2$ is the product in the monoid $A_k$
and $\ell$ is the number of
blocks removed from the the middle row when constructing
the composition $d_1 \circ d_2$.  For example,
$$\hbox{if}\qquad
{\beginpicture
\setcoordinatesystem units <0.5cm,0.2cm> 
\setplotarea x from 0 to 7, y from 0 to 3    
\linethickness=0.5pt
\put{$d_1 = $} at -.5 0
\put{$\bullet$} at 1 -1.5 \put{$\bullet$} at 1 1.5
\put{$\bullet$} at 2 -1.5 \put{$\bullet$} at 2 1.5
\put{$\bullet$} at 3 -1.5 \put{$\bullet$} at 3 1.5
\put{$\bullet$} at 4 -1.5 \put{$\bullet$} at 4 1.5
\put{$\bullet$} at 5 -1.5 \put{$\bullet$} at 5 1.5
\put{$\bullet$} at 6 -1.5 \put{$\bullet$} at 6 1.5
\put{$\bullet$} at 7 -1.5 \put{$\bullet$} at 7 1.5
\setquadratic
\plot 1 1.5 2 0.75 3 1.5 /
\plot 1 1.5 2 0.75 3 1.5 /
\plot 4 1.5 4.5 0.75 5 1.5 /
\plot 5 1.5 5.5 0.75 6 1.5 /
\plot 2 -1.5 2.5 -0.25 3 -1.5 /
\plot 5 -1.5 6 -0.25 7 -1.5 /
\setlinear
\plot 3 1.5 4 -1.5 /
\endpicture}
\quad\hbox{and}\quad
{\beginpicture
\setcoordinatesystem units <0.5cm,0.2cm> 
\setplotarea x from 0 to 7, y from 0 to 3    
\linethickness=0.5pt
\put{$d_2 = $} at -.5 0
\put{$\bullet$} at 1 1.5 \put{$\bullet$} at 1 -1.5
\put{$\bullet$} at 2 1.5 \put{$\bullet$} at 2 -1.5
\put{$\bullet$} at 3 1.5 \put{$\bullet$} at 3 -1.5
\put{$\bullet$} at 4 1.5 \put{$\bullet$} at 4 -1.5
\put{$\bullet$} at 5 1.5 \put{$\bullet$} at 5 -1.5
\put{$\bullet$} at 6 1.5 \put{$\bullet$} at 6 -1.5
\put{$\bullet$} at 7 1.5 \put{$\bullet$} at 7 -1.5
\setquadratic
\plot 2 1.5 3 0.25 4 1.5 /
\plot 5 1.5 6 0.25 7 1.5 /
\plot 2 -1.5 4.5 -0.5 7 -1.5 /
\plot 4 -1.5 4.5 -1 5 -1.5 /
\plot 5 -1.5 5.5 -1 6 -1.5 /
\setlinear
\plot 3 1.5 4 -1.5 /
\plot 6 1.5 7 -1.5 /
\endpicture}
\quad\hbox{then}
$$
\medskip
$$
{\beginpicture
\setcoordinatesystem units <0.4cm,0.15cm> 
\setplotarea x from 0 to 7, y from 0 to 3    
\linethickness=0.5pt
\put{$d_1d_2 = $} at -.75 0
\put{$\bullet$} at 1 1 \put{$\bullet$} at 1 4
\put{$\bullet$} at 2 1 \put{$\bullet$} at 2 4
\put{$\bullet$} at 3 1 \put{$\bullet$} at 3 4
\put{$\bullet$} at 4 1 \put{$\bullet$} at 4 4
\put{$\bullet$} at 5 1 \put{$\bullet$} at 5 4
\put{$\bullet$} at 6 1 \put{$\bullet$} at 6 4
\put{$\bullet$} at 7 1 \put{$\bullet$} at 7 4
\put{$\bullet$} at 1 -1 \put{$\bullet$} at 1 -4
\put{$\bullet$} at 2 -1 \put{$\bullet$} at 2 -4
\put{$\bullet$} at 3 -1 \put{$\bullet$} at 3 -4
\put{$\bullet$} at 4 -1 \put{$\bullet$} at 4 -4
\put{$\bullet$} at 5 -1 \put{$\bullet$} at 5 -4
\put{$\bullet$} at 6 -1 \put{$\bullet$} at 6 -4
\put{$\bullet$} at 7 -1 \put{$\bullet$} at 7 -4
\setquadratic
\plot 1 4 2 3.25 3 4 /
\plot 1 4 2 3.25 3 4 /
\plot 4 4 4.5 3.25 5 4 /
\plot 5 4 5.5 3.25 6 4 /
\plot 2 1 2.5 1.75 3 1 /
\plot 5 1 6 1.75 7 1 /
\plot 2 -1 3 -1.75 4 -1 /
\plot 5 -1 6 -1.75 7 -1 /
\plot 2 -4 4.5 -3 7 -4 /
\plot 4 -4 4.5 -3.5 5 -4 /
\plot 5 -4 5.5 -3.5 6 -4 /
\setlinear
\plot 3 4 4 1 /
\plot 3 -1 4 -4 /
\plot 6 -1 7 -4 /
\setdashpattern <.03cm,.03cm>
\plot 1 1 1 -1 /
\plot 2 1 2 -1 /
\plot 3 1 3 -1 /
\plot 4 1 4 -1 /
\plot 5 1 5 -1 /
\plot 6 1 6 -1 /
\plot 7 1 7 -1 /
\endpicture}
~=
~~n^2~ {\beginpicture
\setcoordinatesystem units <0.4cm,0.15cm> 
\setplotarea x from 1 to 7, y from 0 to 3    
\linethickness=0.5pt                        
\put{$\bullet$} at 1 -1 \put{$\bullet$} at 1 2
\put{$\bullet$} at 2 -1 \put{$\bullet$} at 2 2
\put{$\bullet$} at 3 -1 \put{$\bullet$} at 3 2
\put{$\bullet$} at 4 -1 \put{$\bullet$} at 4 2
\put{$\bullet$} at 5 -1 \put{$\bullet$} at 5 2
\put{$\bullet$} at 6 -1 \put{$\bullet$} at 6 2
\put{$\bullet$} at 7 -1 \put{$\bullet$} at 7 2
\plot 3 2 4 -1 /
\setquadratic
\plot 1 2 2 1.25 3 2 /
\plot 1 2 2 1.25 3 2 /
\plot 4 2 4.5 1.25 5 2 /
\plot 5 2 5.5 1.25 6 2 /
\plot 2 -1 4.5 0 7 -1 /
\plot 4 -1 4.5 -.5 5 -1 /
\plot 5 -1 5.5 -.5 6 -1 /
\endpicture},\formula
$$
since two blocks are removed from the middle row.
There are inclusions of algebras given by
$$\matrix{
\CC \Akm &\hookrightarrow &\CC A_k  \cr\cr
{\beginpicture
\setcoordinatesystem units <0.5cm,0.2cm> 
\setplotarea x from 1 to 4, y from -1 to 3    
\linethickness=0.5pt                        
\setdashes <.04cm>
\plot 1 -1   1  2 /
\plot 4 -1   4  2 /
\plot 1 -1   4 -1 /
\plot 1  2   4  2 /
\put{$\scriptstyle{1}$}[b] at 1 2.5
\put{$d$} at 2.5 .5
\put{$\scriptstyle{k}$}[b] at 4 2.5
\put{$\bullet$} at 1 -1 \put{$\bullet$} at 1 2
\put{$\bullet$} at 4 -1 \put{$\bullet$} at 4 2
\put{$\bullet$} at 1.5 -1 \put{$\bullet$} at 1.5 2
\put{$\bullet$} at 2 -1 \put{$\bullet$} at 2 2
\put{$\bullet$} at 2.5 -1 \put{$\bullet$} at 2.5 2
\put{$\bullet$} at 3 -1 \put{$\bullet$} at 3 2
\put{$\bullet$} at 3.5 -1 \put{$\bullet$} at 3.5 2
\endpicture} & \mapsto & {\beginpicture
\setcoordinatesystem units <0.5cm,0.2cm> 
\setplotarea x from 1 to 4, y from -1 to 3    
\linethickness=0.5pt                        
\setdashes <.04cm>
\plot 1 -1   1  2 /
\plot 4 -1   4  2 /
\plot 1 -1   4 -1 /
\plot 1  2   4  2 /
\put{$\scriptstyle{1}$}[b] at 1 2.5
\put{$d$} at 2.5 .5
\put{$\scriptstyle{k}$}[b] at 4 2.5
\put{$\bullet$} at 1 -1 \put{$\bullet$} at 1 2
\put{$\bullet$} at 4 -1 \put{$\bullet$} at 4 2
\put{$\bullet$} at 1.5 -1 \put{$\bullet$} at 1.5 2
\put{$\bullet$} at 2 -1 \put{$\bullet$} at 2 2
\put{$\bullet$} at 2.5 -1 \put{$\bullet$} at 2.5 2
\put{$\bullet$} at 3 -1 \put{$\bullet$} at 3 2
\put{$\bullet$} at 3.5 -1 \put{$\bullet$} at 3.5 2
\endpicture} \cr }
\qquad\hbox{and}\qquad
\matrix{
\CC A_{k-1} &\hookrightarrow& \CC \Akm \cr\cr
{\beginpicture
\setcoordinatesystem units <0.5cm,0.2cm> 
\setplotarea x from 1 to 4, y from -1 to 3    
\linethickness=0.5pt                        
\setdashes <.04cm>
\plot 1 -1   1  2 /
\plot 3.5 -1   3.5  2 /
\plot 1 -1   3.5 -1 /
\plot 1  2   3.5  2 /
\put{$\scriptstyle{1}$}[b] at 1 2.5
\put{$d$} at 2.5 .5
\put{$\scriptstyle{k-1}$}[b] at 3.7 2.5
\put{$\bullet$} at 1 -1 \put{$\bullet$} at 1 2
\put{$\bullet$} at 1.5 -1 \put{$\bullet$} at 1.5 2
\put{$\bullet$} at 2 -1 \put{$\bullet$} at 2 2
\put{$\bullet$} at 2.5 -1 \put{$\bullet$} at 2.5 2
\put{$\bullet$} at 3 -1 \put{$\bullet$} at 3 2
\put{$\bullet$} at 3.5 -1 \put{$\bullet$} at 3.5 2
\endpicture}
&\mapsto
&{\beginpicture
\setcoordinatesystem units <0.5cm,0.2cm> 
\setplotarea x from 1 to 4, y from -1 to 3    
\linethickness=0.5pt                        
\plot 4 -1   4  2 /
\setdashes <.04cm>
\plot 1 -1   1  2 /
\plot 3.5 -1   3.5  2 /
\plot 4 -1   4  2 /
\plot 1 -1   3.5 -1 /
\plot 1  2   3.5  2 /
\put{$\scriptstyle{1}$}[b] at 1 2.5
\put{$d$} at 2.5 .5
\put{$\scriptstyle{k}$}[b] at 4 2.5
\put{$\bullet$} at 1 -1 \put{$\bullet$} at 1 2
\put{$\bullet$} at 4 -1 \put{$\bullet$} at 4 2
\put{$\bullet$} at 1.5 -1 \put{$\bullet$} at 1.5 2
\put{$\bullet$} at 2 -1 \put{$\bullet$} at 2 2
\put{$\bullet$} at 2.5 -1 \put{$\bullet$} at 2.5 2
\put{$\bullet$} at 3 -1 \put{$\bullet$} at 3 2
\put{$\bullet$} at 3.5 -1 \put{$\bullet$} at 3.5 2
\put{.} at 4.5 -1
\endpicture} \cr } \formula
$$
For $d_1, d_2 \in A_k$, define
$$d_1 \le d_2, \qquad\hbox{
if the set partition $d_2$ is coarser than the set partition $d_1$,}$$
i.e., $i$ and $j$ in the same block of $d_1$
implies that $i$ and $j$ are in the same block of $d_2$.
Let $\{ x_d \in \CC A_k | d \in A_k \}$ be the basis of
$\CC A_k$ uniquely defined by the relation
$$
d = \sum_{d' \le d} x_{d'}, \qquad \hbox{for all $d \in A_k$}. \formula
$$
Under any linear extension of the partial order $\le$
the transition matrix
between the basis $\{d\ |\  d\in A_k\}$ of $\CC A_k(n)$ and
the basis $\{ x_d \  |\  d \in A_k \}$ of $\CC A_k(n)$ is
upper triangular with 1s on the diagonal and so the $x_d$ are well defined.

\bigskip\noindent
{\it The maps
$\varepsilon_{1 \over 2}\colon \CC A_k \to \CC \Akm$,
$\varepsilon^{1 \over 2}\colon \CC \Akm \to \CC A_{k-1}$ and
$tr_k\colon \CC A_k \to \CC$}
\medskip

Let $k\in \ZZ_{>0}$.
Define linear maps
$$\matrix{
\varepsilon_{1 \over 2}: & \CC A_k &\longrightarrow &\CC \Akm  \cr\cr
&{\beginpicture
\setcoordinatesystem units <0.5cm,0.2cm> 
\setplotarea x from 1 to 4, y from -1 to 3    
\linethickness=0.5pt                        
\setdashes <.04cm>
\plot 1 -1   1  2 /
\plot 4 -1   4  2 /
\plot 1 -1   4 -1 /
\plot 1  2   4  2 /
\put{$\scriptstyle{1}$}[b] at 1 2.5
\put{$d$} at 2.5 .5
\put{$\scriptstyle{k}$}[b] at 4 2.5
\put{$\bullet$} at 1 -1 \put{$\bullet$} at 1 2
\put{$\bullet$} at 4 -1 \put{$\bullet$} at 4 2
\put{$\bullet$} at 1.5 -1 \put{$\bullet$} at 1.5 2
\put{$\bullet$} at 2 -1 \put{$\bullet$} at 2 2
\put{$\bullet$} at 2.5 -1 \put{$\bullet$} at 2.5 2
\put{$\bullet$} at 3 -1 \put{$\bullet$} at 3 2
\put{$\bullet$} at 3.5 -1 \put{$\bullet$} at 3.5 2
\endpicture} & \mapsto & {\beginpicture
\setcoordinatesystem units <0.5cm,0.2cm> 
\setplotarea x from 1 to 4, y from -1 to 3    
\linethickness=0.5pt                        
\setquadratic
\plot 4 2 4.25 0.5 4 -1 /
\setlinear
\setdashes <.04cm>
\plot 1 -1   1  2 /
\plot 4 -1   4  2 /
\plot 1 -1   4 -1 /
\plot 1  2   4  2 /
\put{$\scriptstyle{1}$}[b] at 1 2.5
\put{$d$} at 2.5 .5
\put{$\scriptstyle{k}$}[b] at 4 2.5
\put{$\bullet$} at 1 -1 \put{$\bullet$} at 1 2
\put{$\bullet$} at 4 -1 \put{$\bullet$} at 4 2
\put{$\bullet$} at 1.5 -1 \put{$\bullet$} at 1.5 2
\put{$\bullet$} at 2 -1 \put{$\bullet$} at 2 2
\put{$\bullet$} at 2.5 -1 \put{$\bullet$} at 2.5 2
\put{$\bullet$} at 3 -1 \put{$\bullet$} at 3 2
\put{$\bullet$} at 3.5 -1 \put{$\bullet$} at 3.5 2
\endpicture} \cr }
\qquad\hbox{and}\qquad
\matrix{
\varepsilon^{1 \over 2}:& \CC \Akm &\longrightarrow& \CC A_{k-1} \cr\cr
&{\beginpicture
\setcoordinatesystem units <0.5cm,0.2cm> 
\setplotarea x from 1 to 4, y from -1 to 3    
\linethickness=0.5pt                        
\setquadratic
\plot 4 2 4.25 0.5 4 -1 /
\setlinear
\setdashes <.04cm>
\plot 1 -1   1  2 /
\plot 4 -1   4  2 /
\plot 1 -1   4 -1 /
\plot 1  2   4  2 /
\put{$\scriptstyle{1}$}[b] at 1 2.5
\put{$d$} at 2.5 .5
\put{$\scriptstyle{k}$}[b] at 4 2.5
\put{$\bullet$} at 1 -1 \put{$\bullet$} at 1 2
\put{$\bullet$} at 4 -1 \put{$\bullet$} at 4 2
\put{$\bullet$} at 1.5 -1 \put{$\bullet$} at 1.5 2
\put{$\bullet$} at 2 -1 \put{$\bullet$} at 2 2
\put{$\bullet$} at 2.5 -1 \put{$\bullet$} at 2.5 2
\put{$\bullet$} at 3 -1 \put{$\bullet$} at 3 2
\put{$\bullet$} at 3.5 -1 \put{$\bullet$} at 3.5 2
\endpicture}
&\mapsto
&{\beginpicture
\setcoordinatesystem units <0.5cm,0.2cm> 
\setplotarea x from 1 to 4, y from -1 to 3    
\linethickness=0.5pt                        
\setquadratic
\plot 4 2 4.25 0.5 4 -1 /
\setlinear
\setdashes <.04cm>
\plot 1 -1   1  2 /
\plot 4 -1   4  2 /
\plot 1 -1   4 -1 /
\plot 1  2   4  2 /
\put{$\scriptstyle{1}$}[b] at 1 2.5
\put{$d$} at 2.5 .5
\put{$\scriptstyle{k-1}$}[b] at 3.7 2.5
\put{$\bullet$} at 1 -1 \put{$\bullet$} at 1 2
\put{$\bullet$} at 1.5 -1 \put{$\bullet$} at 1.5 2
\put{$\bullet$} at 2 -1 \put{$\bullet$} at 2 2
\put{$\bullet$} at 2.5 -1 \put{$\bullet$} at 2.5 2
\put{$\bullet$} at 3 -1 \put{$\bullet$} at 3 2
\put{$\bullet$} at 3.5 -1 \put{$\bullet$} at 3.5 2
\endpicture} \cr }
$$
so that $\varepsilon_{1\over2}(d)$ is the same as $d$ except
that the block containing $k$ and the block
containing $k'$ are combined, and $\varepsilon^{1\over2}(d)$
has the same blocks as $d$ except with $k$ and $k'$ removed.
There is a factor of $n$ in $\varepsilon^{1\over2}(d)$ if the
removal of $k$ and $k'$ reduces the number of blocks by $1$.
For example,
$$
\varepsilon_{1\over 2}\left({\beginpicture
\setcoordinatesystem units <0.4cm,0.15cm> 
\setplotarea x from 1 to 5, y from -1 to 3    
\linethickness=0.5pt                        
\put{$\bullet$} at 1 -1 \put{$\bullet$} at 1 2
\put{$\bullet$} at 2 -1 \put{$\bullet$} at 2 2
\put{$\bullet$} at 3 -1 \put{$\bullet$} at 3 2
\put{$\bullet$} at 4 -1 \put{$\bullet$} at 4 2
\put{$\bullet$} at 5 -1 \put{$\bullet$} at 5 2
\plot 2 -1 4 2 /
\plot 1 2 3 -1 /
\plot 3 2 3 -1 /
\setquadratic
\plot 1 2 2 1.25 3 2 /
\plot 4 -1 4.5 -.5 5 -1 /
\endpicture} \right) = {\beginpicture
\setcoordinatesystem units <0.4cm,0.15cm> 
\setplotarea x from 1 to 5, y from -1 to 3    
\linethickness=0.5pt                        
\put{$\bullet$} at 1 -1 \put{$\bullet$} at 1 2
\put{$\bullet$} at 2 -1 \put{$\bullet$} at 2 2
\put{$\bullet$} at 3 -1 \put{$\bullet$} at 3 2
\put{$\bullet$} at 4 -1 \put{$\bullet$} at 4 2
\put{$\bullet$} at 5 -1 \put{$\bullet$} at 5 2
\plot 2 -1 4 2 /
\plot 1 2 3 -1 /
\plot 3 2 3 -1 /
\plot 5 2 5 -1 /
\setquadratic
\plot 1 2 2 1.25 3 2 /
\plot 4 -1 4.5 -.5 5 -1 /
\endpicture},
\qquad\qquad
\varepsilon^{1\over 2}\left(
{\beginpicture
\setcoordinatesystem units <0.4cm,0.15cm> 
\setplotarea x from 1 to 5, y from -1 to 3    
\linethickness=0.5pt                        
\put{$\bullet$} at 1 -1 \put{$\bullet$} at 1 2
\put{$\bullet$} at 2 -1 \put{$\bullet$} at 2 2
\put{$\bullet$} at 3 -1 \put{$\bullet$} at 3 2
\put{$\bullet$} at 4 -1 \put{$\bullet$} at 4 2
\put{$\bullet$} at 5 -1 \put{$\bullet$} at 5 2
\plot 2 -1 4 2 /
\plot 1 2 3 -1 /
\plot 3 2 3 -1 /
\plot 5 2 5 -1 /
\setquadratic
\plot 1 2 2 1.25 3 2 /
\plot 4 -1 4.5 -.5 5 -1 /
\endpicture}
\right)
= {\beginpicture
\setcoordinatesystem units <0.4cm,0.15cm> 
\setplotarea x from 1 to 4, y from -1 to 3    
\linethickness=0.5pt                        
\put{$\bullet$} at 1 -1 \put{$\bullet$} at 1 2
\put{$\bullet$} at 2 -1 \put{$\bullet$} at 2 2
\put{$\bullet$} at 3 -1 \put{$\bullet$} at 3 2
\put{$\bullet$} at 4 -1 \put{$\bullet$} at 4 2
\plot 2 -1 4 2 /
\plot 1 2 3 -1 /
\plot 3 2 3 -1 /
\setquadratic
\plot 1 2 2 1.25 3 2 /
\endpicture},
$$
and
$$
\varepsilon_{1\over 2}\left({\beginpicture
\setcoordinatesystem units <0.4cm,0.15cm> 
\setplotarea x from 1 to 5, y from -1 to 3    
\linethickness=0.5pt                        
\put{$\bullet$} at 1 -1 \put{$\bullet$} at 1 2
\put{$\bullet$} at 2 -1 \put{$\bullet$} at 2 2
\put{$\bullet$} at 3 -1 \put{$\bullet$} at 3 2
\put{$\bullet$} at 4 -1 \put{$\bullet$} at 4 2
\put{$\bullet$} at 5 -1 \put{$\bullet$} at 5 2
\plot 2 -1 4 2 /
\plot 1 2 3 -1 /
\plot 3 2 3 -1 /
\setquadratic
\plot 1 2 2 1.25 3 2 /
\endpicture} \right)
= {\beginpicture
\setcoordinatesystem units <0.4cm,0.15cm> 
\setplotarea x from 1 to 5, y from -1 to 3    
\linethickness=0.5pt                        
\put{$\bullet$} at 1 -1 \put{$\bullet$} at 1 2
\put{$\bullet$} at 2 -1 \put{$\bullet$} at 2 2
\put{$\bullet$} at 3 -1 \put{$\bullet$} at 3 2
\put{$\bullet$} at 4 -1 \put{$\bullet$} at 4 2
\put{$\bullet$} at 5 -1 \put{$\bullet$} at 5 2
\plot 2 -1 4 2 /
\plot 1 2 3 -1 /
\plot 3 2 3 -1 /
\plot 5 2 5 -1 /
\setquadratic
\plot 1 2 2 1.25 3 2 /
\endpicture},
\qquad\qquad
\varepsilon^{1\over 2}\left(
{\beginpicture
\setcoordinatesystem units <0.4cm,0.15cm> 
\setplotarea x from 1 to 5, y from -1 to 3    
\linethickness=0.5pt                        
\put{$\bullet$} at 1 -1 \put{$\bullet$} at 1 2
\put{$\bullet$} at 2 -1 \put{$\bullet$} at 2 2
\put{$\bullet$} at 3 -1 \put{$\bullet$} at 3 2
\put{$\bullet$} at 4 -1 \put{$\bullet$} at 4 2
\put{$\bullet$} at 5 -1 \put{$\bullet$} at 5 2
\plot 2 -1 4 2 /
\plot 1 2 3 -1 /
\plot 3 2 3 -1 /
\plot 5 2 5 -1 /
\setquadratic
\plot 1 2 2 1.25 3 2 /
\endpicture}
\right)
= n {\beginpicture
\setcoordinatesystem units <0.4cm,0.15cm> 
\setplotarea x from 1 to 4, y from -1 to 3    
\linethickness=0.5pt                        
\put{$\bullet$} at 1 -1 \put{$\bullet$} at 1 2
\put{$\bullet$} at 2 -1 \put{$\bullet$} at 2 2
\put{$\bullet$} at 3 -1 \put{$\bullet$} at 3 2
\put{$\bullet$} at 4 -1 \put{$\bullet$} at 4 2
\plot 2 -1 4 2 /
\plot 1 2 3 -1 /
\plot 3 2 3 -1 /
\setquadratic
\plot 1 2 2 1.25 3 2 /
\endpicture}.
$$
The map $\varepsilon^{1 \over 2}$ is the composition
$\CC \Akm \hookrightarrow \CC A_k \mapright{\varepsilon_1} \CC A_{k-1}.$
The composition of $\varepsilon_{1 \over 2}$ and
$\varepsilon^{1\over 2}$ is the map
$$
\matrix{
\varepsilon_{1}: & \CC A_k &\longrightarrow &\CC A_{k-1}  \cr\cr
&{\beginpicture
\setcoordinatesystem units <0.5cm,0.2cm> 
\setplotarea x from 1 to 4, y from -1 to 3    
\linethickness=0.5pt                        
\setdashes <.04cm>
\plot 1 -1   1  2 /
\plot 4 -1   4  2 /
\plot 1 -1   4 -1 /
\plot 1  2   4  2 /
\put{$\scriptstyle{1}$}[b] at 1 2.5
\put{$d$} at 2.5 .5
\put{$\scriptstyle{k}$}[b] at 4 2.5
\put{$\bullet$} at 1 -1 \put{$\bullet$} at 1 2
\put{$\bullet$} at 4 -1 \put{$\bullet$} at 4 2
\put{$\bullet$} at 1.5 -1 \put{$\bullet$} at 1.5 2
\put{$\bullet$} at 2 -1 \put{$\bullet$} at 2 2
\put{$\bullet$} at 2.5 -1 \put{$\bullet$} at 2.5 2
\put{$\bullet$} at 3 -1 \put{$\bullet$} at 3 2
\put{$\bullet$} at 3.5 -1 \put{$\bullet$} at 3.5 2
\endpicture} & \mapsto &  {\beginpicture
\setcoordinatesystem units <0.5cm,0.2cm> 
\setplotarea x from 1 to 4, y from -1 to 3    
\linethickness=0.5pt                        
\setquadratic
\plot 4 2 4.25 0.5 4 -1 /
\setlinear
\setdashes <.04cm>
\plot 1 -1   1  2 /
\plot 4 -1   4  2 /
\plot 1 -1   4 -1 /
\plot 1  2   4  2 /
\put{$\scriptstyle{1}$}[b] at 1 2.5
\put{$d$} at 2.5 .5
\put{$\scriptstyle{k-1}$}[b] at 3.7 2.5
\put{$\bullet$} at 1 -1 \put{$\bullet$} at 1 2
\put{$\bullet$} at 1.5 -1 \put{$\bullet$} at 1.5 2
\put{$\bullet$} at 2 -1 \put{$\bullet$} at 2 2
\put{$\bullet$} at 2.5 -1 \put{$\bullet$} at 2.5 2
\put{$\bullet$} at 3 -1 \put{$\bullet$} at 3 2
\put{$\bullet$} at 3.5 -1 \put{$\bullet$} at 3.5 2
\endpicture} \cr }.
\formula
$$
By drawing diagrams it is straightforward to check that, for
$k\in \ZZ_{>0}$,
$$
\eqalign{
\varepsilon_{1 \over 2}(a_1 b a_2) & =  a_1 \varepsilon_{1 \over 2}(b) a_2,
\quad
\hbox{for $a_1,a_2 \in \Akm, b \in A_k$} \cr
\varepsilon^{1 \over 2}(a_1 b a_2) &  =  a_1 \varepsilon^{1 \over 2}(b)
a_2, \quad
\hbox{for $a_1,a_2 \in A_{k-1}, b \in \Akm$} \cr
\varepsilon_1(a_1 b a_2) & = a_1 \varepsilon_1(b) a_2, \quad
\hbox{for $a_1,a_2 \in A_{k-1}, b \in A_k$}. \cr
} \formula
$$
and
$$
\eqalign{
\pkp b \pkp & = \varepsilon_{1 \over 2} (b) \pkp
= \pkp \varepsilon_{1 \over 2} (b),
\quad\hbox{for $b \in A_k$} \cr
p_k  b p_k & = \varepsilon^{1 \over 2} (b) p_k
= p_k \varepsilon^{1 \over 2} (b),
\quad\hbox{for $b \in \Akm$} \cr
e_k b e_k & = \varepsilon_{1} (b) e_k = e_k \varepsilon_{1} (b),
\quad\hbox{for $b \in A_k$}. \cr
} \formula
$$

Define $tr_k\colon \CC A_k \to \CC$
and $tr_{k - {1\over 2}}\colon \CC \Akm \to \CC$ by
the equations
$$\tr_k(b) = \tr_{k-{1\over2}}(\varepsilon_{1\over2}(b)),
\quad\hbox{for $b\in A_k$},
\quad\hbox{ and }\quad
tr_{k - {1 \over 2}}(b)
=\tr_{k-1}(\varepsilon^{1\over2}(b)),
\quad\hbox{for $b \in \Akm$},
\formula
$$
so that
$$
tr_k(b) = \varepsilon_1^k(b), \quad\hbox{for $b \in A_k$},
\qquad\hbox{ and }\qquad
tr_{k - {1 \over 2}}(b) = \varepsilon_1^{k-1} \varepsilon^{1 \over 2}(b),
\quad\hbox{for $b \in \Akm$}.
\formula
$$
Pictorially $tr_k(d) = n^c$ where $c$ is the number of connected
components in the closure of the diagram $d$,
$$
tr_k(d) =
{\beginpicture
\setcoordinatesystem units <0.5cm,0.16cm> 
\setplotarea x from 1 to 8.2, y from -5 to 7    
\linethickness=0.5pt                        
\plot 1 -1   1  2 /
\plot 4 -1   4  2 /
\plot 1 -1   4 -1 /
\plot 1  2   4  2 /
\put{$d$} at 2.5 .5
\put{$\bullet$} at 1 -1 \put{$\bullet$} at 1 2
\put{$\bullet$} at 1.5 -1 \put{$\bullet$} at 1.5 2
\put{$\bullet$} at 2 -1 \put{$\bullet$} at 2 2
\put{$\bullet$} at 2.5 -1 \put{$\bullet$} at 2.5 2
\put{$\bullet$} at 3 -1 \put{$\bullet$} at 3 2
\put{$\bullet$} at 3.5 -1 \put{$\bullet$} at 3.5 2
\put{$\bullet$} at 4 -1 \put{$\bullet$} at 4 2
\ellipticalarc axes ratio 1:1.2 -280 degrees from 4 2 center at 4.5 .5
\ellipticalarc axes ratio 1:1.1 -304 degrees from 3.5 2 center at 4.5 .5
\ellipticalarc axes ratio 1:1.0 -320 degrees from 3 2 center at 4.5 .5
\ellipticalarc axes ratio 1:.9 -325 degrees from 2.5 2 center at 4.5 .5
\ellipticalarc axes ratio 1:.85 -330 degrees from 2 2 center at 4.5 .5
\ellipticalarc axes ratio 1:.8 -330 degrees from 1.5 2 center at 4.5 .5
\ellipticalarc axes ratio 1:.75 -335 degrees from 1 2 center at 4.5 .5
\endpicture},
\qquad\hbox{for $d \in A_k$}.
\formula
$$

\bigskip\noindent
{\it The ideal $\CC I_k(n)$}
\medskip

For $k\in {1\over2}\ZZ_{\ge 0}$ define
$$\CC I_k(n) = \CC\hbox{-span}\{ d\in I_k\}.
\formula$$
By (1.3),
$$\hbox{$\CC I_k(n)$ is an ideal of $\CC A_k(n)$}
\qquad\hbox{and}\qquad
\CC A_k(n)/\CC I_k(n) \cong \CC S_k,
\formula$$
since the set partitions with propagating number $k$
are exactly the permutations in the symmetric group
$S_k$ (by convention $S_{\ell+{1\over2}}=S_\ell$
for $\ell\in \ZZ_{>0}$, see (2.2)).

View $\CC I_k(n)$ as an algebra (without identity).
Since $\CC A_k(n)/\CC I_k\cong \CC S_k$ and $\CC S_k$ is
semisimple, $\Rad(\CC A_k(n))\subseteq \CC I_k(n)$.
Since $\CC I_k(n)/\Rad(\CC A_k(n))$ is an ideal in
$\CC A_k(n)/\Rad(\CC A_k(n))$ the quotient
$\CC I_k(n)/\Rad(\CC A_k(n))$ is semisimple.  Therefore
$\Rad(\CC I_k(n))\subseteq \Rad(\CC A_k(n))$.
On the other hand,
since $\Rad(\CC A_k(n))$ is an ideal of nilpotent elements
in $\CC A_k(n)$,
it is an ideal of nilpotent elements in $\CC I_k(n)$ and so
$\Rad(\CC I_k(n)) \supseteq \Rad(\CC A_k(n))$.
Thus
$$\Rad(\CC A_k(n))=\Rad(\CC I_k(n)).
\formula$$

Let $k\in \ZZ_{\ge 0}$.
By (2.5) the maps
$$
\varepsilon_{1 \over 2}: \CC A_k \longrightarrow \CC \Akm
\qquad\hbox{and}\qquad
\varepsilon^{1 \over 2}: \CC \Akm \longrightarrow \CC A_{k-1}
$$
are $(\CC \Akm, \CC \Akm)$-bimodule
and $(\CC A_{k-1}, \CC A_{k-1})$-bimodule homomorphisms,
respectively.
The corresponding {\it basic constructions} (see Section 4)
are the algebras
$$
\CC A_k(n) \otimes_{\CC \Akm(n)} \CC A_k(n)
\qquad\hbox{and}\qquad
\CC \Akm(n) \otimes_{\CC A_{k-1}(n)} \CC \Akm(n)
\formula
$$
with products given by
$$
(b_1\otimes b_2)(b_3\otimes b_4) =
b_1\otimes \varepsilon_{1\over2}(b_2b_3)b_4,
\qquad\hbox{and}\qquad
(c_1\otimes c_2)(c_3\otimes c_4) =
c_1\otimes \varepsilon^{1\over2}(c_2c_3)c_4,
\formula
$$
for $b_1,b_2,b_3,b_4\in \CC A_k(n)$, and
for $c_1,c_2,c_3,c_4\in \CC \Akm(n)$.

Let $k \in {1 \over 2} \ZZ_{>0}$.
Then, by the relations in (2.6) and the fact that
$$\hbox{every $d\in I_k$}
\qquad\hbox{can be written as}\qquad
d=d_1p_kd_2,\qquad\hbox{with $d_1,d_2\in \Akm$,}
\formula
$$
the maps
$$
\matrix{
\CC \Akm(n) \otimes_{\CC A_{k-1}(n)} \CC \Akm(n)
&\longrightarrow &\CC I_k(n)  \cr
\cr
b_1 \otimes b_2 &\longmapsto  &b_1 p_k b_2 }
\formula
$$
are algebra isomorphisms.
Thus the ideal $\CC I_k(n)$ is always isomorphic to a 
basic construction (in the sense of Section 4).

\bigskip\noindent
{\it Representations of the symmetric group}
\medskip

A partition $\lambda$ is a collection of boxes in a  corner.
We shall conform to the conventions in [Mac] and assume
that gravity goes up and to the left, i.e.,
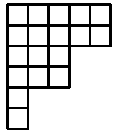
$$
\beginpicture
\setcoordinatesystem units <0.25cm,0.25cm>         
\setplotarea x from 0 to 4, y from 0 to 3    
\linethickness=0.5pt                          
\putrule from 0 6 to 5 6          %
\putrule from 0 5 to 5 5          
\putrule from 0 4 to 5 4          %
\putrule from 0 3 to 3 3          %
\putrule from 0 2 to 3 2          %
\putrule from 0 1 to 1 1          %
\putrule from 0 0 to 1 0          %

\putrule from 0 0 to 0 6        %
\putrule from 1 0 to 1 6        %
\putrule from 2 2 to 2 6        %
\putrule from 3 2 to 3 6        
\putrule from 4 4 to 4 6        %
\putrule from 5 4 to 5 6        %
\endpicture
$$
Numbering the rows and columns in the same way as for matrices, let
$$\eqalign{
\lambda_i &= \hbox{the number of boxes in row $i$ of $\lambda$},
\cr
\lambda_j' &= \hbox{the number of boxes in column $j$ of $\lambda$,
\quad and} \cr
|\lambda| &= \hbox{the total number of boxes in $\lambda$.}\cr}
\formula$$
Any partition $\lambda$ can be identified with
the sequence $\lambda=(\lambda_1\ge \lambda_2\ge \ldots )$ and
the {\it conjugate partition} to $\lambda$ is the
partition $\lambda'=(\lambda_1',\lambda_2',\ldots)$.
The {\it hook length} of the box $b$ of $\lambda$ is
$$h(b) = (\lambda_i-i)+(\lambda_j'-j) + 1,
\quad\hbox{if $b$ is in position $(i,j)$ of $\lambda$.}\formula$$
Write $\lambda\vdash n$ if $\lambda$ is a partition with $n$ boxes.
In the example above $\lambda=(553311)$ and $\lambda\vdash 18$.

See [Mac, \S I.7] for details on the representation theory of the symmetric
group.
The irreducible $\CC S_k$-modules
$S_k^\lambda$ are indexed by the elements of
$$\hat S_k = \{ \lambda \vdash n\}\qquad\hbox{and}\qquad
\dim(S_k^\lambda) = {k!\over \displaystyle{\prod_{b\in \lambda} h(b)}}.
\formula
$$
For $\lambda \in \hat S_k$, and $\mu \in \hat
S_{k-1}$,
$$
\Res^{S_k}_{S_{k-1}} (S_k^\lambda)
\cong
\bigoplus_{\lambda/\nu = \square} S_{k-1}^\nu
\qquad\hbox{and}\qquad
\Ind^{S_k}_{S_{k-1}} (S_{k-1}^\mu)
\cong
\bigoplus_{\nu/\mu = \square} S_{k}^\nu. \formula
$$
where the first sum is over all partitions $\nu$ that are obtained
from $\lambda$ by removing a box, and the second sum is over
all partitions $\nu$ which are obtained from $\mu$ by adding a box (this
result
follows, for example, from [Mac, \S I.7 Ex.\ 22(d)]).

The {\it Young lattice} is the graph $\hat S$ given by setting
$$\matrix{
\hbox{vertices on level $k$:}\hskip.4in
\hat S_k = \{ \hbox{partitions $\lambda$ with $k$ boxes}\},
\quad\hbox{and}\cr
\hbox{an edge $\lambda\to\mu$,
$\lambda\in \hat S_k$, $\mu\in \hat S_{k+1}$
if $\mu$ is obtained from $\lambda$ by adding a box}. \hfill \cr
}
\formula$$
It encodes the decompositions in (2.20).  The first few levels
of $\hat S$ are given by
$$
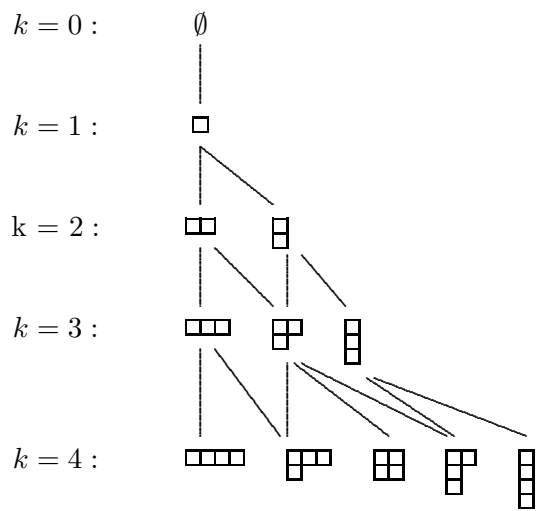
\matrix{
{\beginpicture
\setcoordinatesystem units <0.175cm,0.175cm>         
\setplotarea x from -8 to 27, y from 4 to 31   
\linethickness=0.5pt                          
\put{$k = 4:$} at  -6 -1.5
\putrectangle corners at 3 -2 and 4 -1
\putrectangle corners at 4 -2 and 5 -1
\putrectangle corners at 5 -2 and 6 -1
\putrectangle corners at 6 -2 and 7 -1
\putrectangle corners at 10 -2 and 11 -1
\putrectangle corners at 11 -2 and 12 -1
\putrectangle corners at 12 -2 and 13 -1
\putrectangle corners at 10 -3 and 11 -2
\putrectangle corners at 16 -2 and 17 -1
\putrectangle corners at 17 -2 and 18 -1
\putrectangle corners at 16 -3 and 17 -2
\putrectangle corners at 17 -3 and 18 -2
\putrectangle corners at 21 -2 and 22 -1
\putrectangle corners at 22 -2 and 23 -1
\putrectangle corners at 21 -3 and 22 -2
\putrectangle corners at 21 -4 and 22 -3
\putrectangle corners at 26 -2 and 27 -1
\putrectangle corners at 26 -3 and 27 -2
\putrectangle corners at 26 -4 and 27 -3
\putrectangle corners at 26 -5 and 27 -4
%
\plot 4 6 4 0 /
\plot 5 6 9.5 0 /
\plot 10 5 10 0 /
\plot 10.5 5 17 0 /
\plot 11 5 21 0 /
\plot 15.5 4 21.5 0 /
\plot 16 4 26.5 0 /
\put{$k = 3:$} at  -6 7.5
\putrectangle corners at 3 7 and 4 8
\putrectangle corners at 4 7 and 5 8
\putrectangle corners at 5 7 and 6 8
\putrectangle corners at 9 7 and 10 8
\putrectangle corners at 10 7 and 11 8
\putrectangle corners at 9 6 and 10 7
\putrectangle corners at 14 7 and 15 8
\putrectangle corners at 14 6 and 15 6
\putrectangle corners at 14 5 and 15 7
%
\plot 4 13 4 9  /
\plot 5 13 9 9  /
\plot 10 12.5 10 9  /
\plot 11 12.5 14 9  /
%
%
\put{k = $2:$} at  -6 14.5
\putrectangle corners at 3  14 and 4  15
\putrectangle corners at 4  14 and 5  15
\putrectangle corners at 9  14 and 10  15
\putrectangle corners at 9  13 and 10  14
%
\plot 4 16 4 20  /
\plot 9 16 4 20  /
%
%
\put{$k = 1:$} at  -6 21.5
\putrectangle corners at 3.5  21 and 4.5   22
%
\plot 4 27 4 23  /
%
%
\put{$k = 0:$} at  -6 28.5
\put{$\emptyset$} at  4 28.5
\endpicture} \cr}
$$
For $\mu\in \hat S_k$ define
$$\hat S_k^\mu = \left\{
T=(T^{(0)},T^{(1)},\ldots, T^{(k)})\ \Bigg|\
\matrix{
\hbox{$T^{(0)}=\emptyset\,$,\  $T^{(k)}=\mu\,$, \quad and,
\quad for each $\ell$,}\hfill\cr
\hbox{$T^{(\ell)}\in \hat S_\ell$ and
$T^{(\ell)}\to T^{(\ell+1)}$ is an edge in $\hat S$} \cr}
\right\}$$
so that $\hat S_k^\mu$ is the set of paths from
$\emptyset\in \hat S_0$ to $\mu\in \hat S_k$ in
the graph $\hat S$.  In terms of the Young lattice
$$\dim(S_k^\mu) = \Card(\hat S_k^\mu).
\formula$$
This is a translation of the classical statement (see [Mac, \S I.7.6(ii)])
that $\dim(S_k^\mu)$ is the number of standard Young tableaux
of shape $\lambda$ (the correspondence is obtained by putting
the entry $\ell$ in the box of $\lambda$ which is added at the
$\ell$th step $T^{(\ell-1)}\to T^{(\ell)}$ of the path).

\bigskip\noindent
{\it Structure of the algebra $\CC A_k(n)$}
\medskip

Build a graph $\hat A$ by setting
$$\matrix{
\hbox{vertices on level $k$:}\hskip.4in
\hat A_k = \{ \hbox{partitions $\mu$}\ |\
k-|\mu| \in \ZZ_{\ge 0}\}, \hfill \cr
\hbox{vertices on level $k+{1\over2}$:}\quad
\hat A_{k+{1\over2}} =
\hat A_k = \{ \hbox{partitions $\mu$}\ |\
k-|\mu| \in \ZZ_{\ge 0}\}, \hfill \cr
\hbox{an edge $\lambda\to\mu$,
$\lambda\in \hat A_k$, $\mu\in \hat A_{k+{1\over2}}$
if $\lambda=\mu$ or if $\mu$ is obtained from
$\lambda$ by removing a box}, \hfill \cr
\hbox{an edge $\mu\to \lambda$,
$\mu\in \hat A_{k+{1\over2}}$,
$\lambda\in \hat A_{k+1}$,
if $\lambda=\mu$ or if $\lambda$ is obtained from
$\mu$ by adding a box}. \cr
}
\formula$$
The first few levels of $\hat A$ are given by
$$
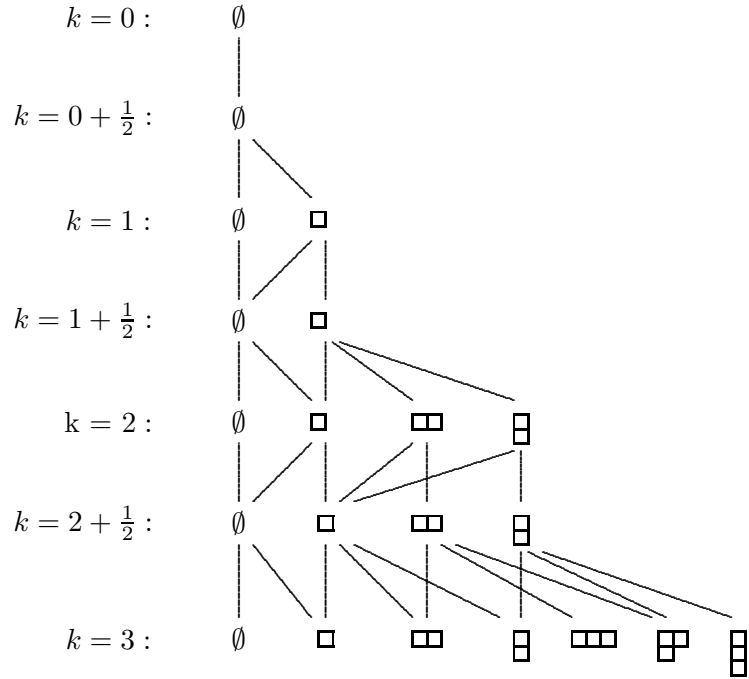
\matrix{
{\beginpicture
\setcoordinatesystem units <0.175cm,0.175cm>         
\setplotarea x from -8 to 40, y from 4 to 31   
\linethickness=0.5pt                          
\put{$k = 3:$} at  -6 -.5
\put{$\emptyset$} at 3 -.5
\putrectangle corners at 8.5 -1 and 9.5 0
\putrectangle corners at 15 -1 and 16 0
\putrectangle corners at 16 -1 and 17 0
\putrectangle corners at 22 -1 and 23 0
\putrectangle corners at 22 -2 and 23 -1
\putrectangle corners at 27 -1 and 28 0
\putrectangle corners at 28 -1 and 29 0
\putrectangle corners at 26 -1 and 27 0
\putrectangle corners at 32 -1 and 33 0
\putrectangle corners at 33 -1 and 34 0
\putrectangle corners at 32 -2 and 33 -1
\putrectangle corners at 37 -1 and 38 0
\putrectangle corners at 37 -2 and 38 -1
\putrectangle corners at 37 -3 and 38 -2
%
\plot 3 1 3 6  /
\plot 8 1 4 6  /
\plot 9 1 9 6  /
\plot 15 1 10 6  /
\plot 16 1 16 6  /
\plot 26 1 17 6  /
\plot 31.5 1 18 6  /
\plot 21 1 11 6  /
\plot 22.5 1 22.5 5.5  /
\plot 32.5 1 23 5.5  /
\plot 37 1 24 5.5  /
%
\put{$k = 2+{1\over2}:$} at  -7.75 7.5
\put{$\emptyset$} at 3 7.5
\putrectangle corners at 8.5 7  and 9.5 8
\putrectangle corners at 15 7 and 16 8
\putrectangle corners at 16 7 and 17 8
\putrectangle corners at 22 7 and 23 8
\putrectangle corners at 22 6 and 23 7
%
\plot 3 9 3 13  /
\plot 8 13 4 9  /
\plot 9 13 9 9  /
\plot 15 13 10 9  /
\plot 16 13 16 9  /
\plot 22 12.5 11 9  /
\plot 22.5 12.5 22.5 9  /
%
%
\put{k = $2:$} at  -6 14.5
\put{$\emptyset$} at 3 14.5
\putrectangle corners at 8  14 and 9   15
\putrectangle corners at 15  14 and 16  15
\putrectangle corners at 16  14 and 17  15
\putrectangle corners at 22  14 and 23  15
\putrectangle corners at 22  13 and 23  14
%
\plot 3 16 3 20  /
\plot 8 16 4 20  /
\plot 9 16 9 20  /
\plot 15 16 9.5 20  /
\plot 22 16 10 20  /
%
%
\put{$k = 1+{1\over2}:$} at  -7.75 21.5
\put{$\emptyset$} at 3 21.5
\putrectangle corners at 8  21 and 9   22
%
\plot 3 23 3 27  /
\plot 8 27 4 23  /
\plot 9 27 9 23  /
%
%
\put{$k = 1:$} at  -6 28.5
\put{$\emptyset$} at  3 28.5
\putrectangle corners at 8  28 and 9  29
%
\put{$k = 0+{1\over2}:$} at  -7.75 35.5
\put{$\emptyset$} at  3 35.5
%
%
\plot 3 30 3 34  /
\plot 8 30 4 34  /
%
%
\put{$k = 0:$} at  -6 42.5
\put{$\emptyset$} at  3 42.5
%
%
\plot 3 37 3 41  /
\endpicture} \cr}
$$

The following result is an immediate consequence of the
Tits deformation theorem, Theorems 5.12 and 5.13 in this paper
(see also [CR (68.17)]).

\thm
\item{(a)}  For all but a finite number of $n \in \CC$ the
algebra $\CC A_k(n)$ is semisimple.
\item{(b)}  If $\CC A_k(n)$ is semisimple
then the irreducible $\CC A_k(n)$-modules,
$$\hbox{$A_k^\mu$\quad are indexed by elements of the set}\quad
\hat A_k = \{\ \hbox{partitions $\mu$}\ |\ k-|\mu| \in \ZZ_{\ge 0} \},
\qquad\hbox{and}
$$
$
\dim(A_k^\mu) = \hbox{(number of paths from $\emptyset\in \hat A_0$
to $\mu\in \hat A_k$ in the graph $\hat A$)}.
$
\endthm

Let
$$\hat A_k^\mu = \left\{
T=(T^{(0)},T^{({1\over2})},\ldots, T^{(k-{1\over2})},T^{(k)})\ \Bigg|\
\matrix{
\hbox{$T^{(0)}=\emptyset\,$,\  $T^{(k)}=\mu\,$, \quad and,
\quad for each $\ell$,}\hfill\cr
\hbox{$T^{(\ell)}\in \hat A_\ell$ and
$T^{(\ell)}\to T^{(\ell+{1\over2})}$ is an edge in $\hat A$} \cr}
\right\}$$
so that $\hat A_k^\mu$ is the set of paths from
$\emptyset\in \hat A_0$ to $\mu\in \hat A_k$ in
the graph $\hat A$.
If $\mu\in \hat S_k$ then $\mu\in \hat A_k$ and $\mu\in \hat
A_{k+{1\over2}}$
and, for notational convenience in the following theorem,
$$\eqalign{
\hbox{identify}\qquad &P=(P^{(0)},P^{(1)},\ldots,P^{(k)}) \in \hat S_k^\mu
\qquad \hbox{with the corresponding} \cr
&P=(P^{(0)},P^{(0)},P^{(1)},P^{(1)},\ldots, P^{(k-1)},P^{(k-1)},P^{(k)})
\qquad\ \in \hat A_k^\mu, \qquad\hbox{and} \cr
&P=(P^{(0)},P^{(0)},P^{(1)},P^{(1)},\ldots,
P^{(k-1)},P^{(k-1)}, P^{(k)},P^{(k)})\in
\hat A_{k+{1\over2}}^\mu. \cr
}$$

For $\ell\in {1\over2}\ZZ_{\ge 0}$
and $n\in \CC$ such that $\CC A_\ell(n)$ is semisimple let
$\chi_{A_\ell(n)}^\mu$, $\mu\in \hat A_\ell$, be the irreducible
characters of $\CC A_\ell(n)$. Let
$\tr_\ell\colon \CC A_\ell(n)\to \CC$ be the traces on
$\CC A_\ell(n)$ defined in (2.8) and
define constants $\tr_\ell^\mu(n)$, $\mu\in \hat A_\ell$, by
$$\tr_\ell = \sum_{\mu\in \hat A_\ell} \tr_\ell^\mu(n)\chi_{A_\ell(n)}^\mu.
\formula
$$

\thm
\item{(a)}
Let $n\in \CC$ and let $k\in \hbox{$1\over2$}\ZZ_{\ge 0}$.
Assume that
$$\tr_\ell^\lambda(n)\ne 0,
\qquad\hbox{
for all $\lambda\in \hat A_\ell$,
$\ell\in {1\over2}\ZZ_{\ge 0}$,
$\ell < k$.}$$
Then the partition algebras
$$\hbox{
$\CC A_{\ell}(n)$ are semisimple
for all $\ell\in \hbox{$1\over2$}\ZZ_{\ge 0}, \ell\le k$}.
\formula
$$
For each $\ell\in {1\over2}\ZZ_{\ge 0}$, $\ell\le k-{1\over2}$,
define
$$\varepsilon_\mu^\lambda
= {\tr_{\ell-{1\over2}}^\lambda(n)\over \tr_{\ell-1}^\mu(n)}
\qquad\hbox{for each edge $\mu\to \lambda$,
$\mu\in \hat A_{\ell-1}$, $\lambda\in \hat A_{\ell-{1\over2}}$,
in the graph $\hat A$.}
$$
Inductively define elements in $\CC A_\ell(n)$ by
$$e_{PQ}^\mu
= {1\over \sqrt{\varepsilon_\mu^\tau \varepsilon_\mu^\gamma}}\,
e^\tau_{P^-T}p_{\ell} e^\gamma_{TQ^-},
\qquad\hbox{for $\mu\in \hat A_\ell$,
$|\mu|\le \ell-1$, $P,Q\in \hat A_{\ell}^\mu$,}
\formula
$$
where $\tau = P^{(\ell-{1\over2})}$, $\gamma = Q^{(\ell-{1\over2})}$,
$R^-=(R^{(0)},\ldots, R^{(\ell-{1\over2})})$ for
$R=(R^{(0)},\ldots, R^{(\ell-{1\over2})}, R^{(\ell)})\in \hat A_\ell^\mu$
and $T$ is an element of $\hat A_{\ell-1}^\mu$ (the element
$e_{PQ}^\lambda$ does not depend on the choice of $T$).
Then define
$$
e^\lambda_{PQ} = (1-z)s^\lambda_{PQ},
\quad\hbox{for}\quad
\lambda\in \hat S_\ell, P,Q\in \hat S_\ell^\lambda,
\qquad\hbox{where}\quad
z = \sum_{\mu\in \hat A_{\ell}\atop |\mu|\le \ell-1}
\sum_{P\in \hat A_\ell^\mu} e^\mu_{PP}
\formula$$
and
$\{ s^\lambda_{PQ}\ |\ \lambda\in \hat S_\ell, P,Q\in \hat S_\ell^\lambda\}
$ is any set of matrix units for the
the group algebra of the symmetric group $\CC S_\ell$.
Together, the elements in (2.28) and (2.29)
form a set of matrix units in $\CC A_\ell(n)$.
\item{(b)}
Let $n\in \ZZ_{\ge 0}$ and let $k\in \hbox{$1\over2$}\ZZ_{>0}$
be minimal such that $\tr_k^\lambda(n)=0$ for some
$\lambda\in \hat A_k$.  Then
$\CC A_{k+{1\over2}}(n)$ is not semisimple.
\item{(c)}
Let $n\in \ZZ_{\ge 0}$ and $k\in \hbox{$1\over2$}\ZZ_{>0}$.
If $\CC A_k(n)$ is not semisimple then $\CC A_{k+j}(n)$ is
not semisimple for $j\in \ZZ_{>0}$.

\pf
(a) Assume that $\CC A_{\ell-1}(n)$ and $\CC A_{\ell-{1\over2}}(n)$
are both semisimple and that
$\tr_{\ell-1}^\mu(n)\ne 0$ for all $\mu\in \hat A_{\ell-1}$.
If $\lambda\in \hat A_{\ell-{1\over2}}$ then
$\varepsilon_\mu^\lambda \ne 0$
if and only if
$\tr_{\ell-{1\over2}}^\lambda(n)\ne 0$, and,
since the ideal $\CC I_\ell(n)$ is isomorphic to the
basic construction
$\CC A_{\ell-{1\over2}}(n)
\otimes_{\CC A_{\ell-1}(n)} \CC A_{\ell-{1\over2}}(n)$
(see (2.13))
it then follows from Theorem 4.28 that
$\CC I_\ell(n)$ is semisimple if and only if
$\tr_{\ell-{1\over2}}^\lambda(n)\ne 0$
for all $\lambda\in\hat A_{\ell-{1\over2}}$.
Thus, by (2.12), {\it if $\CC A_{\ell-1}(n)$ and $\CC A_{\ell-{1\over2}}(n)$
are both semisimple and
$\tr_{\ell-1}^\mu(n)\ne 0$ for all $\mu\in \hat A_{\ell-1}$} then
$$\CC A_\ell(n)\qquad\hbox{is semisimple if and only if}\qquad
\hbox{$\tr_{\ell-{1\over2}}^\lambda(n)\ne 0$
for all $\lambda\in\hat A_{\ell-{1\over2}}$.}
\formula
$$
By Theorem 4.28, when
$\tr_{\ell-{1\over2}}^\lambda(n)\ne 0$ for all
$\lambda\in\hat A_{\ell-{1\over2}}$,
the algebra $\CC I_\ell(n)$ has matrix units given by
the formulas in (2.28).
The element $z$ in (2.29)
is the central idempotent in $\CC A_{\ell}(n)$ such that
$\CC I_{\ell}(n)=z\CC A_{\ell}(n)$.
Hence the complete set of elements in (2.28) and (2.29)
form a set of matrix units for $\CC A_{\ell}(n)$.
This completes the proof of (a) and (b) follows from Theorem 4.28(b).

\smallskip\noindent
(c) Part (g) of Theorem 4.28 shows that
if $\CC A_{\ell-1}(n)$ is not semisimple then
$\CC A_{\ell}(n)$ is not semisimple.
\endpf

\medskip\noindent
{\it Specht modules}

Let $A$ be an algebra. An idempotent is a nonzero element $p\in A$ such
that $p^2=p$. A {\it minimal idempotent} is an idempotent
$p$ which cannot be written as a sum $p = p_1+p_2$ with $p_1p_2=p_2p_1=0$.
If $p$ is an idempotent in $A$ and $pAp = \CC p$ then
$p$ is a minimal idempotent of $A$ since, if
$p=p_1+p_2$ with $p_1^2=p_1$, $p_2^2=p_2$ and $p_1p_2=p_2p_1=0$ then
$pp_1p=kp$ for some constant $p$ and so
$kp_1 = kpp_1 = pp_1pp_1 = p_1$ giving that either $p_1=0$ or
$k=1$, in which case $p_1=pp_1p=p$.

Let $p$ be an idempotent in $A$.  Then the map
$$\matrix{ (pAp)^{\rm op}&\mapright{\sim}& \End_A(Ap) \cr
pbp &\mapsto &\phi_{pbp} \cr},
\qquad\hbox{where}\qquad
\phi_{pbp}(ap) = (ap)(pbp)=apbp,
\quad\hbox{for $ap\in Ap$,}\formula$$
is a ring isomorphism.

If $p$ is a minimal idempotent of $A$ and $Ap$ is a semisimple
$A$-module then $Ap$ must be a simple $A$-module.
To see this suppose that $Ap$ is not simple
so that there are $A$-submodules $V_1$ and $V_2$ of $Ap$ such that
$Ap = V_1\oplus V_2$.
Let $\phi_1, \phi_2\in \End_A(Ap)$ be the $A$-invariant projections
on $V_1$ and $V_2$.  By (2.31) $\phi_1$ and $\phi_2$
are given by right multiplication by $p_1=p\tilde p_1 p$
and $p_2 = p\tilde p_2p$, respectively, and it follows that
$p=p_1+p_2$, $V_1=Ap_1$, $V_2=Ap_2$, and $Ap = Ap_1\oplus Ap_2$.
Then
$p_1^2 = \phi_1(p_1) = \phi_1^2(p) = p_1$ and
$p_1p_2 =\phi_2(p_1) = \phi_2(\phi_1(p))=0$.
Similarly $p_2^2=p_2$ and $p_2p_1=0$.  Thus $p$ is not a minimal idempotent.

If $p$ is an idempotent in $A$ and $Ap$ is a simple $A$-module
then
$$pAp=\End_A(Ap)^{\rm op} = \CC (p\cdot 1\cdot p) = \CC p,$$
by (2.31) and Schur's lemma (Theorem 5.3).

The group algebra of the symmetric group $S_k$ over the ring $\ZZ$ is
$$S_{k,\ZZ} = \ZZ S_k
\qquad\hbox{and}\qquad
\CC S_k = \CC \otimes_\ZZ S_{k,\ZZ}\,,
\formula$$
where the tensor product is defined via the inclusion
$\ZZ\hookrightarrow \CC$.
Let $\lambda=(\lambda_1,\lambda_2,\ldots, \lambda_\ell)$
be a partition of $k$.  Define subgroups of $S_k$ by
$$S_\lambda = S_{\lambda_1}\times\cdots \times S_{\lambda_\ell}
\qquad\hbox{and}\qquad
S_{\lambda'} = S_{\lambda'_1}\times\cdots \times S_{\lambda'_r},
\formula$$
where  $\lambda'=(\lambda'_1,\lambda'_2,\ldots, \lambda'_r)$
is the conjugate partition to $\lambda$, and let
$${\bf 1}_\lambda = \sum_{w\in S_\lambda} w
\qquad\hbox{and}\qquad
\varepsilon_{\lambda'}
=\sum_{w\in S_{\lambda'}} (-1)^{\ell(w)} w.
\formula$$
Let $\tau$ be the permutation in $S_k$ that takes
the row reading tableau of shape $\lambda$ to the column reading
tableau of shape $\lambda$.  For example for $\lambda = (553311)$,
$$\tau =(2,7,8,12,9,16,14,4,15,10,18,6)(3,11)(5,17),
\quad\hbox{since}\quad\tau \cdot
{\beginpicture
\setcoordinatesystem units <0.42cm,0.4cm>
\setplotarea x from 0 to 5, y from -1 to 2
\linethickness=0.5pt
\put{$1$} at 0.5 2.5
\put{$2$} at 1.5 2.5
\put{$3$} at 2.5 2.5
\put{$4$} at 3.5 2.5
\put{$5$} at 4.5 2.5
\put{$6$} at 0.5 1.5
\put{$7$} at 1.5 1.5
\put{$8$} at 2.5 1.5
\put{$9$} at 3.5 1.5
\put{$10$} at 4.5 1.5
\put{$11$} at 0.5 0.5
\put{$12$} at 1.5 0.5
\put{$13$} at 2.5 0.5
\put{$14$} at 0.5 -0.5
\put{$15$} at 1.5 -0.5
\put{$16$} at 2.5 -0.5
\put{$17$} at 0.5 -1.5
\put{$18$} at 0.5 -2.5
\putrectangle corners at 0 -3 and 1 -2
\putrectangle corners at 0 -2 and 1 -1
\putrectangle corners at 0 -1 and 3 0
\putrectangle corners at 0 0 and 3 1
\putrectangle corners at 0 1 and 5 2
\putrectangle corners at 0 2 and 5 3
\putrectangle corners at 0 -3 and  1 3
\putrectangle corners at 1 -1 and  2 3
\putrectangle corners at 2 -1 and  3 3
\putrectangle corners at 3 1 and  4 3
\endpicture} =
{\beginpicture
\setcoordinatesystem units <0.42cm,0.4cm>
\setplotarea x from 0 to 5, y from -1 to 2
\linethickness=0.5pt
\put{$1$} at 0.5 2.5
\put{$7$} at 1.5 2.5
\put{$11$} at 2.5 2.5
\put{$15$} at 3.5 2.5
\put{$17$} at 4.5 2.5
\put{$2$} at 0.5 1.5
\put{$8$} at 1.5 1.5
\put{$12$} at 2.5 1.5
\put{$16$} at 3.5 1.5
\put{$18$} at 4.5 1.5
\put{$3$} at 0.5 0.5
\put{$9$} at 1.5 0.5
\put{$13$} at 2.5 0.5
\put{$4$} at 0.5 -0.5
\put{$10$} at 1.5 -0.5
\put{$14$} at 2.5 -0.5
\put{$5$} at 0.5 -1.5
\put{$6$} at 0.5 -2.5
\putrectangle corners at 0 -3 and 1 -2
\putrectangle corners at 0 -2 and 1 -1
\putrectangle corners at 0 -1 and 3 0
\putrectangle corners at 0 0 and 3 1
\putrectangle corners at 0 1 and 5 2
\putrectangle corners at 0 2 and 5 3
\putrectangle corners at 0 -3 and  1 3
\putrectangle corners at 1 -1 and  2 3
\putrectangle corners at 2 -1 and  3 3
\putrectangle corners at 3 1 and  4 3
\endpicture}.$$
The {\it Specht module} for $S_k$ is the $\ZZ S_k$-module
$$S_{k,\ZZ}^\lambda = \im~\Psi_{S_k} = (\ZZ S_k)\,p_\lambda,
\qquad\hbox{where $p_\lambda= {\bf 1}_\lambda
\tau \varepsilon_{\lambda'}\tau^{-1},$ and}
\formula$$
where $\Psi_{S_k}$ is the $\ZZ S_k$-module homomorphism given by
$$\matrix{
\Psi_{S_k}\colon
&(\ZZ S_k) {\bf 1}_\lambda &\mapright{\iota} &\ZZ S_k &\mapright{\pi}
&(\ZZ S_k) \tau \varepsilon_{\lambda'}\tau^{-1} \cr
\cr
&b\,{\bf 1}_\lambda
&\longmapsto &b\,{\bf 1}_\lambda
&\longmapsto &b\,{\bf 1}_\lambda \tau \varepsilon_{\lambda'}\tau^{-1}\cr
}
\formula$$
By induction and restriction rules for the representations of the
symmetric groups, the $\CC S_k$-modules $(\CC S_k){\bf 1}_\lambda$
and $(\CC S_k)\tau \varepsilon_{\lambda'}\tau^{-1}$ have only
one irreducible component in common and it follows (see [Mac,\S I.7, Ex.\ 
15])
that
$$S_k^\lambda = \CC\otimes_{\ZZ} S_{k,\ZZ}^\lambda
\qquad\hbox{is the irreducible $\CC S_k$-module indexed by $\lambda$,}
\formula$$
once one shows that $\Psi_{S_k}$ is not the zero map.

Let  $k \in {1\over 2} \ZZ_{>0}$. For an indeterminate $x$,
define the $\ZZ[x]$-algebra by
$$A_{k,\ZZ} = \ZZ[x]\hbox{-span}\{ d\in A_k\}
\formula$$
with multiplication given by replacing $n$ with $x$ in (2.1).
For each $n\in \CC$,
$$\CC A_k(n) = \CC\otimes_{\ZZ[x]} A_{k,\ZZ},
\quad\hbox{where the $\ZZ$-module homomorphism}\quad
\matrix{
{\rm ev}_n\colon &\ZZ[x] &\to &\CC \cr
&x &\mapsto &n \cr}
\formula$$
is used to define the tensor product.
Let $\lambda$ be a partition with $\le k$ boxes.
Let $b\otimes p_k^{\otimes(k-|\lambda|)}$ denote the
image of $b\in A_{|\lambda|,\ZZ}$ under the map given by
$$
\matrix{
\matrix{
A_{|\lambda|,\ZZ} &\longrightarrow &A_{k,\ZZ} \cr
\cr
{\beginpicture
\setcoordinatesystem units <0.35cm,0.15cm>
\setplotarea x from 1 to 5, y from -1 to 2
\linethickness=0.5pt
\put{$b$} at 3 .5
\put{$\bullet$} at 1 -1 \put{$\bullet$} at 1 2
\put{$\bullet$} at 5 -1 \put{$\bullet$} at 5 2
\setdashes <.04cm>
\plot 1 2 5 2  5 -1 1 -1 1 2 /
\endpicture} &\longmapsto &
{\beginpicture
\setcoordinatesystem units <0.35cm,0.15cm>
\setplotarea x from 1 to 3, y from -1 to 2
\linethickness=0.5pt
\put{$b$} at 3 .5
\put{$\bullet$} at 1 -1 \put{$\bullet$} at 1 2
\put{$\bullet$} at 5 -1 \put{$\bullet$} at 5 2
\setdashes <.04cm>
\plot 1 2 5 2  5 -1 1 -1 1 2 /
\endpicture}
\underbrace{\beginpicture
\setcoordinatesystem units <0.35cm,0.15cm>
\setplotarea x from 1 to 3, y from -1 to 2
\linethickness=0.5pt
\put{$\bullet$} at 1 -1 \put{$\bullet$} at 1 2
\put{$\bullet$} at 2 -1 \put{$\bullet$} at 2 2
\put{$\cdots$} at 3 -1 \put{$\cdots$} at 3 2
\put{$\bullet$} at 4 -1 \put{$\bullet$} at 4 2
\endpicture}_{k - |\lambda|} \,, \cr}
\hfill &
\qquad\hbox{if $k$ is an integer, and}
\hfill\cr
\matrix{
A_{|\lambda|+{1\over2},\ZZ} &\longrightarrow &A_{k,\ZZ} \cr
\cr
{\beginpicture
\setcoordinatesystem units <0.35cm,0.15cm>
\setplotarea x from 1 to 4, y from -1 to 2
\linethickness=0.5pt
\put{$b$} at 3 .5
\put{$\bullet$} at 1 -1 \put{$\bullet$} at 1 2
\put{$\bullet$} at 5 -1 \put{$\bullet$} at 5 2
\setquadratic
\plot 5 2 5.25 .5 5 -1 /
\setlinear
\setdashes <.04cm>
\plot 1 2 5 2  5 -1 1 -1 1 2 /
\endpicture} &\longmapsto &
{\beginpicture
\setcoordinatesystem units <0.35cm,0.15cm>
\setplotarea x from 1 to 4, y from -1 to 2
\linethickness=0.5pt
\put{$b$} at 3 .5
\put{$\bullet$} at 1 -1 \put{$\bullet$} at 1 2
\put{$\bullet$} at 5 -1 \put{$\bullet$} at 5 2
\setquadratic
\plot 5 2 5.25 .5 5 -1 /
\setlinear
\setdashes <.04cm>
\plot 1 2 5 2  5 -1 1 -1 1 2 /
\endpicture}
\underbrace{\beginpicture
\setcoordinatesystem units <0.35cm,0.15cm>
\setplotarea x from .5 to 3, y from -1 to 2
\linethickness=0.5pt
\put{$\bullet$} at 1 -1 \put{$\bullet$} at 1 2
\put{$\bullet$} at 2 -1 \put{$\bullet$} at 2 2
\put{$\cdots$} at 3.5 -1 \put{$\cdots$} at 3.5 2
\put{$\bullet$} at 5 -1 \put{$\bullet$} at 5 2
\put{$\bullet$} at 6 -1 \put{$\bullet$} at 6 2
\plot 0 -1 2.7 -1 /
\plot 4.3 -1 6 -1 /
\plot 0 2 2.7 2 /
\plot 4.3 2 6 2 /
\endpicture}_{k - |\lambda|-{1\over2}} \,, \cr}
\hfill &
\qquad\hbox{if $k-{1\over2}$ is an integer.}
}
$$

For $k \in {1\over2} \ZZ_{>0}$, define an $A_{k,\ZZ}$-module homomorphism
$$\matrix{
\Psi_{A_k}\colon
&A_{k,\ZZ} t_\lambda
&\mapright{\psi_1}
&A_{k,\ZZ} s_{\lambda'}
&\mapright{\psi_2} &A_{k,\ZZ}/I_{|\lambda|,\ZZ} \cr
\cr
&bt_\lambda
&\longmapsto &bt_\lambda s_{\lambda'}
&\longmapsto &\overline{bt_\lambda s_{\lambda'}}\,,  \cr
}
\formula$$
where $I_{|\lambda|,\ZZ}$ is the ideal
$$I_{|\lambda|,\ZZ} = \ZZ[x]\hbox{-span}\,
\big\{\ d\in A_k\ |\ \hbox{$d$ has propagating number $<|\lambda|$}\
\big\}$$
and $t_\lambda,s_{\lambda'} \in A_{k,\ZZ}$ are defined by
$$
\matrix{
t_\lambda = {\bf 1}_\lambda\otimes p_k^{\otimes(k-|\lambda|)}
&\hbox{and}&
s_{\lambda'}
= \tau\varepsilon_{\lambda'}\tau^{-1}\otimes p_k^{\otimes(k-|\lambda|)}.}
\formula
$$
The {\it Specht module} for $\CC A_k(n)$ is the $A_{k,\ZZ}$-module
$$A_{k,\ZZ}^\lambda = \im~\Psi_{A_k}
= \big(\hbox{image of $A_{k,\ZZ}e_\lambda$ in
$A_{k,\ZZ}/I_{|\lambda|,\ZZ}$}\big),
\qquad\hbox{where}\qquad
e_\lambda = p_\lambda\otimes p_k^{\otimes(k-|\lambda|)}.
\formula$$

\prop Let $k \in {1\over 2}\ZZ_{>0}$, and let $\lambda$ be a partition with
$\le k$ boxes. If $n\in \CC$ such that $\CC A_k(n)$ is semisimple, then
$$
A_k^\lambda(n) = \CC\otimes_{\ZZ[x]} A_{k,\ZZ}^\lambda
\qquad\hbox{is the irreducible $\CC A_k(n)$-module indexed by $\lambda$,}
$$
where the tensor product is defined via the $\ZZ$-module homomorphism in
(2.39).

\pf
Let $r=|\lambda|$.  Since
$$\CC A_r(n)/\CC I_r(n) \cong \CC S_r$$
and $p_\lambda$ is a minimal idempotent of $\CC S_r$, it follows
from (4.20) that $\overline{e_\lambda}$, the image of
$e_\lambda$ in $(\CC A_k(n))/(\CC I_r(n))$,
is a minimal idempotent in $(\CC A_k(n))/(\CC I_r(n))$.
Thus
$$\left({\CC A_k(n)\over \CC I_r(n)}\right)\overline e_\lambda
\qquad\hbox{is a simple $(\CC A_k(n))/(\CC I_r(n))$-module.}
$$
Since the projection $\CC A_k(n)\to (\CC A_k(n))/(\CC I_r(n))$
is surjective, any simple
$(\CC A_k(n))/(\CC I_r(n))$-module is a simple $\CC A_k(n)$-module.
\endpf

\section 3. Schur-Weyl Duality for Partition Algebras

\bigskip
Let $n \in \ZZ_{>0}$ and let $V$ be a vector space with basis
$v_1, \ldots, v_n$. Then the tensor product
$$
V^{\otimes k} = \underbrace{V \otimes V \otimes \cdots \otimes V}_{k~{\rm
factors}}
\quad
\hbox{ has basis }
\quad
\left\{\ v_{i_1} \otimes \cdots \otimes v_{i_k} \ | \ 1 \le i_1, \ldots,
i_k \le n\ \right\}.
$$
For $d \in A_k$ and values $i_1, \ldots, i_k, i_{1'}, \ldots, i_{k'} \in
\{1, \ldots, n\}$
define
$$
(d)^{i_1, \ldots, i_k}_{i_{1'}, \ldots, i_{k'}} = \cases{
1, & if $i_r = i_s$ when $r$ and $s$ are in the same block of $d$, \cr
0, & otherwise.} \formula
$$
For example, viewing $(d)^{i_1, \ldots, i_k}_{i_{1'}, \ldots, i_{k'}}$ as
the diagram
$d$ with vertices labeled by the values $i_1, \ldots, i_k$ and  $i_{1'},
\ldots, i_{k'}$,
we have
$$
{\beginpicture
\setcoordinatesystem units <0.5cm,0.2cm> 
\setplotarea x from 0 to 7, y from 0 to 3    
\linethickness=0.5pt                        
\put{$i_1$} at 0 3.5
\put{$i_2$} at 1 3.5
\put{$i_3$} at 2 3.5
\put{$i_4$} at 3 3.5
\put{$i_5$} at 4 3.5
\put{$i_6$} at 5 3.5
\put{$i_7$} at 6 3.5
\put{$i_8$} at 7 3.5
\put{$i_{1'}$} at  0 -2.75
\put{$i_{2'}$} at 1 -2.75
\put{$i_{3'}$} at 2 -2.75
\put{$i_{4'}$} at 3 -2.75
\put{$i_{5'}$} at 4 -2.75
\put{$i_{6'}$} at 5 -2.75
\put{$i_{7'}$} at 6 -2.75
\put{$i_{8'}$} at 7 -2.75
\put{$\bullet$} at 0 -1 \put{$\bullet$} at 0 2
\put{$\bullet$} at 1 -1 \put{$\bullet$} at 1 2
\put{$\bullet$} at 2 -1 \put{$\bullet$} at 2 2
\put{$\bullet$} at 3 -1 \put{$\bullet$} at 3 2
\put{$\bullet$} at 4 -1 \put{$\bullet$} at 4 2
\put{$\bullet$} at 5 -1 \put{$\bullet$} at 5 2
\put{$\bullet$} at 6 -1 \put{$\bullet$} at 6 2
\put{$\bullet$} at 7 -1 \put{$\bullet$} at 7 2
\plot 1 -1 1 2  /
\plot 5 -1 5 2  /
\plot 7 -1 7 2  /
\setquadratic
\plot 0  2  .5 1.25 1 2 /
\plot 4  2 4.5 1.25 5 2 /
\plot 5  2 5.5 1.25 6 2 /
\plot 0  2  .5 .25 1 -1 /
\plot 1 -1 2.5 0 4 -1 /
\plot 1  2  2  1.25 3 2 /
\plot 2 -1 2.5 0.5 3 -1 /
\plot 3 -1 4 1 5 -1 /
\plot 5 -1 5.5 .5 6 -1 /
\endpicture} =
\delta_{i_{1}i_{2}} \delta_{i_{1}i_{4}}
\delta_{i_{1}i_{2'}} \delta_{i_{1}i_{5'}}
\delta_{i_{5}i_{6}} \delta_{i_{5}i_{7}}
\delta_{i_{5}i_{3'}}\delta_{i_{5}i_{4'}}
\delta_{i_{5}i_{6'}}\delta_{i_{5}i_{7'}}
\delta_{i_{8}i_{8'}}.
$$
With this notation, the formula
$$
d(v_{i_1} \otimes \cdots \otimes v_{i_k})
= \sum_{1 \le i_{1'}, \ldots i_{k'} \le n}
(d)^{i_1, \ldots, i_k}_{i_{1'}, \ldots, i_{k'}}
v_{i_{1'}} \otimes \cdots\otimes v_{i_{k'}} \formula
$$
defines actions
$$
\Phi_k: \CC A_k \longrightarrow \End(V^{\otimes k})
\qquad\hbox{and}\qquad
\Phi_{k + {1 \over 2}}: \CC \Akp \longrightarrow \End(V^{\otimes k})
\formula
$$
of $\CC A_k$  and $\CC \Akp$ on $V^{\otimes k}$, where the
second map $\Phi_{k+{1\over2}}$ comes from the fact that
if $d \in \Akp$, then $d$ acts on the subspace
$$
V^{\otimes k} \cong
V^{\otimes k} \otimes v_n
= \CC\hbox{-span}\left\{ v_{i_1} \otimes \cdots \otimes v_{i_k} \otimes v_n
\ | \ 1 \le i_1, \ldots, i_k \le n \right\}
\quad\subseteq V^{\otimes(k+1)}. \formula
$$
In other words, the map $\Phi_{k+{1\over2}}$ is
obtained from $\Phi_{k+1}$ by restricting to the subspace
$V^{\otimes k} \otimes v_n$ and identifying $V^{\otimes k}$
with $V^{\otimes k} \otimes v_n$.

The group $GL_n(\CC)$ acts on the vector spaces $V$ and
$V^{\otimes k}$ by
$$gv_i= \sum_{j=1}^n g_{ji}v_j,
\qquad\hbox{and}\qquad
g (v_{i_1} \otimes v_{i_2} \otimes \cdots\otimes v_{i_k})
= g v_{i_1} \otimes g v_{i_2} \otimes \cdots \otimes g v_{i_k},
\formula $$
for $g = (g_{ij})\in GL_n(\CC)$.
View $S_n \subseteq GL_n(\CC)$ as the subgoup of permutation matrices
and let
$$\End_{S_n}( V^{\otimes k} )
= \left\{ b \in \End(V^{\otimes k})\ | \ b \sigma v = \sigma b v
\hbox{ for all } \sigma \in S_n \hbox{ and } v \in V^{\otimes k} \right\}.
$$

\thm Let $n \in \ZZ_{>0}$ and let $\{x_d\ |\ d\in A_k\}$ be the basis
of $\CC A_k(n)$ defined in (2.3).  Then
\smallskip
\item{(a)}
$\Phi_k: \CC A_k(n) \to \End(V^{\otimes k})$ has
$$
\hbox{\rm im}~\Phi_k
= \End_{S_n}(V^{\otimes k}) \qquad\hbox{and}\qquad
\ker \Phi_k
= \CC\hbox{-span}\{ x_d\ |\ \hbox{$d$ has more than $n$ blocks} \},
\qquad\hbox{and}
$$
\item{(b)}
$\Phi_{k + {1 \over 2}}: \CC A_{k + {1 \over 2}}(n) \to \End(V^{\otimes
k})$ has
$$
\hbox{\rm im}~\Phi_{k + {1 \over 2}}
= \End_{S_{n-1}}(V^{\otimes k}) \qquad\hbox{and}\qquad
\ker \Phi_{k + {1 \over 2}}
= \CC\hbox{-span}\{ x_d\ |\ \hbox{$d$ has more than $n$ blocks}
\}.
$$

\pf (a)
As a subgroup of $GL_n(\CC)$, $S_n$ acts on $V$ via
its permutation representation and $S_n$ acts on $V^{\otimes k}$ by
$$
\sigma (v_{i_1} \otimes v_{i_2} \otimes \cdots\otimes v_{i_k})
=  v_{\sigma(i_1)} \otimes
v_{\sigma(i_2)} \otimes \cdots \otimes v_{\sigma(i_k)}.  \formula
$$
Then $b \in \End_{S_n}(V^{\otimes k})$ if and only if
$\sigma^{-1} b \sigma = b$ (as endomorphisms on $V^{\otimes k}$)
for all $\sigma \in S_n$. Thus, using the notation of (3.1),
$b \in \End_{S_n}(V^{\otimes k})$ if and only if
$$
b^{i_1, \ldots, i_k}_{i_{1'}, \ldots, i_{k'}} =
(\sigma^{-1} b \sigma)^{i_1, \ldots, i_k}_{i_{1'}, \ldots, i_{k'}} =
b^{\sigma(i_1), \ldots, \sigma(i_k)}_{\sigma(i_{1'}),
\ldots, \sigma(i_{k'})},
\qquad\hbox{ for all } \sigma \in S_n.
$$
It follows that the matrix entries of $b$ are constant on the
$S_n$-orbits of its matrix coordinates. These orbits decompose
$\{1, \ldots, k, 1', \ldots, k'\}$ into
subsets and thus correspond to set partitions $d \in A_k$.
It follows from (2.3) and (3.1) that for all $d \in A_k$,
$$
(\Phi_k(x_d))^{i_1, \ldots, i_k}_{i_{1'}, \ldots, i_{k'}}
= \cases{
1, &if $i_r = i_s$ if and only if $r$ and $s$ are in the same block of $d$,
\cr 0, & otherwise.} \formula
$$
Thus
$\Phi_k(x_d)$ has  1s in the matrix positions corresponding to $d$
and 0s elsewhere, and so $b$ is a linear combination of
$\Phi_k(x_d), d \in A_k$.  Since
$x_d, d\in A_k$, form a basis of $\CC A_k$,
$\hbox{\rm im}~\Phi_k = \End_{S_n}(V^{\otimes k})$.

If $d$ has more than $n$ blocks, then by (3.8) the matrix entry
$(\Phi_k(x_d))^{i_1, \ldots, i_k}_{i_{1'}, \ldots, i_{k'}} = 0$
for all indices ${i_1, \ldots, i_k, i_{1'}, \ldots, i_{k'}}$,
since we need a distinct $i_j \in \{1, \ldots, n\}$ for each
block of $d$. Thus, $x_d \in \ker\Phi_k$.
If $d$ has $\le n$ blocks, then we can find an index set
${i_1, \ldots, i_k, i_{1'}, \ldots, i_{k'}}$ with
$(\Phi_k(x_d))^{i_1, \ldots, i_k}_{i_{1'}, \ldots, i_{k'}} = 1$
simply by choosing a distinct index
from $\{1, \ldots, n\}$ for each block
of $d$.  Thus, if $d$ has $\le n$ blocks then $x_d\not\in \ker\Phi_k$,
and so
$\ker \Phi_k
= \CC\hbox{-span}\{ x_d | \hbox{$d$ has more than $n$ blocks} \}$.

(b) The vector space $V^{\otimes k} \otimes v_n \subseteq
V^{\otimes (k+1)}$ is a submodule both for
$\CC \Akp \subseteq \CC A_{k+1}$ and $\CC S_{n-1} \subseteq \CC S_n$.
If $\sigma \in S_{n-1}$, then
$\sigma(v_{i_1} \otimes
\cdots \otimes v_{i_k} \otimes v_n)
= v_{\sigma(i_1)} \otimes \cdots \otimes v_{\sigma(i_k)} \otimes v_n$.
Then as above $b \in \End_{S_{n-1}}(V^{\otimes k})$ if and only if
$$
b^{i_1, \ldots, i_k,n}_{i_{1'}, \ldots, i_{k'},n} =
b^{\sigma(i_1), \ldots, \sigma(i_k),n}_{\sigma(i_{1'}),
\ldots, \sigma(i_{k'}),n},
\qquad\hbox{ for all } \sigma \in S_{n-1}.
$$
The $S_{n-1}$ orbits of the matrix coordinates of $b$ correspond to set
partitions  $d \in \Akp$; that is vertices $i_{k+1}$ and $i_{(k+1)'}$
must be in the same block of $d$.  The same argument as part (a)
can be used to show that $\ker \Phi_{k+ {1 \over 2}}$ is the span of
$x_d$ with $d \in \Akp$ having more than $n$ blocks.  We always choose
the index $n$ for the block containing $k+1$ and $(k+1)'$.
\endpf

\bigskip\noindent
{\it The maps
$\varepsilon_{1 \over  2}: \End(V^{\otimes k}) \to \End(V^{\otimes k})$
and
$\varepsilon^{1 \over  2}: \End(V^{\otimes k}) \to \End(V^{\otimes (k-1)})$
}
\medskip

If $b \in \End(V^{\otimes k})$ let $b^{i_1, \ldots, i_k}_{i_{1'}, \ldots,
i_{k'}} \in \CC$ be
the coefficients in the expansion
$$
b(v_{i_1} \otimes\cdots\otimes v_{i_k}) = \sum_{1 \le i_{1'}, \ldots i_{k'}
\le n}
b^{i_1, \ldots, i_k}_{i_{1'}, \ldots, i_{k'}}
v_{i_{1'}} \otimes\cdots\otimes v_{i_{k'}}. \formula
$$
Define linear maps
$$
\varepsilon_{1 \over  2}: \End(V^{\otimes k}) \to \End(V^{\otimes k})
\qquad\hbox{ and }\qquad
\varepsilon^{1 \over  2}: \End(V^{\otimes k}) \to \End(V^{\otimes (k-1)})
\qquad\hbox{by}
$$
$$
\varepsilon_{1 \over 2}(b)^{i_1, \ldots, i_k}_{i_{1'}, \ldots, i_{k'}} =
b^{i_1, \ldots, i_k}_{i_{1'}, \ldots, i_{k'}} \delta_{i_k i_{k'}}
\qquad\hbox{ and }\qquad
\varepsilon^{1 \over 2}(b)^{i_1, \ldots, i_{k-1}}_{i_{1'}, \ldots,
i_{(k-1)'}} =
\sum_{j,\ell = 1}^n
b^{i_1, \ldots, i_{k-1},j}_{i_{1'}, \ldots, i_{(k-1)', \ell}}.
\formula
$$
The composition of $\varepsilon_{1 \over 2}$ and $\varepsilon^{1\over 2}$
is the map
$$
\varepsilon_1: \End(V^{\otimes k}) \to \End(V^{\otimes (k-1)})
\qquad\hbox{given by}\qquad
\varepsilon_1(b)^{i_1, \ldots, i_{k-1}}_{i_{1'}, \ldots, i_{(k-1)'}} =
\sum_{j = 1}^n
b^{i_1, \ldots, i_{k-1},j}_{i_{1'}, \ldots, i_{(k-1)', j}},
\formula
$$
and
$$
\Tr(b)
= \varepsilon_1^k(b), \qquad\hbox{ for $b \in \End(V^{\otimes k})$}.
\formula
$$

The relation between the maps
$\varepsilon^{1\over2}$, $\varepsilon_{1\over2}$ in
(3.10) and the maps
$\varepsilon^{1\over2}$, $\varepsilon_{1\over2}$ in Section 2
is given by
$$\matrix{
\Phi_{k-{1\over2}}(\varepsilon_{1\over2}(b))
=\varepsilon_{1\over2}(\Phi_k(b))
\big\vert_{V^{\otimes(k-1)}\otimes v_n},
\hfill &\hbox{for $b\in \CC A_k(n)$,} \hfill \cr
\Phi_{k-1}(\varepsilon^{1\over2}(b))
={1\over n}\varepsilon^{1\over2}(\Phi_k(b)), \hfill
&\hbox{for $b\in \CC A_{k-{1\over2}}(n)$,\quad and} \hfill \cr
\Phi_{k-1}(\varepsilon_1(b)) =\varepsilon_1(\Phi_k(b)), \hfill
&\hbox{for $b\in \CC A_k(n)$,} \hfill \cr
}\formula$$
where, on the right hand side of the middle
equality $b$ is viewed as an element of $\CC A_k$ via the natural
inclusion $\CC A_{k-{1\over2}}(n)\subseteq \CC A_k(n)$.
Then
$$\Tr(\Phi_k(b)) = \varepsilon_1^k(\Phi_k(b))
=\Phi_0(\varepsilon_1^k(b))=\varepsilon_1^k(b)
=\tr_k(b),
\formula$$
and, by (3.4), if $b\in \CC A_{k-{1\over2}}(n)$ then
$$\Tr(\Phi_{k-{1\over2}}(b))
=\Tr(\Phi_k(b)\big\vert_{V^{\otimes(k-1)}\otimes v_n})
=\hbox{$1\over n$}\,\Tr(\Phi_k(b))
= \hbox{$1\over n$}\,\tr_k(b) = \hbox{$1\over n$}\, \tr_{k-{1\over2}}(b).
\formula$$

\bigskip\noindent
{\it The representations
{\rm $ (\Ind^{S_n}_{S_{n-1}} \Res^{S_n}_{S_{n-1}})^k ({\bf 1}_n) $}
and
{\rm
$\Res^{S_n}_{S_{n-1}}
(\Ind^{S_n}_{S_{n-1}} \Res^{S_n}_{S_{n-1}})^k ({\bf 1}_n)$
}
}
\medskip

Let ${\bf 1}_n = S_n^{(n)}$ be the trivial representation of $S_n$
and let $V = \CC\hbox{-span}\{ v_1,\ldots, v_n\}$ be the permutation
representation of $S_n$ given in (3.5).
Then
$$V \cong \Ind^{S_n}_{S_{n-1}}\Res^{S_n}_{S_{n-1}}({\bf 1}_n).
\formula$$
More generally, for any $S_n$-module $M$,
$$\eqalign{
\Ind^{S_n}_{S_{n-1}} \Res^{S_n}_{S_{n-1}} (M)
&\cong
\Ind^{S_n}_{S_{n-1}} (\Res^{S_n}_{S_{n-1}} (M) \otimes {\bf 1}_{n-1} )  \cr
&\cong \Ind^{S_n}_{S_{n-1}}
(\Res^{S_n}_{S_{n-1}} (M) \otimes \Res^{S_n}_{S_{n-1}}({\bf 1}_n) )  \cr
&\cong M \otimes \Ind^{S_n}_{S_{n-1}} \Res^{S_n}_{S_{n-1}} ({\bf 1}_n)
\cong M \otimes V,  \cr
}\formula$$
where the third isomorphism comes from the ``tensor identity,"
$$
\matrix{
\Ind^{S_n}_{S_{n-1}}
(\Res^{S_n}_{S_{n-1}} (M) \otimes N)
&\mapright{\sim}& M \otimes \Ind^{S_n}_{S_{n-1}} N  \cr
g \otimes (m \otimes n) & \mapsto & gm \otimes (g \otimes n) \cr
},
\formula
$$
for $g\in S_n$, $m \in M$, $n \in N$, and
the fact that $\Ind_{S_{n-1}}^{S_n} (W) = \CC S_n \otimes_{S_{n-1}} W$.
By iterating (3.17) it follows that
$$
(\Ind^{S_n}_{S_{n-1}} \Res^{S_n}_{S_{n-1}})^k ({\bf 1}) \cong V^{\otimes k}
\qquad\hbox{and}\qquad
\Res^{S_n}_{S_{n-1}}
(\Ind^{S_n}_{S_{n-1}} \Res^{S_n}_{S_{n-1}})^k ({\bf 1}) \cong V^{\otimes k}
\formula
$$
as $S_n$-modules and $S_{n-1}$-modules, respectively.

If
$$\lambda = (\lambda_1,\lambda_2, \ldots, \lambda_\ell)
\qquad\hbox{define}\qquad
\lambda_{> 1} = (\lambda_2, \ldots, \lambda_\ell)
\formula$$
to be the same partition as $\lambda$ except with the first row removed.
Build a graph $\hat A(n)$ which encodes the decomposition of
$V^{\otimes k}$, $k\in \ZZ_{\ge 0}$, by letting
$$\matrix{
\hbox{vertices on level $k$:}\hskip.4in
\hat A_k(n) = \{ \lambda\vdash n\ |\
k-|\lambda_{> 1}| \in \ZZ_{\ge 0}\}, \hfill \cr
\hbox{vertices on level $k+{1\over2}$:}\quad
\hat A_{k+{1\over2}}(n) = \{ \lambda\vdash n-1\ |\
k-|\lambda_{> 1}| \in \ZZ_{\ge 0}\},
\quad\hbox{and}  \hfill \cr
\hbox{an edge $\lambda\to\mu$,
if $\mu\in \hat A_{k+{1\over2}}(n)$ is obtained from
$\lambda\in \hat A_k(n)$ by removing a box}, \hfill \cr
\hbox{an edge $\mu\to \lambda$,
if $\lambda\in \hat A_{k+1}(n)$ is obtained from
$\mu\in \hat A_{k+{1\over2}}(n)$ by adding a box}. \cr
}
\formula$$
For example, if $n=5$ then the first few levels of $\hat A(n)$ are
$$
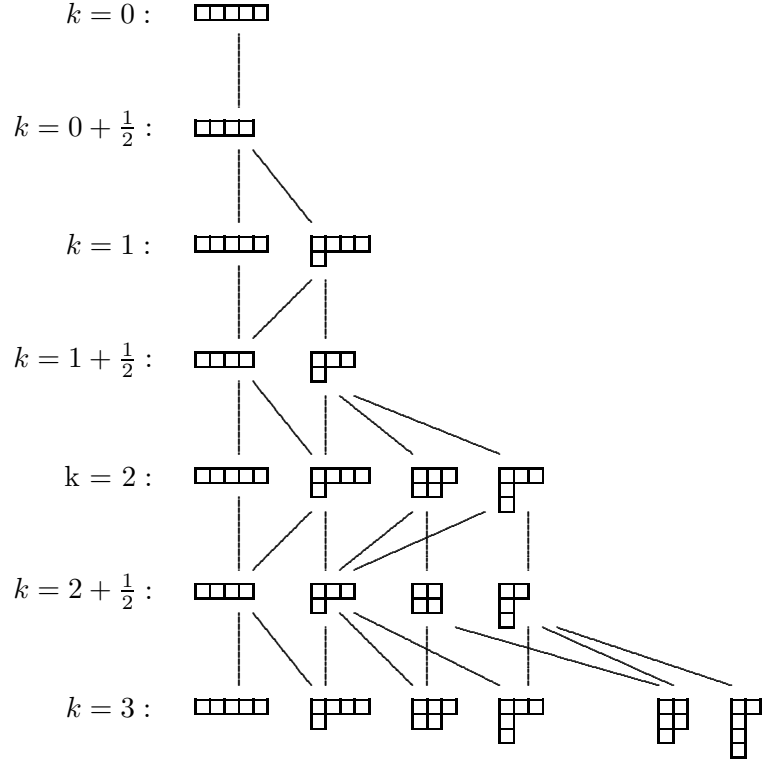
\matrix{
{\beginpicture
\setcoordinatesystem units <0.175cm,0.175cm>         
\setplotarea x from -8 to 40, y from 4 to 31   
\linethickness=0.5pt                          
\put{$k = 3:$} at  -6 -.5
\putrectangle corners at 0 -1 and 1 0
\putrectangle corners at 1 -1 and 2 0
\putrectangle corners at 2 -1 and 3 0
\putrectangle corners at 3 -1 and 4 0
\putrectangle corners at 4 -1 and 5 0
%
\putrectangle corners at 8 -1 and 9 0
\putrectangle corners at 9 -1 and 10 0
\putrectangle corners at 10 -1 and 11 0
\putrectangle corners at 11 -1 and 12 0
\putrectangle corners at 8 -2 and 9 -1
\putrectangle corners at 15 -1 and 16 0
\putrectangle corners at 16 -1 and 17 0
\putrectangle corners at 17 -1 and 18 0
\putrectangle corners at 15 -2 and 16 -1
\putrectangle corners at 16 -2 and 17 -1
\putrectangle corners at 21 -1 and 22 0
\putrectangle corners at 22 -1 and 23 0
\putrectangle corners at 23 -1 and 24 0
\putrectangle corners at 21 -2 and 22 -1
\putrectangle corners at 21 -3 and 22 -2
%
%
\putrectangle corners at 32 -1 and 33 0
\putrectangle corners at 33 -1 and 34 0
\putrectangle corners at 32 -2 and 33 -1
\putrectangle corners at 33 -2 and 34 -1
\putrectangle corners at 32 -3 and 33 -2
\putrectangle corners at 37 -1 and 38 0
\putrectangle corners at 38 -1 and 39 0
\putrectangle corners at 37 -2 and 38 -1
\putrectangle corners at 37 -3 and 38 -2
\putrectangle corners at 37 -4 and 38 -3
%
\plot 3 1 3 6  /
\plot 8 1 4 6  /
\plot 9 1 9 5  /
\plot 15 1 10 6  /
\plot 16 1 16 5  /
\plot 32 1 18 5  /
\plot 21 1 11 6  /
\plot 23 1 23 5  /
\plot 33 1 24 5  /
\plot 37 1 25 5  /
%
\put{$k = 2+{1\over2}:$} at  -7.75 7.5
\putrectangle corners at 0 7 and 1 8
\putrectangle corners at 1 7 and 2 8
\putrectangle corners at 2 7 and 3 8
\putrectangle corners at 3 7 and 4 8
%
\putrectangle corners at 8 7  and 9 8
\putrectangle corners at 9 7  and 10 8
\putrectangle corners at 10 7 and 11 8
\putrectangle corners at 8  6 and 9 7
\putrectangle corners at 15 7 and 16 8
\putrectangle corners at 16 7 and 17 8
\putrectangle corners at 15 6 and 16 7
\putrectangle corners at 16 6 and 17 7
\putrectangle corners at 21 7 and 22 8
\putrectangle corners at 22 7 and 23 8
\putrectangle corners at 21 6 and 22 7
\putrectangle corners at 21 5 and 22 6
%
\plot 3 9 3 14  /
\plot 8 13 4 9  /
\plot 9 13 9 9  /
\plot 15 13 10 9  /
\plot 16 13 16 9  /
\plot 20 13 11 9  /
\plot 23 13 23 9  /
%
%
\put{k = $2:$} at  -6 15.5
\putrectangle corners at 0  15 and 1  16
\putrectangle corners at 1  15 and 2  16
\putrectangle corners at 2  15 and 3  16
\putrectangle corners at 3  15 and 4  16
\putrectangle corners at 4  15 and 5  16
%
\putrectangle corners at 8  15 and 9   16
\putrectangle corners at 9  15 and 10  16
\putrectangle corners at 10 15 and 11  16
\putrectangle corners at 11 15 and 12  16
\putrectangle corners at 8  14 and 9  15
\putrectangle corners at 15  15 and 16  16
\putrectangle corners at 16  15 and 17  16
\putrectangle corners at 17  15 and 18  16
\putrectangle corners at 15  14 and 16  15
\putrectangle corners at 16  14 and 17  15
\putrectangle corners at 21  15 and 22  16
\putrectangle corners at 22  15 and 23  16
\putrectangle corners at 23  15 and 24  16
\putrectangle corners at 21  14 and 22  15
\putrectangle corners at 21  13 and 22 14
%
\plot 3 17 3 22  /
\plot 8 17 4 22  /
\plot 9 17 9 21  /
\plot 15 17 10 21  /
\plot 21 17 11 21  /
%
%
\put{$k = 1+{1\over2}:$} at  -7.75 23.5
\putrectangle corners at 0  23 and 1  24
\putrectangle corners at 1  23 and 2  24
\putrectangle corners at 2  23 and 3  24
\putrectangle corners at 3  23 and 4  24
%
%
\putrectangle corners at 8  23 and 9   24
\putrectangle corners at 9  23 and 10  24
\putrectangle corners at 10 23 and 11  24
\putrectangle corners at 8  22 and 9  23
%
\plot 3 25 3 30  /
\plot 8 29 4 25  /
\plot 9 29 9 25  /
%
%
\put{$k = 1:$} at  -6 31.5
\putrectangle corners at 0  31 and 1  32
\putrectangle corners at 1  31 and 2  32
\putrectangle corners at 2  31 and 3  32
\putrectangle corners at 3  31 and 4  32
\putrectangle corners at 4  31 and 5  32
%
\putrectangle corners at 8  31 and 9   32
\putrectangle corners at 9  31 and 10  32
\putrectangle corners at 10 31 and 11  32
\putrectangle corners at 11 31 and 12  32
\putrectangle corners at 8  30 and 9  31
%
\put{$k = 0+{1\over2}:$} at  -7.75 39.5
\putrectangle corners at 0  39 and 1  40
\putrectangle corners at 1  39 and 2  40
\putrectangle corners at 2  39 and 3  40
\putrectangle corners at 3  39 and 4  40
%
%
\plot 3 33 3 38  /
\plot 8 33 4 38  /
%
%
\put{$k = 0:$} at  -6 47.5
\putrectangle corners at 0  47 and 1  48
\putrectangle corners at 1  47 and 2  48
\putrectangle corners at 2  47 and 3  48
\putrectangle corners at 3  47 and 4  48
\putrectangle corners at 4  47 and 5  48
%
%
\plot 3 41 3 46  /
\endpicture} \cr
}
$$

The following theorem is a consequence of Theorem 3.5 and the Centralizer
Theorem, Theorem 5.4, (see also [GW, Theorem 3.3.7]).

\thm
Let $n, k \in \ZZ_{\ge 0}$.  Let
$S_n^\lambda$ denote the irreducible $S_n$-module indexed by $\lambda$.
\item{(a)} As $(\CC S_n,\CC A_k(n))$-bimodules,
$$
V^{\otimes k} \cong \bigoplus_{\lambda \in \hat A_k(n)}
S_n^\lambda \otimes A_k^\lambda(n),
$$
where the vector spaces $A_k^\lambda(n)$ are irreducible
$\CC A_k(n)$-modules and
$$\dim(A_k^\lambda(n)) =
(\hbox{number of paths from $(n)\in \hat A_0(n)$ to
$\lambda\in \hat A_k(n)$ in the graph $\hat A(n)$}).
$$
\item{(b)}
As $(\CC S_{n-1}, \CC A_{k+{1\over2}}(n))$-bimodules,
$$V^{\otimes k} \cong \bigoplus_{\mu \in \hat A_{k+{1\over2}}(n)}
S_{n-1}^\mu \otimes A_{k+{1\over2}}^\mu(n),
$$
where the vector spaces $A_{k+{1\over2}}^\mu(n)$ are irreducible
$\CC A_{k+{1\over2}}(n)$-modules and
$$\dim(A_{k+{1\over2}}^\mu(n)) =
(\hbox{number of paths from $(n)\in \hat A_0(n)$ to
$\mu\in \hat A_{k+{1\over2}}(n)$ in the graph $\hat A(n)$}).
$$
\endthm

\bigskip\noindent
{\it Determination of the polynomials $\tr^\mu(n)$}
\medskip

Let $n\in \ZZ_{>0}$.  For a partition $\lambda$, let
$$
\lambda_{>1} = (\lambda_2, \ldots, \lambda_\ell),
\qquad\hbox{ if }\qquad
\lambda = (\lambda_1,\lambda_2, \ldots, \lambda_\ell),
$$
i.e., remove the first row of $\lambda$ to get $\lambda_{>1}$.  Then,
for
$n\ge 2k$, the maps
$$\matrix{
\hat A_k(n) &\longleftrightarrow &\hat A_k \cr
\lambda &\longmapsto &\lambda_{> 1} \cr}
\qquad\hbox{are bijections}
\formula
$$
which provide an isomorphism between levels $0$ to $n$ of the graphs
$\hat A(n)$ and $\hat A$.

\prop
For $k\in {1\over2}\ZZ_{\ge 0}$ and $n\in \CC$ such that
$\CC A_k(n)$ is semisimple, let $\chi_{A_k(n)}^\mu$, $\mu\in
\hat A_k$, be the irreducible characters of $\CC A_\ell(n)$ and let
$\tr_k\colon \CC A_k(n)\to \CC$ be the trace on
$\CC A_k(n)$ defined in (2.25).  Use the notations for partitions
in (2.17). For $k > 0$ the coefficients in the expansion
$$\tr_k = \sum_{\mu\in \hat A_k} \tr^\mu(n)\chi_{A_k(n)}^\mu,
\qquad\hbox{are}\qquad
\tr^\mu(n) = \Big(\prod_{b\in \mu} {1\over h(b)} \Big)
\prod_{j=1}^{|\mu|} (n-|\mu|-(\mu_{j}-j)).$$
If $n\in \CC$ is such that $\CC A_{k+{1\over2}}(n)$
is semisimple then for $k > 0$ the coefficients in the expansion
$$
\tr_{k+{1\over2}}
= \sum_{\mu\in \hat A_{k+{1\over2}}}
\tr_{1\over2}^\mu(n)\chi_{A_{k+{1\over2}}(n)}^\mu,
\qquad\hbox{are}\qquad
\tr_{{1\over2}}^\mu(n) = \Big(\prod_{b\in \mu} {1\over h(b)} \Big)
\cdot n \cdot
\prod_{j=1}^{|\mu|} (n-1-|\mu|-(\mu_{j}-j)).
$$
\endprop

\pf
Let $\lambda$ be a partition with $n$ boxes.
Beginning with the vertical edge at the end of the
first row, label the boundary edges of $\lambda$
sequentially with $0,1,2\ldots, n$.  Then the
$$\eqalign{
\hbox{vertical edge label for row $i$}
&=
\hbox{(number of horizontal steps)}+
\hbox{(number of vertical steps)} \cr
&=(\lambda_1-\lambda_i)+(i-1) = (\lambda_1-1)-(\lambda_i-i),
\qquad\hbox{and the} \cr
\hbox{horizontal edge label for column $j$}
&=
\hbox{(number of horizontal steps)}+
\hbox{(number of vertical steps)} \cr
&=(\lambda_1-j+1)+(\lambda_j'-1)
=(\lambda_1-1)+(\lambda_j'-j)+1. \cr
}$$
Hence
$$\eqalign{
\{1,2,\ldots,n\}
&= \{(\lambda_1-1)-(\lambda'_j-j)+1\ |\ 1\le j\le \lambda_1\}
\sqcup
\{(\lambda_1-1)-(\lambda_i-i)\ |\ 2\le i\le n-\lambda_1+1\} \cr
&= \{h(b)\ |\ \hbox{$b$ is in row $1$ of $\lambda$}\}
\sqcup
\{(\lambda_1-1)-(\lambda_i-i)\ |\ 2\le i\le n-\lambda_1+1\}. \cr
}$$
For example, if $\lambda = (10,7,3,3,1)\vdash 24$, then the boundary
labels of $\lambda$ and the hook numbers in the first row of $\lambda$
are
$$
{\beginpicture
\setcoordinatesystem units <0.5cm,0.5cm> 
\setplotarea x from 0 to 16, y from 3 to 10    
\linethickness=0.5pt
\putrectangle corners at 13 9 and 14 10
\putrectangle corners at 14 9 and 15 10
\putrectangle corners at 12 9 and 13 10
\putrectangle corners at 11 9 and 12 10
\putrectangle corners at 11 8 and 12 10
\putrectangle corners at 10 8 and 11 10
\putrectangle corners at 9 8 and 10 10
\putrectangle corners at 8 8 and 9 10
\putrectangle corners at 7 8 and 8 10
\putrectangle corners at 5 7 and 8 8
\putrectangle corners at 5 6 and 8 7
\putrectangle corners at 6 6 and 7 10
\putrectangle corners at 5 6 and 6 10
\putrectangle corners at 5 5 and 6 10
\putrectangle corners at 5 8 and 8 9
\putrectangle corners at 5 9 and 12 8
\putrectangle corners at 5 -.5 and 5 9
\put{$\scriptstyle{24}$} at 5.3 -.2
\put{$\scriptstyle{23}$} at 5.3 .3
\put{$\scriptstyle{22}$} at 5.3 .8
\put{$\scriptstyle{21}$} at 5.3 1.3
\put{$\scriptstyle{20}$} at 5.3 1.8
\put{$\scriptstyle{19}$} at 5.3 2.3
\put{$\scriptstyle{18}$} at 5.3 2.8
\put{$\scriptstyle{17}$} at 5.3 3.3
\put{$\scriptstyle{16}$} at 5.3 3.8
\put{$\scriptstyle{15}$} at 5.3 4.3
\put{$\scriptstyle{14}$} at 5.5 4.75
\put{$\scriptstyle{13}$} at 6.3 5.3
\put{$\scriptstyle{12}$} at 6.5 5.75
\put{$\scriptstyle{11}$} at 7.5 5.75
\put{$\scriptstyle{10}$} at 8.3 6.5
\put{$\scriptstyle{9}$} at 8.25 7.35
\put{$\scriptstyle{8}$} at 8.6 7.75
\put{$\scriptstyle{7}$} at 9.5 7.75
\put{$\scriptstyle{6}$} at 10.5 7.75
\put{$\scriptstyle{5}$} at 11.5 7.75
\put{$\scriptstyle{4}$} at 12.25 8.35
\put{$\scriptstyle{3}$} at 12.6 8.75
\put{$\scriptstyle{2}$} at 13.5 8.75
\put{$\scriptstyle{1}$} at 14.5 8.75
\put{$\scriptstyle{0}$} at 15.25 9.5
\put{${14}$} at 5.5 9.5
\put{${12}$} at 6.5 9.5
\put{${11}$} at 7.5 9.5
\put{${8}$} at 8.5 9.5
\put{${7}$} at 9.5 9.5
\put{${6}$} at 10.5 9.5
\put{${5}$} at 11.5 9.5
\put{${3}$} at 12.5 9.5
\put{${2}$} at 13.5 9.5
\put{${1}$} at 14.5 9.5
\endpicture}.
$$
Thus,
since $\lambda_1 = n-|\lambda_{> 1}|$,
$$\dim(S_n^\lambda) = {n!\over \prod_{b\in \lambda} h(b)}
=
\Big(\prod_{b\in \lambda_{> 1}} {1\over h(b)} \Big)
\prod_{i=2}^{|\lambda_{> 1}|+1} (n-|\lambda_{> 1}|-(\lambda_i-(i-1))).
\formula
$$

Let $n\in \ZZ_{>0}$ and let $\chi_{S_n}^\lambda$ denote the
irreducible characters of the symmetric group $S_n$.
By taking the trace on both sides of the equality in Theorem 3.22,
$$\Tr(b,V^{\otimes k})
= \sum_{\lambda\in \hat A_k(n)}
\chi_{S_n}^\lambda(1)\chi_{A_k(n)}(b)
= \sum_{\lambda\in \hat A_k(n)}
\dim(S_n^\lambda)\chi_{A_k(n)}(b),
\qquad\hbox{for $b\in \CC A_k(n)$.}
$$
Thus the equality in (3.25) and the
bijection in (3.23) provide the expansion of $\tr_k$
for all $n\in \ZZ_{\ge 0}$ such that
$n\ge 2k$.  The statement for all $n\in \CC$ such
that $\CC A_k(n)$ is semisimple is then a consequence of
the fact that any polynomial is determined by its
evaluations at an infinite number of values of the parameter.
The proof of the expansion of $\tr_{k+{1\over2}}$ is exactly
analogous.
\endpf

Note that the polynomials $\tr^\mu(n)$ and $\tr_{{1\over2}}^\mu(n)$
(of degrees $|\mu|$ and $|\mu|+1$, respectively)
do not depend on $k$.  By Proposition 3.24,
$$\eqalign{
&\{ \hbox{roots of $\tr_{1\over2}^\mu(n)$}\ |\ \mu\in \hat
A_{1\over2}\}=\{0\}, \cr
&\{ \hbox{roots of $\tr_1^\mu(n)$}\ |\ \mu\in \hat A_1\}=\{1\}, \cr
&\{ \hbox{roots of $\tr_{1{1\over2}}^\mu(n)$}\ |\
\mu\in \hat A_{1{1\over2}}\}=\{0,2\}, \qquad\hbox{and}
\cr
&\{ \hbox{roots of $\tr_k^\mu(n)$}\ |\ \mu\in \hat A_k\}
=\{ 0, 1,\ldots,2k-1\},
\qquad\hbox{for $k\in {1\over2}\ZZ_{\ge 0}$, $k\ge 2$.}
\cr} \formula
$$
For example, the first few values of
$\tr^\mu$ and $\tr_{1\over2}^\mu$ are
$$
\matrix{
\tr^\emptyset(n) = 1,
\hfill&\quad&
\tr^\emptyset_{1 \over 2}(n) = n
\hfill\cr
\tr^{{\beginpicture\setcoordinatesystem units <0.1cm,0.1cm>
\putrectangle corners at 0 0 and 1 1  \endpicture}}(n) = n-1,
\hfill&\quad&
\tr_{1 \over 2}^{{\beginpicture\setcoordinatesystem units <0.1cm,0.1cm>
\putrectangle corners at 0 0 and 1 1  \endpicture}}(n) = n(n-2),
\hfill\cr
\tr^{{\beginpicture\setcoordinatesystem units <0.1cm,0.1cm>
\putrectangle corners at 0 0 and 1 1
\putrectangle corners at 1 0 and 2 1 \endpicture}}(n) = {1 \over 2}n(n-3),
\hfill&\quad&
\tr_{1 \over 2}^{{\beginpicture\setcoordinatesystem units <0.1cm,0.1cm>
\putrectangle corners at 0 0 and 1 1
\putrectangle corners at 1 0 and 2 1 \endpicture}}(n) = {1 \over
2}n(n-1)(n-4),
\hfill\cr
\tr^{{\beginpicture\setcoordinatesystem units <0.1cm,0.1cm>
\putrectangle corners at 0 0 and 1 1
\putrectangle corners at 0 1 and 1 2 \endpicture}}(n) = {1 \over
2}(n-1)(n-2),
\hfill&\quad&
\tr_{1 \over 2}^{{\beginpicture\setcoordinatesystem units <0.1cm,0.1cm>
\putrectangle corners at 0 0 and 1 1
\putrectangle corners at 0 1 and 1 2 \endpicture}}(n) = {1 \over
2}n(n-2)(n-3),
\hfill\cr
\tr^{{\beginpicture\setcoordinatesystem units <0.1cm,0.1cm>
\putrectangle corners at 0 0 and 1 1
\putrectangle corners at 1 0 and 2 1
\putrectangle corners at 2 0 and 3 1
\endpicture}}(n) = {1 \over 6}n(n-1)(n-5),
\hfill&\quad&
\tr_{1 \over 2}^{{\beginpicture\setcoordinatesystem units <0.1cm,0.1cm>
\putrectangle corners at 0 0 and 1 1 \putrectangle corners at 2 0 and 3 1
\putrectangle corners at 1 0 and 2 1 \endpicture}}(n) = {1 \over
6}n(n-1)(n-2)(n-6),
\hfill\cr
\tr^{{\beginpicture\setcoordinatesystem units <0.1cm,0.1cm>
\putrectangle corners at 0 0 and 1 1
\putrectangle corners at 1 0 and 2 1
\putrectangle corners at 0 -1 and 1 0
\endpicture}}(n) = {1 \over 6}n(n-2)(n-4),
\hfill&\quad&
\tr_{1 \over 2}^{{\beginpicture\setcoordinatesystem units <0.1cm,0.1cm>
\putrectangle corners at 0 0 and 1 1 \putrectangle corners at 0 -1 and 1 0
\putrectangle corners at 1 0 and 2 1 \endpicture}}(n) = {1 \over
6}n(n-1)(n-3)(n-5),
\hfill\cr
\tr^{{\beginpicture\setcoordinatesystem units <0.1cm,0.1cm>
\putrectangle corners at 0 0 and 1 1
\putrectangle corners at 0 -2 and 1 -1
\putrectangle corners at 0 -1 and 1 0
\endpicture}}(n) = {1 \over 6}(n-1)(n-2)(n-3),
\hfill&\quad&
\tr_{1 \over 2}^{{\beginpicture\setcoordinatesystem units <0.1cm,0.1cm>
\putrectangle corners at 0 0 and 1 1
\putrectangle corners at 0 -2 and 1 -1
\putrectangle corners at 0 -1 and 1 0 \endpicture}}(n) = {1 \over
6}n(n-2)(n-3)(n-4),
\hfill\cr
}
$$

\thm Let $n\in \ZZ_{\ge2}$ and $k\in {1\over2}\ZZ_{\ge 0}$.  Then
$$\hbox{$\CC A_k(n)$ is semisimple
\quad if and only if\quad
$k\le {n+1\over2}$.}
$$

\pf
By Theorem 2.26(a) and the observation  (3.26) it follows that
$\CC A_k(n)$ is semisimple if $n\ge 2k-1$.

Suppose $n$ is even.  Then Theorems 2.26(a) and 2.26(b) imply that
$$
\hbox{ $\CC A_{{n\over2}+{1\over2}}(n)$ is semisimple
\quad and\quad
$\CC A_{{n\over2}+1}(n)$ is not semisimple,}$$
since $(n/2)\in \hat A_{{n\over2}+{1\over2}}$
and $\tr_{1\over2}^{(n/2)}(n)=0$.
Since $(n/2)\in \hat A_{{n\over2}+1}(n)$, the $A_{{n\over2}+1}(n)$-module
$A_{{n\over2}+1}^{(n/2)}(n) \ne 0$.
Since the path $(\emptyset,\ldots,(n/2),(n/2),(n/2))
\in \hat A_{{n\over2}+1}^{(n/2)}$
does not correspond to an element of
$\hat A_{{n\over2}+1}^{(n/2)}(n)$,
$$\Card (\hat A_{{n\over2}+1}^{(n/2)})\ne
\Card(\hat A_{{n\over2}+1}^{(n/2)}(n)).$$
Thus, Tits deformation theorem (Theorem 5.13) implies that $\CC
A_{{n\over2}+1}(n)$ is cannot be semisimple.  Now it follows from Theorem
Theorem 2.26(c) that
$\CC A_k(n)$ is not semisimple for $k\ge {n\over2}+{1\over2}$.

If $n$ is odd then Theorems 2.26(a) and 2.26(b) imply that
$$
\hbox{$\CC A_{{n\over2}+{1\over2}}(n)$ is semisimple
\quad and\quad
$\CC A_{{n\over2}+1}(n)$ is not semisimple,}$$
since $(n/2)\in \hat A_{{n\over2}+{1\over2}}$
and $\tr^{(n/2)}(n)=0$.
Since $({n\over2}-{1\over2})\in \hat A_{{n\over2}+1}(n)$,
the $A_{{n\over2}+1}(n)$-module
$A_{{n\over2}+1}^{({n\over2}-{1\over2})}(n) \ne 0$.
Since the path $(\emptyset,\ldots,({n\over2}-{1\over2}),
({n\over2}+{1\over2}),({n\over2}-{1\over2}))
\in \hat A_{{n\over2}+1}^{({n\over2}-{1\over2})}$
does not correspond to an element
of $\hat A_{{n\over2}+1}^{({n\over2}-{1\over2})}(n)$,
and since
$$\Card(\hat A_{{n\over2}+1}^{({n\over2}-{1\over2})})\ne
\Card(\hat A_{{n\over2}+1}^{({n\over2}-{1\over2})}(n))$$
the Tits deformation theorem implies that $\CC A_{{n\over2}+1}(n)$ is
not semisimple.  Now it follows from Theorem 2.26(c) that
$\CC A_k(n)$ is not semisimple for $k\ge {n\over2}+{1\over2}$.
\endpf

\bigskip\noindent
{\it Murphy elements for $\CC A_k(n)$}
\medskip

Let $\kappa_n$ be the element of $\CC S_n$ given by
$$
\kappa_n = \sum_{1 \le \ell < m \le n} s_{\ell m}\,,
\formula
$$
where $s_{\ell m}$ is the transposition in $S_n$ which switches $\ell$
and $m$. Let $S\subseteq\{1,2,\ldots, k\}$  and let $I\subseteq S\cup S'$.
Define $b_S, d_I \in A_k$ by
$$b_S = \{ S\cup S', \{\ell,\ell'\}_{\ell\not\in S}\}
\qquad\hbox{and}\qquad
d_{I\subseteq S} = \{ I, I^c, \{\ell,\ell'\}_{\ell\not\in S}\}.
\formula$$
For example, in $A_9$, if $S= \{2,4,5,8\}$ and
$I = \{2, 4, 4', 5, 8\}$ then
$$b_S
={\beginpicture
\setcoordinatesystem units <0.5cm,0.2cm> 
\setplotarea x from 1 to 9, y from -1 to 2   
\linethickness=0.5pt                        
\put{$\bullet$} at 1 -1 \put{$\bullet$} at 1 2
\put{$\bullet$} at 2 -1 \put{$\bullet$} at 2 2
\put{$\bullet$} at 3 -1 \put{$\bullet$} at 3 2
\put{$\bullet$} at 4 -1 \put{$\bullet$} at 4 2
\put{$\bullet$} at 5 -1 \put{$\bullet$} at 5 2
\put{$\bullet$} at 6 -1 \put{$\bullet$} at 6 2
\put{$\bullet$} at 7 -1 \put{$\bullet$} at 7 2
\put{$\bullet$} at 8 -1 \put{$\bullet$} at 8 2
\put{$\bullet$} at 9 -1 \put{$\bullet$} at 9 2
\plot 1 -1 1 2 /
\plot 3 -1 3 2 /
\plot 6 -1 6 2 /
\plot 7 -1 7 2 /
\plot 9 -1 9 2 /
\plot 2 -1 2 2 /
\plot 8 -1 8 2 /
\setquadratic
\plot 2 2   3 1   4 2 /
\plot 4 2   4.5 1   5 2 /
\plot 5 2   6.5 1   8 2 /
\plot 2 -1   3 0   4 -1 /
\plot 4 -1   4.5 0   5 -1 /
\plot 5 -1   6.5 0   8 -1 /
\endpicture}
\qquad\hbox{and}\qquad
d_{I\subseteq S}
={\beginpicture
\setcoordinatesystem units <0.5cm,0.2cm> 
\setplotarea x from 1 to 9, y from -1 to 2   
\linethickness=0.5pt                        
\put{$\bullet$} at 1 -1 \put{$\bullet$} at 1 2
\put{$\bullet$} at 2 -1 \put{$\bullet$} at 2 2
\put{$\bullet$} at 3 -1 \put{$\bullet$} at 3 2
\put{$\bullet$} at 4 -1 \put{$\bullet$} at 4 2
\put{$\bullet$} at 5 -1 \put{$\bullet$} at 5 2
\put{$\bullet$} at 6 -1 \put{$\bullet$} at 6 2
\put{$\bullet$} at 7 -1 \put{$\bullet$} at 7 2
\put{$\bullet$} at 8 -1 \put{$\bullet$} at 8 2
\put{$\bullet$} at 9 -1 \put{$\bullet$} at 9 2
\plot 1 -1 1 2 /
\plot 3 -1 3 2 /
\plot 6 -1 6 2 /
\plot 7 -1 7 2 /
\plot 9 -1 9 2 /
\plot 4 -1 4 2 /
\setquadratic
\plot 2 2   3 1   4 2 /
\plot 4 2   4.5 1   5 2 /
\plot 5 2   6.5 1   8 2 /
\plot 2 -1   3.5 0   5 -1 /
\plot 5 -1   6.5 0   8 -1 /
\endpicture}.$$
For $S\subseteq\{1,2,\ldots, k\}$ define
$$p_S = \sum_I {1 \over 2}
(-1)^{\#(\{\ell,\ell'\}\subseteq I)+\#(\{\ell,\ell'\}\subseteq I^c)}
d_I,
\formula$$
where the sum is over $I\subseteq S\cup S'$ such that
$I\ne \emptyset$, $I\ne S\cup S'$, $I\ne \{\ell,\ell'\}$
and $I\ne \{\ell,\ell'\}^c$.
For $S\subseteq \{1,\ldots, k+1\}$ such that $k+1\in S$ define
$$\tilde p_S = \sum_I {1 \over 2}
(-1)^{\#(\{\ell,\ell'\}\subseteq I)+\#(\{\ell,\ell'\}\subseteq I^c)}
d_I,
\formula$$
where the sum is over all $I\subseteq S\cup S'$ such that
$\{k+1,(k+1)'\}\subseteq I$ or $\{k+1,(k+1)'\}\subseteq I^c$,
$I\ne S\cup S'$, $I\ne \{k+1,(k+1)'\}$ and $I\ne \{k+1,(k+1)'\}^c$.

Let $Z_1 = 1$ and, for $k \in \ZZ_{>1}$, let
$$
Z_k = {k \choose 2} + \sum_{S\subseteq \{1,\ldots, k\}\atop |S|\ge 1} p_S
+ \sum_{S\subseteq \{1,\ldots, k\}\atop |S|\ge 2}
(n-k+|S|))(-1)^{|S|} b_S .
\formula
$$
View $Z_k \in \CC A_k \subseteq \CC \Akp$ using the embedding in
(2.2), and define $Z_{1 \over 2} = 1$ and
$$
Z_{k + {1 \over 2}} = k + Z_k +
\sum_{|S|\ge 2\atop k+1\in S}
\tilde p_S +(n-(k+1)+|S|)(-1)^{|S|}b_S,
\formula$$
where the sum is over $S\subseteq \{1,\ldots, k+1\}$ such
that $k+1\in S$ and $|S|\ge 2$.
Define
$$M_{1\over 2} = 1,\qquad\hbox{and}\qquad
M_k = Z_k - Z_{k - {1 \over 2}},\quad
\hbox{for $k \in {1 \over 2} \ZZ_{>0}$.}
\formula
$$

For example, the first few $Z_k$ are
$$
\eqalign{
Z_0  &= 1, \qquad\qquad Z_{1 \over 2} = 1, \qquad\qquad Z_{1} =
\underbrace{{\beginpicture
\setcoordinatesystem units <0.35cm,0.12cm> 
\setplotarea x from 1 to 1, y from -1 to 2   
\linethickness=0.5pt                        
\put{$\bullet$} at 1 -1 \put{$\bullet$} at 1 2
\endpicture}}_{p_{\{1\}}},
\cr
Z_{1 {1\over 2}} &=
{\beginpicture
\setcoordinatesystem units <0.35cm,0.12cm> 
\setplotarea x from 1 to 2, y from -1 to 2   
\linethickness=0.5pt                        
\put{$\bullet$} at 1 -1 \put{$\bullet$} at 1 2
\put{$\bullet$} at 2 -1 \put{$\bullet$} at 2 2
\plot 1 -1 1 2 /
\plot 2 -1 2 2 /
\endpicture}
+\underbrace{
{\beginpicture
\setcoordinatesystem units <0.35cm,0.12cm> 
\setplotarea x from 1 to 2, y from -1 to 2   
\linethickness=0.5pt                        
\put{$\bullet$} at 1 -1 \put{$\bullet$} at 1 2
\put{$\bullet$} at 2 -1 \put{$\bullet$} at 2 2
\plot 2 -1 2 2 /
\endpicture}}_{Z_1}
-\underbrace{
{\beginpicture
\setcoordinatesystem units <0.35cm,0.12cm> 
\setplotarea x from 1 to 2, y from -1 to 2   
\linethickness=0.5pt                        
\put{$\bullet$} at 1 -1 \put{$\bullet$} at 1 2
\put{$\bullet$} at 2 -1 \put{$\bullet$} at 2 2
\plot 1 -1 2 -1 /
\plot 2 -1 2 2 /
\endpicture}
-{\beginpicture
\setcoordinatesystem units <0.35cm,0.12cm> 
\setplotarea x from 1 to 2, y from -1 to 2   
\linethickness=0.5pt                        
\put{$\bullet$} at 1 -1 \put{$\bullet$} at 1 2
\put{$\bullet$} at 2 -1 \put{$\bullet$} at 2 2
\plot 1 2 2 2 /
\plot 2 -1 2 2 /
\endpicture}
}_{\tilde p_{\{1,2\}}}
+ n
\underbrace{\beginpicture
\setcoordinatesystem units <0.35cm,0.12cm> 
\setplotarea x from 1 to 2, y from -1 to 2   
\linethickness=0.5pt                        
\put{$\bullet$} at 1 -1 \put{$\bullet$} at 1 2
\put{$\bullet$} at 2 -1 \put{$\bullet$} at 2 2
\plot 1 -1 1 2 2 2 2 -1 1 -1 /
\endpicture}_{b_{\{1,2\}}},\qquad\hbox{and} \cr
Z_2 &={\beginpicture
\setcoordinatesystem units <0.35cm,0.12cm> 
\setplotarea x from 1 to 2, y from -1 to 2   
\linethickness=0.5pt                        
\put{$\bullet$} at 1 -1 \put{$\bullet$} at 1 2
\put{$\bullet$} at 2 -1 \put{$\bullet$} at 2 2
\plot 1 -1 1 2 /
\plot 2 -1 2 2 /
\endpicture}
+
\underbrace{
{\beginpicture
\setcoordinatesystem units <0.35cm,0.12cm> 
\setplotarea x from 1 to 2, y from -1 to 2   
\linethickness=0.5pt                        
\put{$\bullet$} at 1 -1 \put{$\bullet$} at 1 2
\put{$\bullet$} at 2 -1 \put{$\bullet$} at 2 2
\plot 2 -1 2 2 /
\endpicture}}_{p_{\{1\}}}
+\underbrace{
{\beginpicture
\setcoordinatesystem units <0.35cm,0.12cm> 
\setplotarea x from 1 to 2, y from -1 to 2   
\linethickness=0.5pt                        
\put{$\bullet$} at 1 -1 \put{$\bullet$} at 1 2
\put{$\bullet$} at 2 -1 \put{$\bullet$} at 2 2
\plot 1 -1 1 2 /
\endpicture}}_{p_{\{2\}}}
-\underbrace{
{\beginpicture
\setcoordinatesystem units <0.35cm,0.12cm> 
\setplotarea x from 1 to 2, y from -1 to 2   
\linethickness=0.5pt                        
\put{$\bullet$} at 1 -1 \put{$\bullet$} at 1 2
\put{$\bullet$} at 2 -1 \put{$\bullet$} at 2 2
\plot 1 -1 2 -1 /
\plot 2 -1 2 2 /
\endpicture}
-{\beginpicture
\setcoordinatesystem units <0.35cm,0.12cm> 
\setplotarea x from 1 to 2, y from -1 to 2   
\linethickness=0.5pt                        
\put{$\bullet$} at 1 -1 \put{$\bullet$} at 1 2
\put{$\bullet$} at 2 -1 \put{$\bullet$} at 2 2
\plot 1 2 1 -1 /
\plot 1 -1 2 -1 /
\endpicture}
-{\beginpicture
\setcoordinatesystem units <0.35cm,0.12cm> 
\setplotarea x from 1 to 2, y from -1 to 2   
\linethickness=0.5pt                        
\put{$\bullet$} at 1 -1 \put{$\bullet$} at 1 2
\put{$\bullet$} at 2 -1 \put{$\bullet$} at 2 2
\plot 1 2 1 -1 /
\plot 1 2 2 2 /
\endpicture}
-{\beginpicture
\setcoordinatesystem units <0.35cm,0.12cm> 
\setplotarea x from 1 to 2, y from -1 to 2   
\linethickness=0.5pt                        
\put{$\bullet$} at 1 -1 \put{$\bullet$} at 1 2
\put{$\bullet$} at 2 -1 \put{$\bullet$} at 2 2
\plot 1 2 2 2 /
\plot 2 -1 2 2 /
\endpicture}
+
{\beginpicture
\setcoordinatesystem units <0.35cm,0.12cm> 
\setplotarea x from 1 to 2, y from -1 to 2   
\linethickness=0.5pt                        
\put{$\bullet$} at 1 -1 \put{$\bullet$} at 1 2
\put{$\bullet$} at 2 -1 \put{$\bullet$} at 2 2
\plot 1 -1 2 2 /
\plot 2 -1 1 2 /
\endpicture}
+
{\beginpicture
\setcoordinatesystem units <0.35cm,0.12cm> 
\setplotarea x from 1 to 2, y from -1 to 2   
\linethickness=0.5pt                        
\put{$\bullet$} at 1 -1 \put{$\bullet$} at 1 2
\put{$\bullet$} at 2 -1 \put{$\bullet$} at 2 2
\plot 1 -1 2 -1 /
\plot 1 2 2 2 /
\endpicture}}_{p_{\{1,2\}}}
+ n
\underbrace{\beginpicture
\setcoordinatesystem units <0.35cm,0.12cm> 
\setplotarea x from 1 to 2, y from -1 to 2   
\linethickness=0.5pt                        
\put{$\bullet$} at 1 -1 \put{$\bullet$} at 1 2
\put{$\bullet$} at 2 -1 \put{$\bullet$} at 2 2
\plot 1 -1 1 2 2 2 2 -1 1 -1 /
\endpicture}_{b_{\{1,2\}}},
}$$
and the first few $M_k$ are
$$
\eqalign{
M_0 & = 1, \qquad
M_{1 \over 2} = 1, \qquad
M_1 =
 {\beginpicture
\setcoordinatesystem units <0.35cm,0.12cm> 
\setplotarea x from 1 to 1, y from -1 to 2   
\linethickness=0.5pt                        
\put{$\bullet$} at 1 -1 \put{$\bullet$} at 1 2
\endpicture}
-
{\beginpicture
\setcoordinatesystem units <0.35cm,0.12cm> 
\setplotarea x from 1 to 1, y from -1 to 2   
\linethickness=0.5pt                        
\put{$\bullet$} at 1 -1 \put{$\bullet$} at 1 2
\plot 1 -1 1 2 /
\endpicture}\,,
\qquad
M_{1 {1 \over 2}}  =
{\beginpicture
\setcoordinatesystem units <0.35cm,0.12cm> 
\setplotarea x from 1 to 2, y from -1 to 2   
\linethickness=0.5pt                        
\put{$\bullet$} at 1 -1 \put{$\bullet$} at 1 2
\put{$\bullet$} at 2 -1 \put{$\bullet$} at 2 2
\plot 1 -1 1 2 /
\plot 2 -1 2 2 /
\endpicture}
-
{\beginpicture
\setcoordinatesystem units <0.35cm,0.12cm> 
\setplotarea x from 1 to 2, y from -1 to 2   
\linethickness=0.5pt                        
\put{$\bullet$} at 1 -1 \put{$\bullet$} at 1 2
\put{$\bullet$} at 2 -1 \put{$\bullet$} at 2 2
\plot 1 -1 2 -1 /
\plot 2 -1 2 2 /
\endpicture}
-
{\beginpicture
\setcoordinatesystem units <0.35cm,0.12cm> 
\setplotarea x from 1 to 2, y from -1 to 2   
\linethickness=0.5pt                        
\put{$\bullet$} at 1 -1 \put{$\bullet$} at 1 2
\put{$\bullet$} at 2 -1 \put{$\bullet$} at 2 2
\plot 1 2 2 2 /
\plot 2 -1 2 2 /
\endpicture}
+ n
 {\beginpicture
\setcoordinatesystem units <0.35cm,0.12cm> 
\setplotarea x from 1 to 2, y from -1 to 2   
\linethickness=0.5pt                        
\put{$\bullet$} at 1 -1 \put{$\bullet$} at 1 2
\put{$\bullet$} at 2 -1 \put{$\bullet$} at 2 2
\plot 1 -1 1 2 2 2 2 -1 1 -1 /
\endpicture}\,,
\cr
\cr
M_2 & =
{\beginpicture
\setcoordinatesystem units <0.35cm,0.12cm> 
\setplotarea x from 1 to 2, y from -1 to 2   
\linethickness=0.5pt                        
\put{$\bullet$} at 1 -1 \put{$\bullet$} at 1 2
\put{$\bullet$} at 2 -1 \put{$\bullet$} at 2 2
\plot 1 -1 1 2 /
\endpicture}
-{\beginpicture
\setcoordinatesystem units <0.35cm,0.12cm> 
\setplotarea x from 1 to 2, y from -1 to 2   
\linethickness=0.5pt                        
\put{$\bullet$} at 1 -1 \put{$\bullet$} at 1 2
\put{$\bullet$} at 2 -1 \put{$\bullet$} at 2 2
\plot 1 2 1 -1 /
\plot 1 -1 2 -1 /
\endpicture}
-{\beginpicture
\setcoordinatesystem units <0.35cm,0.12cm> 
\setplotarea x from 1 to 2, y from -1 to 2   
\linethickness=0.5pt                        
\put{$\bullet$} at 1 -1 \put{$\bullet$} at 1 2
\put{$\bullet$} at 2 -1 \put{$\bullet$} at 2 2
\plot 1 2 1 -1 /
\plot 1 2 2 2 /
\endpicture}
+
{\beginpicture
\setcoordinatesystem units <0.35cm,0.12cm> 
\setplotarea x from 1 to 2, y from -1 to 2   
\linethickness=0.5pt                        
\put{$\bullet$} at 1 -1 \put{$\bullet$} at 1 2
\put{$\bullet$} at 2 -1 \put{$\bullet$} at 2 2
\plot 1 -1 2 2 /
\plot 2 -1 1 2 /
\endpicture}
+
{\beginpicture
\setcoordinatesystem units <0.35cm,0.12cm> 
\setplotarea x from 1 to 2, y from -1 to 2   
\linethickness=0.5pt                        
\put{$\bullet$} at 1 -1 \put{$\bullet$} at 1 2
\put{$\bullet$} at 2 -1 \put{$\bullet$} at 2 2
\plot 1 -1 2 -1 /
\plot 1 2 2 2 /
\endpicture}\,, \qquad\hbox{and}
\cr
\cr
M_{2 {1 \over 2}} &=
2 {\beginpicture
\setcoordinatesystem units <0.35cm,0.12cm> 
\setplotarea x from 1 to 3, y from -1 to 2   
\linethickness=0.5pt                        
\put{$\bullet$} at 1 -1 \put{$\bullet$} at 1 2
\put{$\bullet$} at 2 -1 \put{$\bullet$} at 2 2
\put{$\bullet$} at 3 -1 \put{$\bullet$} at 3 2
\plot 1 -1 1 2 /
\plot 2 -1 2 2 /
\plot 3 -1 3 2 /
\endpicture}
+
{\beginpicture
\setcoordinatesystem units <0.35cm,0.12cm> 
\setplotarea x from 1 to 3, y from -1 to 2   
\linethickness=0.5pt                        
\put{$\bullet$} at 1 -1 \put{$\bullet$} at 1 2
\put{$\bullet$} at 2 -1 \put{$\bullet$} at 2 2
\put{$\bullet$} at 3 -1 \put{$\bullet$} at 3 2
\plot 2 -1 2 2 /
\plot 3 -1 3 2 /
\setquadratic
\plot 1 -1 2 0 3 -1 /
\endpicture}
+
{\beginpicture
\setcoordinatesystem units <0.35cm,0.12cm> 
\setplotarea x from 1 to 3, y from -1 to 2   
\linethickness=0.5pt                        
\put{$\bullet$} at 1 -1 \put{$\bullet$} at 1 2
\put{$\bullet$} at 2 -1 \put{$\bullet$} at 2 2
\put{$\bullet$} at 3 -1 \put{$\bullet$} at 3 2
\plot 2 -1 2 2 /
\plot 3 -1 3 2 /
\setquadratic
\plot 1 2 2 1 3 2 /
\endpicture}
+
{\beginpicture
\setcoordinatesystem units <0.35cm,0.12cm> 
\setplotarea x from 1 to 3, y from -1 to 2   
\linethickness=0.5pt                        
\put{$\bullet$} at 1 -1 \put{$\bullet$} at 1 2
\put{$\bullet$} at 2 -1 \put{$\bullet$} at 2 2
\put{$\bullet$} at 3 -1 \put{$\bullet$} at 3 2
\plot 1 -1 1 2 /
\plot 2 -1 3 -1 3 2 /
\endpicture}
+
{\beginpicture
\setcoordinatesystem units <0.35cm,0.12cm> 
\setplotarea x from 1 to 3, y from -1 to 2   
\linethickness=0.5pt                        
\put{$\bullet$} at 1 -1 \put{$\bullet$} at 1 2
\put{$\bullet$} at 2 -1 \put{$\bullet$} at 2 2
\put{$\bullet$} at 3 -1 \put{$\bullet$} at 3 2
\plot 1 -1 1 2 /
\plot 3 -1 3 2 2 2 /
\endpicture}
+ (n-1)
{\beginpicture
\setcoordinatesystem units <0.35cm,0.12cm> 
\setplotarea x from 1 to 3, y from -1 to 2   
\linethickness=0.5pt                        
\put{$\bullet$} at 1 -1 \put{$\bullet$} at 1 2
\put{$\bullet$} at 2 -1 \put{$\bullet$} at 2 2
\put{$\bullet$} at 3 -1 \put{$\bullet$} at 3 2
\plot 1 -1 1 2 /
\plot 3 -1 3 2 2 2 2 -1 3 -1 /
\endpicture}
+ (n-1)
{\beginpicture
\setcoordinatesystem units <0.35cm,0.12cm> 
\setplotarea x from 1 to 3, y from -1 to 2   
\linethickness=0.5pt                        
\put{$\bullet$} at 1 -1 \put{$\bullet$} at 1 2
\put{$\bullet$} at 2 -1 \put{$\bullet$} at 2 2
\put{$\bullet$} at 3 -1 \put{$\bullet$} at 3 2
\plot 1 -1 1 2 /
\plot 2 -1 2 2 /
\plot 3 -1 3 2 /
\setquadratic
\plot 1 2 2 1 3 2 /
\plot 1 -1 2 0 3 -1 /
\endpicture}
\cr & \qquad
+
{\beginpicture
\setcoordinatesystem units <0.35cm,0.12cm> 
\setplotarea x from 1 to 3, y from -1 to 2   
\linethickness=0.5pt                        
\put{$\bullet$} at 1 -1 \put{$\bullet$} at 1 2
\put{$\bullet$} at 2 -1 \put{$\bullet$} at 2 2
\put{$\bullet$} at 3 -1 \put{$\bullet$} at 3 2
\plot 2 2 3 2 3 -1 1 -1 /
\endpicture}
+
{\beginpicture
\setcoordinatesystem units <0.35cm,0.12cm> 
\setplotarea x from 1 to 3, y from -1 to 2   
\linethickness=0.5pt                        
\put{$\bullet$} at 1 -1 \put{$\bullet$} at 1 2
\put{$\bullet$} at 2 -1 \put{$\bullet$} at 2 2
\put{$\bullet$} at 3 -1 \put{$\bullet$} at 3 2
\plot 1 2 3 2 3 -1 2 -1 /
\endpicture}
+
{\beginpicture
\setcoordinatesystem units <0.35cm,0.12cm> 
\setplotarea x from 1 to 3, y from -1 to 2   
\linethickness=0.5pt                        
\put{$\bullet$} at 1 -1 \put{$\bullet$} at 1 2
\put{$\bullet$} at 2 -1 \put{$\bullet$} at 2 2
\put{$\bullet$} at 3 -1 \put{$\bullet$} at 3 2
\plot 1 2 1 -1 3 -1 3 2 /
\endpicture}
+
{\beginpicture
\setcoordinatesystem units <0.35cm,0.12cm> 
\setplotarea x from 1 to 3, y from -1 to 2   
\linethickness=0.5pt                        
\put{$\bullet$} at 1 -1 \put{$\bullet$} at 1 2
\put{$\bullet$} at 2 -1 \put{$\bullet$} at 2 2
\put{$\bullet$} at 3 -1 \put{$\bullet$} at 3 2
\plot 1 -1 1 2  3 2 3 -1 /
\endpicture}
-
{\beginpicture
\setcoordinatesystem units <0.35cm,0.12cm> 
\setplotarea x from 1 to 3, y from -1 to 2   
\linethickness=0.5pt                        
\put{$\bullet$} at 1 -1 \put{$\bullet$} at 1 2
\put{$\bullet$} at 2 -1 \put{$\bullet$} at 2 2
\put{$\bullet$} at 3 -1 \put{$\bullet$} at 3 2
\plot 1 2 3 2  3 -1 /
\plot 1 -1 2 -1 /
\endpicture}
-
{\beginpicture
\setcoordinatesystem units <0.35cm,0.12cm> 
\setplotarea x from 1 to 3, y from -1 to 2   
\linethickness=0.5pt                        
\put{$\bullet$} at 1 -1 \put{$\bullet$} at 1 2
\put{$\bullet$} at 2 -1 \put{$\bullet$} at 2 2
\put{$\bullet$} at 3 -1 \put{$\bullet$} at 3 2
\plot 1 -1 3 -1  3 2 /
\plot 1 2 2 2 /
\endpicture}
-
{\beginpicture
\setcoordinatesystem units <0.35cm,0.12cm> 
\setplotarea x from 1 to 3, y from -1 to 2   
\linethickness=0.5pt                        
\put{$\bullet$} at 1 -1 \put{$\bullet$} at 1 2
\put{$\bullet$} at 2 -1 \put{$\bullet$} at 2 2
\put{$\bullet$} at 3 -1 \put{$\bullet$} at 3 2
\plot 1 -1 2 2  3 2 3 -1 /
\plot 1 2 2 -1 /
\endpicture}
-
{\beginpicture
\setcoordinatesystem units <0.35cm,0.12cm> 
\setplotarea x from 1 to 3, y from -1 to 2   
\linethickness=0.5pt                        
\put{$\bullet$} at 1 -1 \put{$\bullet$} at 1 2
\put{$\bullet$} at 2 -1 \put{$\bullet$} at 2 2
\put{$\bullet$} at 3 -1 \put{$\bullet$} at 3 2
\plot 1 2 2 -1  3 -1 3 2 /
\plot 1 -1 2 2 /
\endpicture}
\cr & \qquad
+
{\beginpicture
\setcoordinatesystem units <0.35cm,0.12cm> 
\setplotarea x from 1 to 3, y from -1 to 2   
\linethickness=0.5pt                        
\put{$\bullet$} at 1 -1 \put{$\bullet$} at 1 2
\put{$\bullet$} at 2 -1 \put{$\bullet$} at 2 2
\put{$\bullet$} at 3 -1 \put{$\bullet$} at 3 2
\plot 1 -1 1 2  2 2 /
\plot 2 -1 3 -1 3 2 /
\endpicture}
+
{\beginpicture
\setcoordinatesystem units <0.35cm,0.12cm> 
\setplotarea x from 1 to 3, y from -1 to 2   
\linethickness=0.5pt                        
\put{$\bullet$} at 1 -1 \put{$\bullet$} at 1 2
\put{$\bullet$} at 2 -1 \put{$\bullet$} at 2 2
\put{$\bullet$} at 3 -1 \put{$\bullet$} at 3 2
\plot 1 2 1 -1  2 -1 /
\plot 2 2 3 2 3 -1 /
\endpicture}
+
{\beginpicture
\setcoordinatesystem units <0.35cm,0.12cm> 
\setplotarea x from 1 to 3, y from -1 to 2   
\linethickness=0.5pt                        
\put{$\bullet$} at 1 -1 \put{$\bullet$} at 1 2
\put{$\bullet$} at 2 -1 \put{$\bullet$} at 2 2
\put{$\bullet$} at 3 -1 \put{$\bullet$} at 3 2
\plot 1 -1 2 -1  2 2 /
\plot 3 2 3 -1 /
\setquadratic
\plot 1 2 2 1 3 2 /
\endpicture}
+
{\beginpicture
\setcoordinatesystem units <0.35cm,0.12cm> 
\setplotarea x from 1 to 3, y from -1 to 2   
\linethickness=0.5pt                        
\put{$\bullet$} at 1 -1 \put{$\bullet$} at 1 2
\put{$\bullet$} at 2 -1 \put{$\bullet$} at 2 2
\put{$\bullet$} at 3 -1 \put{$\bullet$} at 3 2
\plot 1 2 2 2  2 -1 /
\plot 3 -1 3 2 /
\setquadratic
\plot 1 -1 2 0 3 -1 /
\endpicture}
- n
{\beginpicture
\setcoordinatesystem units <0.35cm,0.12cm> 
\setplotarea x from 1 to 3, y from -1 to 2   
\linethickness=0.5pt                        
\put{$\bullet$} at 1 -1 \put{$\bullet$} at 1 2
\put{$\bullet$} at 2 -1 \put{$\bullet$} at 2 2
\put{$\bullet$} at 3 -1 \put{$\bullet$} at 3 2
\plot 1 2 3 2  3 -1 1 -1 1 2 /
\plot 2 -1 2 2 /
\endpicture}
}$$

Part (a) of the following theorem is well known.

\thm
\item{(a)} For $n \in \ZZ_{\ge 0}$, $\kappa_n$ is a
central element of $\CC S_n$.
If $\lambda$ is a partition with $n$ boxes and $S_n^\lambda$ is
the irreducible $S_n$-module indexed by the partition $\lambda$,
$$
\kappa_n = \sum_{b \in \lambda} c(b),
\qquad\hbox{as operators on $S_n^\lambda$}.
$$
\smallskip\noindent
\item{(b)} Let $n,k \in \ZZ_{\ge 0}$.  Then,
as operators on $V^{\otimes k}$, where $\dim(V)=n$,
$$
Z_k = \kappa_n - {n \choose 2}  + kn,
\qquad\hbox{and}\qquad
Z_{k+{1 \over 2}} = \kappa_{n-1} - {n \choose 2} + (k+1) n - 1.
$$
\item{(c)} Let $n \in \CC$, $k\in \ZZ_{\ge 0}$.  Then
$Z_k$ is a central element of $\CC A_k(n)$,
and, if $n\in \CC$ is such that $\CC A_k(n)$ is semisimple and
$\lambda \vdash n$ with $|\lambda_{>1}|\le k$ boxes then
$$
Z_k = kn - {n \choose 2} + \sum_{b \in \lambda} c(b),
\qquad\hbox{as operators on $A_k^\lambda$,}
$$
where $A_k^\lambda$ is the irreducible $\CC A_k(n)$-module indexed by
the partition $\lambda$.  Furthermore,
$Z_{k+{1\over2}}$ is a central element of $\CC A_{k+{1\over2}}(n)$,
and, if $n$ is such that $\CC A_{k+{1\over2}}(n)$ is semisimple and
$\lambda\vdash n$ is a partition with $|\lambda_{>1}|\le k$ boxes then
$$
Z_{k+{1\over2}} = kn+n-1 - {n \choose 2} + \sum_{b \in \lambda} c(b),
\qquad\hbox{as operators on $A_{k+{1\over2}}^\lambda$,}
$$
where $A_{k+{1\over2}}^\lambda$ is the
irreducible $\CC A_{k+{1\over2}}(n)$-module indexed by
the partition $\lambda$.
\endthm

\pf (a) The element $\kappa_n$ is the class sum corresponding to the
conjugacy class of transpositions and thus $\kappa_n$ is a central
element of $\CC S_n$.  The constant by which $\kappa_n$ acts on
$S_n^\lambda$ is computed in [Mac, Ch.\ 1 \S 7 Ex.\ 7].

\smallskip\noindent
(c) The first statement follows from parts (a) and (b) and Theorems 3.6 and
3.22 as follows.
By Theorem 3.6, $\CC A_k(n)\cong \End_{S_n}(V^{\otimes k})$ if
$n\ge 2k$.  Thus, by Theorem 3.22, if $n\ge 2k$ then
$Z_k$ acts on the irreducible $\CC A_k(n)$-module $A_k^\lambda(n)$
by the constant given in the statement.  This means that
$Z_k$ is a central element of $\CC A_k(n)$ for all $n\ge k$.
Thus, for $n\ge 2k$, $d Z_k = Z_k d$ for all diagrams $d\in A_k$.
Since the coefficients in $d Z_k$ (in terms of the basis of diagrams)
are polynomials in $n$, it follows that $d Z_k = Z_k d$ for all
$n\in \CC$.

If $n\in \CC$ is such that $\CC A_k(n)$ is semisimple let
$\chi^\lambda_{\CC A_k(n)}$ be the irreducible
characters.  Then $Z_k$ acts on $A_k^\lambda(n)$ by the constant
$\chi^\lambda_{\CC A_k(n)}(Z_k)/\dim(A_k^\lambda(n))$.  If $n\ge k$
this is the constant in the statement, and therefore it is a polynomial
in $n$, determined by its values for $n\ge 2k$.

The proof of the second statement is completely analogous using $\CC \Akp$,
$S_{n-1}$, and the second statement in part (b).

\medskip\noindent
(b) Let $s_{ii} = 1$ so that
$$2 \kappa_n+n = n+2\sum_{1\le i<j\le n} s_{ij}
=\sum_{i=1}^n s_{ii} +\sum_{1\le i<j\le n} (s_{ij}+s_{ji})
=\sum_{i=j} s_{ij} +\sum_{i\ne j} s_{ij}
=\sum_{i,j=1}^n s_{ij}.$$
Then
$$\eqalign{
(2 \kappa_n+n)&(v_{i_1}\otimes \cdots \otimes v_{i_k})
=
\left(\sum_{i,j=1}^n s_{ij}\right)
(v_{i_1}\otimes \cdots \otimes v_{i_k})
=
\sum_{i,j=1}^n s_{ij}v_{i_1}\otimes \cdots \otimes s_{ij}v_{i_k} \cr
&=
\sum_{i,j=1}^n (1-E_{ii}-E_{jj}+E_{ij}+E_{ji})v_{i_1}\otimes
\cdots \otimes (1-E_{ii}-E_{jj}+E_{ij}+E_{ji})v_{i_k} \cr
}
$$
and expanding this sum gives that
$(2 \kappa_n+n)(v_{i_1}\otimes \cdots \otimes v_{i_k})$ is equal to
$$\eqalign{
\sum_{S\subseteq \{1,\ldots, k\}}
&\sum_{i_{1'},\ldots, i_{k'}}
\sum_{i,j=1}^n
\left(\prod_{\ell\in S^c} \delta_{i_\ell i_{\ell'}} \right)\cdot \cr
&\cdot\sum_{I\subseteq S\cup S'}
(-1)^{\#(\{\ell,\ell'\}\subseteq I)+\#(\{\ell,\ell'\}\subseteq I^c)}
\left(\prod_{\ell\in I} \delta_{i_\ell i} \right)
\left(\prod_{\ell\in I^c} \delta_{i_\ell j}\right)
(v_{i_{1'}}\otimes \cdots \otimes v_{i_{k'}}) \cr
}
\formula$$
where $S^c \subseteq \{1, \ldots, k\}$ corresponds to the tensor positions
where
1 is acting, and  where $I \subseteq S \cup S'$ corresponds to the tensor
positions that
must equal $i$ and $I^c$ corresponds to the tensor positions that must
equal $j$.

When $|S|=0$ the set $I$ is empty and the term corresponding to $S$ in
(3.36) is
$$
\sum_{i,j=1}^n
\sum_{i_{1'},\ldots, i_{k'}}
\big(\prod_{\ell\in \{1,\ldots,k\}} \delta_{i_\ell i_{\ell'}} \big)
(v_{i_{1'}}\otimes\cdots\otimes v_{i_{k'}})
= n^2 (v_{i_1}\otimes\cdots\otimes v_{i_k}).
$$
Assume $|S|\ge 1$ and separate the sum according to
the cardinality of $I$.  Note that the sum for $I$
is equal to the sum for $I^c$, since the whole sum
is symmetric in $i$ and $j$.
The sum of the terms in (3.36) which come from
$I=S\cup S'$ is equal to
$$
\sum_{i_{1'},\ldots, i_{k'}}
n \sum_{i=1}^n
\left(\prod_{\ell\in S^c} \delta_{i_\ell i_{\ell'}} \right)
(-1)^{|S|}
\left(\prod_{\ell\in S\cup S'} \delta_{i_\ell i} \right)
(v_{i_{1'}}\otimes\cdots\otimes v_{i_{k'}})
= n (-1)^{|S|}b_S (v_{i_1}\otimes \cdots \otimes v_{i_k}).
$$
We get a similar contribution from the sum of the terms
with $I = \emptyset$.

If $|S|>1$ then the sum of the terms in (3.36) which come from
$I = \{ \ell ,\ell'\}$ is equal to
$$\eqalign{
\sum_{i_{1'},\ldots, i_{k'}}
\sum_{i,j=1}^n
&\left(\prod_{r\in S^c} \delta_{i_r i_{r'}} \right)
(-1)^{|S|}
\delta_{i_\ell i}\delta_{i_{\ell'} i}
\left(\prod_{r\ne \ell} \delta_{i_r j}\delta_{i_{r'} j}\right)
(v_{i_{1'}}\otimes\cdots\otimes v_{i_{k'}})  \cr
&=
(-1)^{|S|} b_{S-\{\ell\}} (v_{i_1}\otimes \cdots \otimes v_{i_k}). \cr
}$$
and there is a corresponding contribution from $I = \{\ell,\ell'\}^c$.
The remaining terms can be written as
$$\eqalign{
\sum_{i_{1'},\ldots, i_{k'}}
\sum_{i,j=1}^n
\big(\prod_{\ell\in S^c} \delta_{i_\ell i_{\ell'}} \big)
\sum_{I\subseteq S\cup S'}
&(-1)^{\#(\{\ell,\ell'\}\subseteq I)+\#(\{\ell,\ell'\}\subseteq I^c)}
\big(\prod_{\ell\in I} \delta_{i_\ell i} \big)
\big(\prod_{\ell\in I^c} \delta_{i_\ell j}\big)
(v_{i_{1'}}\otimes v_{i_{k'}})  \cr
&\qquad=
2p_S (v_{i_1}\otimes \cdots\otimes v_{i_k}). \cr
}$$
Putting these cases together gives that
$2 \kappa_n+n$ acts on $v_{i_1}\otimes \cdots \otimes v_{i_k}$
the same way that
$$\eqalign{
\sum_{|S|=0} n^2
+\sum_{|S|=1} \big(2n(-1)^1b_S + 2p_S\big)
&+\sum_{|S|=2}
\Big(2n(-1)^2 b_S + 2p_S + \sum_{\ell\in S} (-1)^2 2b_{S-\{\ell\}}
\Big) \cr
&+\sum_{|S|>2}
\Big(2n(-1)^{|S|} b_S
+2p_S+\sum_{\ell\in S} (-1)^{|S|} 2b_{S-\{\ell\}}
\Big) \cr
}
$$
acts on $v_{i_1}\otimes \cdots \otimes v_{i_k}$.
Note that $b_S = 1$ if $|S|=1$.
Hence
$2 \kappa_n+n$ acts on $v_{i_1}\otimes \cdots \otimes v_{i_k}$
the same way that
$$\eqalign{
n^2
&+ \sum_{|S|=1} (-2n+2p_S)
+\sum_{|S|=2} (2nb_S + 2 + 2p_S)
+\sum_{|S|>2} \Big( (-1)^{|S|} 2n b_S + 2p_S
+\sum_{\ell\in S} (-1)^{|S|} 2b_{S-\{\ell\}} \Big) \cr
&= n^2-2nk+2{k\choose 2}
+\sum_{|S|\ge 1} 2p_S
+\sum_{|S|\ge 2} 2(n-k+|S|)(-1)^{|S|}b_S \cr
}
$$
acts on $v_{i_1}\otimes \cdots \otimes v_{i_k}$, and so $Z_k = \kappa_n
+(n-n^2+2nk)/2$ as operators on $V^{\otimes k}$.
This proves the first statement.

For the second statement, since
$(1-\delta_{in})(1-\delta_{jn})=\cases{ 0, &if $i=n$ or $j=n$,\cr
1, &otherwise,}$
$$\eqalign{
(2&\kappa_{n-1}+(n-1))(v_{i_1}\otimes\cdots\otimes v_{i_k}\otimes v_n)
=\left(\sum_{i,j=1}^{n-1} s_{ij}\right)
(v_{i_1}\otimes\cdots\otimes v_{i_k}\otimes v_n) \cr
&=\left(\sum_{i,j=1}^n s_{ij}(1-\delta_{in})(1-\delta_{jn})\right)
(v_{i_1}\otimes\cdots\otimes v_{i_k}\otimes v_n) \cr
&=\sum_{i,j=1}^n
s_{ij}v_{i_1}\otimes \cdots \otimes s_{ij}v_{i_k}\otimes
(1-\delta_{in})(1-\delta_{jn})v_n, \cr
&=\sum_{i,j=1}^n
(1-E_{ii}-E_{jj}+E_{ij}+E_{ji})v_{i_1}\otimes \cdots \otimes
(1-E_{ii}-E_{jj}+E_{ij}+E_{ji})v_{i_k} \cr
&\hskip3in \otimes
(1-E_{ii}-E_{jj}+E_{ii}E_{jj})v_n \cr
&=
\big(\sum_{i,j} s_{ij}\big)(v_{i_1}\otimes \cdots\otimes v_{i_k})
\otimes v_n \cr
&\qquad+\sum_{i,j=1}^n
(1-E_{ii}-E_{jj}+E_{ij}+E_{ji})v_{i_1}\otimes \cdots \otimes
(1-E_{ii}-E_{jj}+E_{ij}+E_{ji})v_{i_k}\otimes
(-E_{ii}-E_{jj})v_n \cr
&\qquad+\sum_{i,j=1}^n
(1-E_{ii}-E_{jj}+E_{ij}+E_{ji})v_{i_1}\otimes \cdots \otimes
(1-E_{ii}-E_{jj}+E_{ij}+E_{ji})v_{i_k}\otimes
E_{ii}E_{jj}v_n \cr
}$$
The first sum is known to equal $(2 \kappa_n + n)(v_{i_1} \otimes \cdots
\otimes v_{i_k})$
by the computation proving the first statement,
and the last sum has only one nonzero term, the term corresponding
to $i=j=n$.  Expanding the middle sum gives
$$\eqalign{
\sum_{S\subseteq \{1,\ldots, k+1\}\atop k+1\in S}
&\sum_{i_{1'},\ldots, i_{k'}}
\sum_{i,j=1}^n
\Big(\prod_{\ell\in S} \delta_{i_\ell, i_{\ell'}}\Big) \cdot \cr
&\cdot
\sum_I
(-1)^{\#(\{\ell,\ell'\}\subseteq I)+\#(\{\ell,\ell'\}\subseteq I^c)}
\left(\prod_{\ell\in I} \delta_{i_\ell i} \right)
\left(\prod_{\ell\in I^c} \delta_{i_\ell j}\right)
(v_{i_{1'}}\otimes \cdots \otimes v_{i_{k'}})  \cr
}
$$
where the inner sum is over all $I\subseteq \{1,\ldots, k+1\}$
such that $\{k+1,(k+1)'\}\subseteq I$ or $\{k+1,(k+1)'\}\subseteq I^c$.
As in part (a) this sum is treated in four cases:
(1) when $|S|=0$, (2) when $I = S\cup S'$ or $I=\emptyset$,
(3) when $I = \{\ell,\ell'\}$ or $I= \{\ell,\ell'\}^c$,
and (4) the remaining cases.
Since $k+1\in S$, the first case does not occur,
and cases (2), (3) and (4) are as in
part (a) giving
$$\sum_{|S|=1\atop k+1\in S} -2n
+ \sum_{|S|=2\atop k+1\in S} (2n b_S + 2p_S + 2)
+\sum_{|S|>2\atop k+1\in S}
\Big( 2n(-1)^{|S|}b_S + 2p_S +
2\sum_{\ell\in S} (-1)^{|S|}b_{S-\{\ell\}}\Big).$$
Combining this with the terms
$(2\kappa_n+n)(v_{i_1}\otimes \cdots\otimes v_{i_k})\otimes v_n$
and
$1\otimes(v_{i_1}\otimes \cdots\otimes v_{i_k}\otimes v_n)$
gives that $2\kappa_{n-1}+(n-1)$ acts on
$v_{i_1}\otimes \cdots\otimes v_{i_k}$
as
$$
(2\kappa_n+n)+1-2n+2k+\sum_{|S|\ge 2\atop k+1\in S}
2\tilde p_S +2(n-(k+1)+|S|)(-1)^{|S|}b_S.$$
Thus $\kappa_{n-1}-\kappa_n$ acts on $v_{i_1}\otimes \cdots\otimes v_{i_k}$
as
$$
{1\over2}\left(n-(n-1)+1-2n+2k+\sum_{|S|\ge 2\atop k+1\in S}
2\tilde p_S +2(n-(k+1)+|S|)(-1)^{|S|}b_S\right),
$$
so, as  operators on $V^{\otimes k}$, we have $Z_{k+ {1\over 2}} = k + Z_k +
(\kappa_{n-1} - \kappa_n) -1 + n - k = Z_k + (\kappa_{n-1} - \kappa_n) + n 
- 1$.
By the first statement in part (c) of this theorem we get $Z_{k+
{1\over 2}} = (\kappa_n - {n
\choose 2} + kn) + (\kappa_{n-1} - \kappa_n) + n - 1
= \kappa_{n-1}    - {n \choose 2} + kn  + n - 1
$.
\pfend

\thm
Let $k\in {1\over2}\ZZ_{\ge 0}$ and let $n\in \CC$.
\smallskip\noindent
(a) The elements
$M_{1\over2}, M_1,\ldots, M_{k-{1\over2}}, M_k$,
all commute with each other in $\CC A_k(n)$.
\smallskip\noindent
(b) Assume that $\CC A_k(n)$ is semisimple.
Let $\mu\in \hat A_k$ so that $\mu$ is a partition with $\le k$ boxes,
and let $A_k^\mu(n)$ be the irreducible $\CC A_k(n)$-module
indexed by $\mu$.  Then there is a unique, up to
multiplication by constants, basis
$\{ v_T\ |\ T\in \hat A_k^\mu\}$
of $A_k^\mu(n)$ such that, for all
$T = (T^{(0)},T^{({1 \over 2})},\ldots,T^{(k)})\in \hat A_k^\mu$, and
$\ell\in \ZZ_{\ge 0}$ such that $\ell\le k$,
$$
M_\ell v_T =
\cases{
c( T^{(\ell)}/ T^{ (\ell-{1\over2})}) v_T, &if
$T^{(\ell)}/T^{(\ell-{1\over2})} = \square$, \cr
(n-|T^{(\ell)}|)v_T, &if $T^{(\ell)}=T^{(\ell-{1\over2})}$, \cr
}
$$
and
$$
M_{\ell+{1\over2}} v_T =
\cases{
\big(n-c( T^{(\ell)}/ T^{(\ell+{1\over2})})\big) v_T, &if
$T^{(\ell)}/T^{(\ell+{1\over2})} = \square$, \cr
|T^{(\ell)}|v_T, &if $T^{(\ell)}=T^{(\ell+{1\over2})}$, \cr
}
$$
where $\lambda/\mu$ denotes the box where $\lambda$ and $\mu$ differ.

\pf (a) View $Z_0, Z_{1 \over 2}, \ldots, Z_k \in \CC A_k$. Then $Z_k
\in Z(\CC A_k)$, so $Z_k Z_\ell = Z_\ell Z_k$ for all $0 \le \ell \le k$.
Since $M_\ell = Z_\ell - Z_{\ell - {1 \over 2}}$, we see that the
$M_\ell$ commute with each other in $\CC A_k$.

(b) The basis is defined inductively. If $k = 0, {1\over 2}$ or $1$, then
$\dim(A_k^\lambda(n)) = 1$, so up to a constant there is a unique choice 
for the basis.
For $k > 1$, we consider the restriction
$\Res^{\CC A_k(n)}_{\CC A_{k-{1 \over 2}}(n)}
(A_k^\lambda(n))$. The branching rules for this restriction are
multiplicity free,
meaning that each $\CC A_{k-{1 \over 2}}(n)$-irreducible that shows up in
$A_k^\lambda(n)$ does so exactly once.  By induction, we can choose a
basis for each
$\CC A_{k-{1 \over 2}}(n)$-irreducible, and the union of these bases
forms a basis for
$A_k^\lambda(n)$.  For $\ell < k$,  $M_\ell  \in \CC A_{k-{1 \over 2}}(n)$,
so $M_\ell$
acts on this basis as in the statment of the theorem. It remains only
to check the
statement for $M_k$.  Let $k$ be an integer, and let $\lambda \vdash
n$ and $\gamma \vdash (n-1)$ such that $\lambda_{>1} = T^{(k)}$ and
$\gamma_{>1} = T^{(k-{1 \over 2})}$.  Then by Theorem 3.35(c),
$M_k =
Z_k - Z_{k - {1 \over 2}}$ acts on $v_T$ by the constant,
$$
\eqalign{
&\left(\sum_{b \in \lambda} c(b)-  {n \choose 2}  + kn \right)
-
\left(\sum_{b \in \gamma} c(b) - {n \choose 2}  + kn-1 \right)
= c(\lambda/\gamma)+1, \cr
}
$$
and
$M_{k+{1 \over 2}} = Z_{k+{1 \over 2}} - Z_k$
acts on $v_T \in A_{k + {1 \over 2}}^\lambda(n)$ by the constant
$$
\eqalign{
&\left(\sum_{b \in \gamma} c(b)-  {n \choose 2} + kn + n -1 \right)
-
\left(\sum_{b \in \lambda} c(b) - {n \choose 2} + kn \right)
= - c(\lambda/\gamma) + n  - 1. \cr
}
$$
The result now follows from (3.23) and the observation that
$$\hskip1.25truein
c(\lambda/\gamma)
= \cases{
c(T^{(k)}/T^{(k-{1\over2})})-1, & if $T^{(k)} = T^{(k+{1\over 2})} +
\square$, \cr
n - |T^{(k)}| -1, & if $T^{(k)} = T^{(k+{1\over 2})}$.
} \hskip1.25truein\qed
$$

\section 4. The Basic Construction

\bigskip
In this section we shall assume that all algebras are finite dimensional
algebras over an algebraically closed field $\FF$.  The fact that $\FF$ is
algebraically closed is only for convenience, to avoid the
division rings that could arise in the decomposition of
$\bar A$ just before (4.8) below.

\smallskip
Let $A \subseteq B$ be an inclusion of algebras. Then $B \otimes_{\FF} B$
is an
$(A,A)$-bimodule where $A$ acts on the left by left multiplication and on
the
right by right multiplication.  Fix an
$(A,A)$-bimodule homomorphism
$$
\varepsilon: B \otimes_\FF B \longrightarrow A.
\formula
$$
The {\it basic construction\/} is the algebra $B \otimes_A B$ with
product given by
$$
(b_1 \otimes b_2)(b_3 \otimes b_4)
= b_1 \otimes \varepsilon(b_2 \otimes b_3) b_4,
\qquad\hbox{for } b_1, b_2, b_3, b_4 \in B.
\formula
$$
More generally, let $A$ be an algebra and let $L$ be a left $A$-module and
$R$ a right
$A$-module. Let
$$
\varepsilon: L \otimes_\FF R \longrightarrow A,
\formula
$$
be an $(A,A)$-bimodule homomorphism. The {\it basic construction\/}
is the algebra $R \otimes_A L$ with product given by
$$
(r_1 \otimes \ell_1)(r_2 \otimes \ell_2)
= r_1\otimes\varepsilon(\ell_1 \otimes r_2)\ell_2,
\qquad\hbox{for } r_1, r_2 \in R \hbox{ and } \ell_1, \ell_2 \in L.
\formula
$$
Theorem 4.18 below determines, explicitly, the structure of the algebra
$R\otimes_A L$.

Let $N = \Rad(A)$ and let
$$
\bar A = A/N, \qquad \bar L = L/NL, \qquad\hbox{and}\qquad \bar R = R/RN
\formula
$$
Define an $(\bar A,\bar A)$-bimodule homomorphism
$$
\matrix{
\bar\varepsilon: & \bar L \otimes_\FF\bar R & \longrightarrow & \bar A\cr\cr
& \bar \ell \otimes \bar r & \mapsto & \overline{\varepsilon(\ell\otimes
r)}
\cr }
\formula
$$
where $\bar \ell = \ell + NL$, $\bar r = r + RN$, and $\bar a = a +
N$, for $\ell \in L, r \in R$ and $a \in A$.  Then by basic tensor product
relations [Bou1, Ch.\ II \S 3.3 Cor.\ to Prop.\ 2 and \S 3.6 Cor.\ to
Prop.\ 6], the surjective algebra homomorphism
$$
\matrix{
\pi: & R \otimes_A L & \longrightarrow & \bar R \otimes_{\bar A} \bar
L\cr\cr
& r  \otimes \ell & \mapsto & \bar r \otimes
\bar\ell\cr }
\qquad\hbox{has}\qquad
\ker (\pi) = R \otimes_A NL.
\formula
$$

The algebra $\bar A$ is a split semisimple algebra (an algebra isomorphic to
a direct sum of matrix algebras).   Fix an algebra isomorphism
$$
\matrix{
\bar A  & \mapright{\sim}
& \displaystyle{\bigoplus_{\mu \in \hat A}  M_{d_\mu}(\FF) }\cr \cr
a_{PQ}^\mu & \leftarrow & E_{PQ}^\mu \cr
}
$$
where $\hat A$ is an index set for the components and 
$E_{PQ}^\mu$ is the matrix with 1 in the $(P,Q)$ entry of the $\mu$th
block and 0 in all other entries.  Also, fix isomorphisms
$$
\bar L \cong
\bigoplus_{\mu \in \hat A} \overrightarrow{A}^\mu \otimes L^\mu
\qquad \hbox{ and } \qquad
\bar R \cong
\bigoplus_{\mu \in \hat A} R^\mu \otimes \overleftarrow A^\mu
\formula
$$
where $\overrightarrow{A}^\mu$, $\mu \in \hat A$, are the simple left $\bar
A$-modules,  $\overleftarrow A^\mu$, $\mu \in \hat A$, are the simple right
$\bar A$-modules, and $L^\mu$, $R^\mu$, $\mu \in \hat A$ are vector
spaces.
The practical effect of this setup is that if
$\hat R^\mu$ is an index set for a basis $\{r_Y^\mu | Y \in \hat R^\mu\}$
of $R^\mu$,  $\hat L^\mu$ is an index set for a basis
$\{\ell_X^\mu | X \in \hat L^\mu\}$ of $L^\mu$, and $\hat A^\mu$ is
an index set for bases
$$
\{\overrightarrow{a}^\mu_Q\ |\ Q \in \hat A^\mu\}
\hbox{ of $\overrightarrow{A}^\mu$ }
\qquad\hbox{ and }\qquad
\{\overleftarrow{a}^\mu_P\ |\ P \in \hat A^\mu\}
\hbox{ of $\overleftarrow{A}^\mu$ }
\formula
$$
such that
$$
a_{ST}^\lambda \overrightarrow{a}_Q^\mu
= \delta_{\lambda\mu} \delta_{TQ} \overrightarrow{a}_S^\mu
\qquad\hbox{ and }\qquad
\overleftarrow{a}_P^\mu a_{ST}^\lambda
= \delta_{\lambda\mu} \delta_{PS} \overleftarrow{a}_T^\mu
\formula
$$
then
$$\matrix{
\hbox{
$\bar L$ has basis
$\{ \overrightarrow{a}_P^\mu \otimes \ell_X^\mu\ |\
\mu \in \hat A, P \in \hat A^\mu, X \in \hat L^\mu  \}$ }
\qquad\hbox{and}\hfill \cr
\cr
\hbox{ $\bar R$ has basis $\{ r_Y^\mu \otimes \overleftarrow{a}_Q^\mu\ |\
\mu \in \hat A, Q \in \hat A^\mu, Y \in \hat R^\mu  \}$.} \hfill \cr}
\formula
$$

With notations as in (4.9) and (4.11) the map
$\bar \varepsilon: \bar L \otimes_\FF \bar R \to \bar A$
is determined by the constants $\varepsilon_{XY}^\mu \in \FF$
given by
$$
\varepsilon(\overrightarrow{a}_Q^\mu \otimes \ell_X^\mu
\otimes r_Y^\mu \otimes \overleftarrow{a}_P^\mu)
= \varepsilon_{XY}^\mu a_{QP}^\mu
\formula
$$
and $\varepsilon_{XY}^\mu$ does not depend on $Q$ and $P$ since
$$\eqalign{
\varepsilon(\overrightarrow{a}_S^\lambda \otimes \ell_X^\lambda
\otimes r_Y^\mu \otimes \overleftarrow{a}_T^\mu)
&=
\varepsilon(a_{SQ}^\lambda\overrightarrow{a}_Q^\lambda \otimes
\ell_X^\lambda
\otimes r_Y^\mu \otimes \overleftarrow{a}_P^\mu a_{PT}^\mu)  \cr
&= a_{SQ}^\lambda
\varepsilon(\overrightarrow{a}_Q^\lambda \otimes \ell_X^\lambda
\otimes r_Y^\mu \otimes \overleftarrow{a}_P^\mu) a_{PT}^\mu   \cr
&= \delta_{\lambda\mu} a_{SQ}^\mu
\varepsilon_{XY}^\mu a_{QP}^\mu a_{PT}^\mu
= \varepsilon_{XY}^\mu a_{ST}^\mu.    \cr
}\formula$$
For each $\mu \in \hat A$ construct a matrix
$$
{\cal E}^\mu = (\varepsilon_{XY}^\mu)
\formula
$$
and let $D^\mu = (D_{ST}^\mu)$ and $C^\mu = (C_{ZW}^\mu)$ be invertible
matrices such that $D^\mu {\cal E}^\mu C^\mu$ is a diagonal matrix with
diagonal entries denoted $\varepsilon_X^\mu$,
$$
D^\mu {\cal E}^\mu C^\mu = {\rm diag}(\varepsilon_X^\mu).
\formula
$$
In practice $D^\mu$ and $C^\mu$ are found by row reducing $\cal E^\mu$ to
its Smith normal form. The $\varepsilon_P^\mu$ are the
{\it invariant factors} of $\cal E^\mu$.

For $\mu \in \hat A, X \in \hat R^\mu, Y \in \hat L^\mu$, define the
following elements of $\bar R \otimes_{\bar A} \bar L$,
$$
\bar m_{XY}^\mu
= r_X^\mu \otimes \overrightarrow{a}_P^\mu \otimes \overleftarrow{a}_P^\mu
\otimes \ell_Y^\mu,
\qquad\hbox{and}\qquad
\bar{n}_{XY}^\mu
= \sum_{Q_1, Q_2} C_{Q_1 X}^\mu D_{Y Q_2}^\mu \bar m_{Q_1 Q_2}^\mu.
\formula
$$
Since
$$\eqalign{
(r_S^\lambda \otimes \overrightarrow{a}_W^\lambda
\otimes \overleftarrow{a}_Z^\mu \otimes \ell_T^\mu)
&=
(r_S^\lambda \otimes \overrightarrow{a}_P^\lambda a_{PW}^\lambda
\otimes \overleftarrow{a}_Z^\mu \otimes \ell_T^\mu) \cr
&=
(r_S^\lambda \otimes \overrightarrow{a}_P^\lambda
\otimes a_{PW}^\lambda\overleftarrow{a}_Z^\mu \otimes \ell_T^\mu) \cr
&=\delta_{\lambda\mu}\delta_{WZ}
(r_S^\lambda \otimes \overrightarrow{a}_P^\lambda
\otimes \overleftarrow{a}_P^\lambda \otimes \ell_T^\lambda) \cr
}
\formula$$
the element $\bar m_{XY}^\mu$ does not depend on $P$ and
$\{ \bar m_{XY}^\mu\ |\
\mu\in \hat A, X\in \hat R^\mu, Y\in \hat L^\mu\}$ is a basis of
$\bar R\otimes_{\bar A} \bar L$.

The following theorem determines the structure of the algebras
$R\otimes_A L$ and $\bar R\otimes_{\bar A}\bar L$.
This theorem is used by W.P.\ Brown in the study of the Brauer
algebra. Part (a) is implicit in [Bro1,\S2.2] and part (b) is
proved in [Bro2].

\thm Let $\pi\colon R\otimes_A L \to \bar R\otimes_{\bar A} \bar L$
be as in (4.7) and let $\{ k_i\}$ be a basis of
$\ker(\pi) = R \otimes_A NL$.  Let
$$n_{Y T}^\mu \in R \otimes_A L
\quad\hbox{be such that}\quad
\pi(n_{YT}^\mu) = \bar n_{Y T}^\mu,$$
where the elements $\bar n_{Y T}^\mu \in \bar R \otimes_{\bar A} \bar L$ are
as defined in (4.16).
\medskip
\item{(a)} The sets
$\{\bar m_{XY}^\mu\ |\ \mu \in \hat A, X \in \hat R^\mu, Y \in \hat
L^\mu\}$
and
$\{\bar n_{XY}^\mu\ |\ \mu \in \hat A, X \in \hat R^\mu, Y \in \hat
L^\mu\}$ (see (4.16))
are bases of $\bar R \otimes_{\bar A} \bar L$, which satisfy
$$
\bar m_{ST}^\lambda \bar m_{QP}^\mu
= \delta_{\lambda \mu} \varepsilon_{TQ}^\mu \bar m_{SP}^\mu
\qquad\hbox{ and }\qquad
\bar n_{ST}^\lambda \bar n_{QP}^\mu
= \delta_{\lambda \mu} \delta_{TQ} \varepsilon_{T}^\mu \bar n_{SP}^\mu,
$$
where $\varepsilon_{TQ}^\mu$ and $\varepsilon_{T}^\mu$ are as defined
in (4.12) and (4.15).
\medskip
\item{(b)} The radical of the algebra $R \otimes_A L$ is
$$
\Rad(R \otimes_A L)
= \FF\hbox{-span}\{ k_i, n_{Y T}^\mu\ |\ \varepsilon_Y^\mu = 0
\hbox{ or } \varepsilon_T^\mu = 0 \}
$$
and the images of the elements
$$e_{Y T}^\mu = {1 \over \varepsilon_T^\mu}  n_{Y T}^\mu,
\qquad\hbox{for $\varepsilon_Y^\mu \not= 0$ and $\varepsilon_T^\mu \not=0$,}
$$
are a set of matrix units in
$(R \otimes_A L)/\Rad(R \otimes_A L)$.
\endthm

\pf
The first statement in (a) follows from the equations in (4.17).
If $(C^{-1})^\mu$ and $(D^{-1})^\mu$ are the inverses of
the matrices $C^\mu$ and $D^\mu$ then
$$\eqalign{
\sum_{X,Y} (C^{-1})^\mu_{XS} (D^{-1})^\mu_{TY} \bar n_{XY}
&=
\sum_{X,Y,Q_1,Q_2} (C^{-1})^\mu_{XS} C_{Q_1X}^\mu
\bar m_{Q_1Q_2}D^\mu_{YQ_2}(D^{-1})^\mu_{TY} \cr
&=
\sum_{Q_1,Q_2} \delta_{SQ_1}\delta_{Q_2T}\bar m_{Q_1Q_2}^\mu
 = \bar m_{ST}^\mu, \cr
}$$
and so the elements $\bar m_{ST}^\mu$ can be written as linear
combinations of the $\bar n_{XY}^\mu$.  This establishes
the second statement in (a).  By direct computation, using (4.10) and
(4.12),
$$\eqalign{
\bar m_{ST}^\lambda \bar m_{QP}^\mu
&=
(r_S^\lambda \otimes \overrightarrow{a}_W^\lambda
\otimes \overleftarrow{a}_W^\lambda \otimes \ell_T^\lambda)
(r_Q^\mu \otimes \overrightarrow{a}_Z^\mu
\otimes \overleftarrow{a}_Z^\mu \otimes \ell_P^\mu) \cr
&=
r_S^\lambda \otimes \overrightarrow{a}_W^\lambda
\otimes
\varepsilon(\overleftarrow{a}_W^\lambda \otimes \ell_T^\lambda
\otimes
r_Q^\mu \otimes \overrightarrow{a}_Z^\mu)
\overleftarrow{a}_Z^\mu \otimes \ell_P^\mu \cr
&= \delta_{\lambda\mu}
(r_S^\lambda \otimes \overrightarrow{a}_W^\lambda
\otimes \varepsilon_{TQ}^\lambda \bar a_{WZ}^\lambda
\overleftarrow{a}_Z^\lambda \otimes \ell_P^\lambda) \cr
&= \delta_{\lambda\mu} \varepsilon_{TQ}^\lambda
(r_S^\lambda \otimes \overrightarrow{a}_W^\lambda
\otimes \overleftarrow{a}_W^\lambda \otimes \ell_P^\lambda)
= \delta_{\lambda\mu} \varepsilon_{TQ}^\lambda \bar m_{SP}^\lambda, \cr
}$$
and
$$\eqalign{
\bar n_{ST}^\lambda \bar n_{UV}^\mu
&= \sum_{Q_1,Q_2,Q_3,Q_4} C_{Q_1S}^\lambda D_{TQ_2}^\lambda
\bar m_{Q_1Q_2}^\lambda C_{Q_3U}^\mu D_{VQ_4}^\mu \bar m_{Q_3Q_4}^\mu \cr
&= \sum_{Q_1,Q_2,Q_3,Q_4} \delta_{\lambda\mu}
C_{Q_1S}^\lambda D_{TQ_2}^\lambda
\varepsilon_{Q_2Q_3}^\mu C_{Q_3U}^\mu D_{VQ_4}^\mu
\bar m_{Q_1Q_4}^\mu \cr
&= \delta_{\lambda\mu}
\sum_{Q_1,Q_4} \delta_{TU}\varepsilon_T^\mu
C_{Q_1S}^\mu D_{VQ_4}^\mu \bar m_{Q_1Q_4}^\mu
= \delta_{\lambda\mu}\delta_{TU}\varepsilon_T^\mu \bar n_{SV}^\mu. \cr
}$$

\smallskip\noindent
(b) Let $N=\Rad(A)$ as in (4.5).
If $r_1\otimes n_1\ell_1, r_2\otimes n_2\ell_2\in R\otimes_A NL$
with $n_1\in N^i$ for some $i\in \ZZ_{>0}$ then
$$(r_1\otimes n_1\ell_1)(r_2\otimes n_2\ell_2)
= r_1\otimes \varepsilon(n_1\ell_1\otimes r_2)n_2\ell_2
= r_1\otimes n_1\varepsilon(\ell_1\otimes r_2)n_2\ell_2
\in R\otimes_A N^{i+1}L.$$
Since $N$ is a nilpotent ideal
of $A$ it follows that $\ker(\pi) = R\otimes_A NL$ is a
nilpotent ideal of $R\otimes_A L$.
So $\ker(\pi)\subseteq \Rad(R\otimes_A L)$.

Let
$$I = \FF\hbox{-span}\{ k_i, n_{Y T}^\mu\ |\ \varepsilon_Y^\mu = 0
\hbox{ or } \varepsilon_T^\mu = 0 \}.$$
The multiplication rule for the $\bar n_{YT}$ implies that
$\pi(I)$ is an ideal of
$\bar R\otimes_{\bar A} \bar L$ and
thus, by the correspondence between ideals of
$\bar R\otimes_{\bar A} \bar L$ and ideals of
$R\otimes_A L$ which contain $\ker(\pi)$, $I$ is an
ideal of $R\otimes_A L$.

If $\bar n_{Y_1T_1}^\mu,\bar n_{Y_2T_2}^\mu, \bar n_{Y_3T_3}^\mu\in
\{ \bar n_{YT}^\mu
\ |\ \hbox{$\varepsilon_Y^\mu=0$ or $\varepsilon_T^\mu=0$}\}$
then
$$\bar n_{Y_1T_1}^\mu\bar n_{Y_2T_2}^\mu\bar n_{Y_3T_3}^\mu
= \delta_{T_1Y_2}\varepsilon_{Y_2}^\mu \bar n_{Y_1T_2}^\mu
\bar n_{Y_3T_3}^\mu
= \delta_{T_1Y_2}\delta_{T_2Y_3}
\varepsilon_{Y_2}^\mu\varepsilon_{T_2}^\mu \bar n_{Y_1 T_3}^\mu = 0,$$
since
$\varepsilon_{Y_2}^\mu=0$ or $\varepsilon_{T_2}^\mu=0$.
Thus any product
$n_{Y_1T_1}^\mu n_{Y_2T_2}^\mu n_{Y_3T_3}^\mu$ of three basis elements
of $I$ is in $\ker(\pi)$.  Since $\ker(\pi)$ is
a nilpotent ideal of $R\otimes_A L$ it follows that
$I$ is an ideal of $R\otimes_A L$ consisting of nilpotent elements.
So $I\subseteq \Rad(R\otimes_A L)$.

Since
$$e_{YT}^\lambda e_{UV}^\mu
= {1\over \varepsilon_T^\lambda} {1\over \varepsilon_V^\mu}
n_{YT}^\lambda n_{UV}^\mu
= \delta_{\lambda\mu}\delta_{TU}
{1\over \varepsilon_T^\lambda \varepsilon_V^\lambda}
\varepsilon_T^\lambda n_{YV}^\lambda
=
\delta_{\lambda\mu}\delta_{TU}e_{YV}^\lambda
\qquad\hbox{mod $I$},$$
the images of the elements $e_{YT}^\lambda$ in (4.7) form
a set of matrix units in the algebra $(R\otimes_A L)/I$.  Thus
$(R\otimes_A L)/I$ is a split semisimple algebra and
so $I \supseteq \Rad(R\otimes_A L)$.
\endpf

\bigskip\noindent
\eject
{\it Basic constructions for $A\subseteq B$}
\medskip

Let $A \subseteq B$ be an inclusion of algebras.
Let $\varepsilon_1: B\to A$ be an $(A,A)$ bimodule homomorphism
and use the $(A,A)$-bimodule homomorphism
$$\matrix{ \varepsilon : &B\otimes_\FF B &\longrightarrow & A \cr
&b_1\otimes b_2 &\longmapsto &\varepsilon_1(b_1b_2) \cr}
\formula
$$
and (4.2) to define the basic construction $B\otimes_A B$.
Theorem 4.28 below provides the structure of $B\otimes_A B$ in the
case that both $A$ and $B$ are split semisimple.

Let us record the following facts,
\global\advance\resultno by 1
\itemitem{(4.20a)}
If $p\in A$ and $pAp=\FF p$ then
$(p\otimes 1)(B\otimes_A B)(p\otimes 1)
=\FF \cdot (p\otimes 1)$,
\smallskip\noindent
\itemitem{(4.20b)}
If $p$ is an idempotent of $A$ and $pAp = \FF p$ then
$\varepsilon_1(1)\in \FF$,
\smallskip\noindent
\itemitem{(4.20c)}
If $p\in A$, $pAp = \FF p$ and
if $\varepsilon_1(1)\ne 0$, then 
${1\over\varepsilon(1)}(p\otimes 1)$
is a  minimal idempotent in $B\otimes_A B$,
\smallskip\noindent
which are justified as follows.
If $p\in A$ and $pAp=\FF p$ 
and $b_1,b_2\in B$ then
$(p\otimes 1)(b_1\otimes b_2)(p\otimes 1)
=(p\otimes \varepsilon_1(b_1)b_2)(p\otimes 1)
=p\otimes \varepsilon_1(b_1)\varepsilon_1(b_2p)
=p\varepsilon_1(b_1)\varepsilon_1(b_2)p\otimes 1
=\xi p\otimes 1$,
for some constant $\xi\in \FF$.  This establishes (a).
If $p$ is an idempotent of $A$ and $pAp = \FF p$ then
$p\varepsilon_1(1)p = \varepsilon_1(p^2)=
\varepsilon_1(1\cdot p)=\varepsilon_1(1)p$ and so (b) holds.
If $p\in A$ and $pAp = \FF p$ then
$(p\otimes 1)^2 = \varepsilon_1(1)(p\otimes 1)$
and so, if $\varepsilon_1(1)\ne 0$, then
${1\over\varepsilon(1)}(p\otimes 1)$
is a  minimal idempotent in $B\otimes_A B$.

Assume $A$ and $B$ are split semisimple.  Let
\smallskip\noindent
\itemitem{} $\hat A$ be an index set for the irreducible
$A$-modules $A^\mu$,
\smallskip\noindent
\itemitem{} $\hat B$ be an index set for the irreducible
$B$-modules $B^\lambda$,\quad and let
\smallskip\noindent
\itemitem{}
$\hat A^\mu= \{\  {\scriptstyle{P\to}} \mu\ \}$
be an index set for a basis of the simple $A$-module $A^\mu$,
\smallskip\noindent
for each $\mu\in \hat A$ (the composite ${\scriptstyle{P\to}}\mu$
is viewed as a single symbol).  We think of $\hat A^\mu$ as the set of 
``paths to $\mu$'' in the two level graph
$$\matrix{
\Gamma \qquad\hbox{with}\hfill \cr \cr}
\qquad
\matrix{
\hbox{vertices on level A: \quad $\hat A$, \qquad
vertices on level B: \quad $\hat B$,\qquad and} \hfill \cr
\hbox{$m_\mu^\lambda$ edges $\mu\to\lambda$ if $A^\mu$
appears with multiplicity $m_\mu^\lambda$ in $\Res^B_A(B^\lambda)$.}
\hfill\cr}
\formula
$$
For example, the graph $\Gamma$ for the   symmetric
group algebras $A =\CC S_3$ and $B = \CC S_4$ is
$$
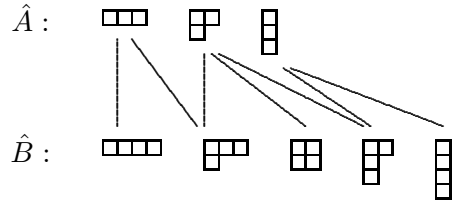
\matrix{
{\beginpicture
\setcoordinatesystem units <0.175cm,0.175cm>         
\setplotarea x from -2 to 27, y from -2 to -4   
\linethickness=0.5pt                          
\put{$\hat B:$} at  -2 -1.5
\putrectangle corners at 3 -2 and 4 -1
\putrectangle corners at 4 -2 and 5 -1
\putrectangle corners at 5 -2 and 6 -1
\putrectangle corners at 6 -2 and 7 -1
\putrectangle corners at 10 -2 and 11 -1
\putrectangle corners at 11 -2 and 12 -1
\putrectangle corners at 12 -2 and 13 -1
\putrectangle corners at 10 -3 and 11 -2
\putrectangle corners at 16 -2 and 17 -1
\putrectangle corners at 17 -2 and 18 -1
\putrectangle corners at 16 -3 and 17 -2
\putrectangle corners at 17 -3 and 18 -2
\putrectangle corners at 21 -2 and 22 -1
\putrectangle corners at 22 -2 and 23 -1
\putrectangle corners at 21 -3 and 22 -2
\putrectangle corners at 21 -4 and 22 -3
\putrectangle corners at 26 -2 and 27 -1
\putrectangle corners at 26 -3 and 27 -2
\putrectangle corners at 26 -4 and 27 -3
\putrectangle corners at 26 -5 and 27 -4
%
\plot 4 6 4 0 /
\plot 5 6 9.5 0 /
\plot 10 5 10 0 /
\plot 10.5 5 17 0 /
\plot 11 5 21 0 /
\plot 15.5 4 21.5 0 /
\plot 16 4 26.5 0 /
\put{$\hat A:$} at  -2 7.5
\putrectangle corners at 3 7 and 4 8
\putrectangle corners at 4 7 and 5 8
\putrectangle corners at 5 7 and 6 8
\putrectangle corners at 9 7 and 10 8
\putrectangle corners at 10 7 and 11 8
\putrectangle corners at 9 6 and 10 7
\putrectangle corners at 14 7 and 15 8
\putrectangle corners at 14 6 and 15 6
\putrectangle corners at 14 5 and 15 7
\endpicture} \cr}
$$
If $\lambda\in \hat B$ then
$$\hat B^\lambda = \{ {\scriptstyle{P\to}} \mu \to \lambda
\ |\
\hbox{$\mu\in \hat A$, ${\scriptstyle{P\to}} \mu \in \hat A^\mu$
and $\mu\to \lambda$ is an edge in $\Gamma$} \}
\formula
$$
is an index set for a basis of the irreducible $B$-module $B^\lambda$.
We think of $\hat B^\lambda$ as the set of paths to $\lambda$ in the 
graph $\Gamma$.
Let 
$$\{ a_\trisub{P}{Q}{\mu}
\ |\ \mu\in \hat A, {\scriptstyle{P\to}}\mu,
{\scriptstyle{Q\to}}\mu \in \hat A^\mu\}
\qquad\hbox{and}\qquad
\{ b_\dtrisub{P}{Q}{\mu}{\nu}{\lambda}
\ |\
\lambda\in \hat B,  {\scriptstyle{P\to}}\mu\to \lambda,
{\scriptstyle{Q\to}}\nu\to \lambda\in \hat B^\lambda\},
\formula
$$
be sets of matrix units in the algebras $A$ and $B$,
respectively, so that
$$
a_\trisub{P}{Q}{\mu} a_\trisub{S}{T}{\nu}
=\delta_{\mu\nu}\delta_{QS} a_\trisub{P}{T}{\mu}
\qquad\hbox{and}\qquad
b_\dtrisub{P}{Q}{\mu}{\gamma}{\lambda} b_\dtrisub{S}{T}{\tau}{\nu}{\sigma}
=\delta_{\lambda\sigma}\delta_{QS}\delta_{\gamma\tau}
b_\dtrisub{P}{T}{\mu}{\nu}{\lambda},
\formula
$$
and such that, for all $\mu\in \hat A$, $P,Q\in \hat A^\mu$,
$$
a^\mu_\trisub{P}{Q}{\mu}
= \sum_{\mu\to \lambda} b^\lambda_\dtrisub{P}{Q}{\mu}{\mu}{\lambda}
\formula
$$
where the sum is over all edges $\mu\to \lambda$ in the graph $\Gamma$.

Though is not necessary for the following it is conceptually
helpful to let $C= B\otimes_A B$, let  $\hat C = \hat A$ and
extend the graph $\Gamma$ to a graph $\hat \Gamma$ with
three levels, so that the edges between level
B and level C are the reflections of the
edges between level A and level B.  In other words,
$$\matrix{
\hbox{$\hat \Gamma$\qquad has}\cr \cr}
\qquad
\matrix{
\hbox{vertices on level $C$: \qquad $\hat C$,
\qquad and} \hfill \cr
\hbox{an edge $\lambda\to \mu$, $\lambda\in \hat B$, $\mu\in\hat C$,
for each edge $\mu\to \lambda$, $\mu\in \hat A$,
$\lambda\in \hat B$.} \hfill \cr}
\formula
$$
For each $\nu\in \hat C$ define
$$\hat C^\nu = \left\{ {\scriptstyle{P\to}}\mu\to\lambda\to\nu\
\Big|\
\matrix{
\hbox{$\mu\in \hat A$, $\lambda\in \hat B$, $\nu\in \hat C$,
${\scriptstyle{P\to}}\mu\in \hat A^\mu$ and} \cr
\hbox{$\mu\to \lambda$ and $\lambda\to \nu$ are edges in $\hat \Gamma$} \cr
}
\right\},
\formula
$$
so that $\hat C^\nu$ is the set of ``paths to $\nu$''
in the graph $\hat \Gamma$.  Continuing with our previous
example,  $\hat \Gamma$ is
$$
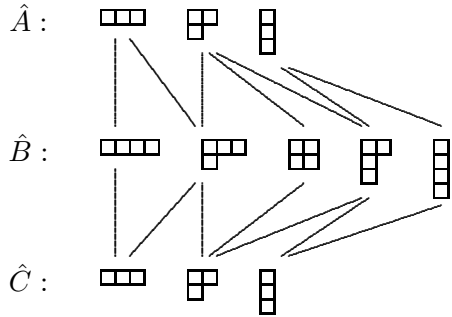
\matrix{
{\beginpicture
\setcoordinatesystem units <0.175cm,0.175cm>         
\setplotarea x from -2 to 27, y from -2 to -4   
\linethickness=0.5pt                          
\put{$\hat B:$} at  -2 -1.5
\putrectangle corners at 3 -2 and 4 -1
\putrectangle corners at 4 -2 and 5 -1
\putrectangle corners at 5 -2 and 6 -1
\putrectangle corners at 6 -2 and 7 -1
\putrectangle corners at 10 -2 and 11 -1
\putrectangle corners at 11 -2 and 12 -1
\putrectangle corners at 12 -2 and 13 -1
\putrectangle corners at 10 -3 and 11 -2
\putrectangle corners at 16 -2 and 17 -1
\putrectangle corners at 17 -2 and 18 -1
\putrectangle corners at 16 -3 and 17 -2
\putrectangle corners at 17 -3 and 18 -2
\putrectangle corners at 21 -2 and 22 -1
\putrectangle corners at 22 -2 and 23 -1
\putrectangle corners at 21 -3 and 22 -2
\putrectangle corners at 21 -4 and 22 -3
\putrectangle corners at 26 -2 and 27 -1
\putrectangle corners at 26 -3 and 27 -2
\putrectangle corners at 26 -4 and 27 -3
\putrectangle corners at 26 -5 and 27 -4
%
\plot 4 6 4 0 /
\plot 5 6 9.5 0 /
\plot 10 5 10 0 /
\plot 10.5 5 17 0 /
\plot 11 5 21 0 /
\plot 15.5 4 21.5 0 /
\plot 16 4 26.5 0 /
\put{$\hat A:$} at  -2 7.5
\putrectangle corners at 3 7 and 4 8
\putrectangle corners at 4 7 and 5 8
\putrectangle corners at 5 7 and 6 8
\putrectangle corners at 9 7 and 10 8
\putrectangle corners at 10 7 and 11 8
\putrectangle corners at 9 6 and 10 7
\putrectangle corners at 14 7 and 15 8
\putrectangle corners at 14 6 and 15 6
\putrectangle corners at 14 5 and 15 7
\put{$\hat C:$} at  -2 -10.5
\putrectangle corners at 3 -10 and 4 -11
\putrectangle corners at 4 -10 and 5 -11
\putrectangle corners at 5 -10 and 6 -11
\putrectangle corners at 9 -10 and 10 -11
\putrectangle corners at 10 -10 and 11 -11
\putrectangle corners at 9 -11 and 10 -12
\putrectangle corners at 14 -10 and 15 -11
\putrectangle corners at 14 -11 and 15 -12
\putrectangle corners at 14 -12 and 15 -13
\plot 4 -9 4 -3 /
\plot 5 -9 9.5 -4 /
\plot 10 -9 10 -4 /
\plot 10.5 -9 17 -4 /
\plot 11 -9 21 -5 /
\plot 15.5 -9 21.5 -5 /
\plot 16 -9 26.5 -5.5 /
\endpicture} \cr}
$$

\thm Assume $A$ and $B$ are split semisimple, and let the notations and
assumption be as in (4.21-4.25).
\smallskip\noindent
\item{(a)}  The elements of $B\otimes_A B$ given by
$$
b_\dtrisub{P}{T}{\mu}{\gamma}{\lambda}
\otimes
b_\dtrisub{T}{Q}{\gamma}{\nu}{\sigma}
$$
do not depend on the choice of
${\scriptstyle{T\to}}\gamma\in \hat A^\gamma$ and
form a basis of $B\otimes_A B$.
\item{(b)} For each edge $\mu\to\lambda$ in $\Gamma$ define a
constant $\varepsilon_\mu^\lambda\in \FF$ by
$$
\varepsilon_1\bigg(
b_\dtrisub{P}{P}{\mu}{\mu}{\lambda}\bigg)
=\varepsilon_\mu^\lambda\,\,
a_\trisub{P}{P}{\mu}
\formula
$$
Then $\varepsilon_\mu^\lambda$ is independent of
the choice of ${\scriptstyle{P\to}}\mu\in \hat A^\mu$ and
$$
\Big(
b_\dtrisub{P}{T}{\mu}{\gamma}{\lambda}
\otimes
b_\dtrisub{T}{Q}{\gamma}{\nu}{\sigma}
\Big)
\Big(
b_\dtrisub{R}{X}{\tau}{\pi}{\rho}
\otimes
b_\dtrisub{X}{S}{\pi}{\xi}{\eta}
\Big)
=\delta_{\gamma\pi}
\delta_{QR}\delta_{\nu\tau}
\delta_{\sigma\rho}
\varepsilon_{\gamma}^\sigma
\Big(
b_\dtrisub{P}{T}{\pi}{\mu}{\gamma}
\otimes
b_\dtrisub{T}{S}{\gamma}{\xi}{\eta}
\Big).
$$
$$\Rad(B\otimes_A B)\qquad\hbox{has basis}\qquad
\Big\{
b_\dtrisub{P}{T}{\mu}{\gamma}{\lambda}
\otimes
b_\dtrisub{T}{Q}{\gamma}{\nu}{\sigma}
\ |\
\hbox{$\varepsilon_\mu^\lambda=0$ or $\varepsilon_\nu^\sigma=0$}\Big\},
$$
and the images of the elements
$$
e_\ddtrisub{P}{Q}{\mu}{\nu}{\lambda}{\sigma}{\gamma}
=\left({1\over \varepsilon_\gamma^\sigma}\right)
\Big(
b_\dtrisub{P}{T}{\mu}{\gamma}{\lambda}
\otimes
b_\dtrisub{T}{Q}{\gamma}{\nu}{\sigma}
\Big),
\qquad\hbox{such that}\quad
\hbox{$\varepsilon_\mu^\lambda\ne 0$ and
$\varepsilon_\nu^\sigma\ne 0$},
$$
form a set of matrix units in $(B\otimes_A B)/\Rad(B\otimes_A B)$.
\smallskip\noindent
\item{(c)}
Let $tr_B: B \to \FF$ and
$tr_A: A \to \FF$ be traces on $B$ and $A$, respectively, such that
$$
tr_A(\varepsilon_1(b)) = tr_B(b), \qquad\hbox{ for all $b \in B$}.
\formula
$$
Let $\chi^\mu_A$, $\mu\in \hat A$, and
$\chi^\lambda_B$, $\lambda\in \hat B$, be the irreducible characters
of the algebras $A$ and $B$, respectively.  Define constants
$\tr_A^\mu$, $\mu\in \hat A$, and $\tr_B^\lambda$, $\lambda\in \hat B$,
by the equations
$$
\tr_A = \sum_{\mu\in \hat A} \tr_A^\mu \chi_A^\mu
\qquad\hbox{and}\qquad
\tr_B = \sum_{\lambda\in \hat B} \tr_B^\lambda \chi_B^\lambda,
\formula
$$
respectively.  Then the constants
$\varepsilon_\mu^\lambda$ defined in (4.29) satisfy
$$\tr_B^\lambda = \varepsilon_\mu^\lambda\,\, \tr_A^\mu.
$$
\item{(d)}
In the algebra $B\otimes_A B$,
$$
1\otimes 1
=
\sum_\starsub{P}{\mu}{\lambda}{\gamma}
b_\dtrisub{P}{P}{\mu}{\mu}{\lambda}
\otimes
b_\dtrisub{P}{P}{\mu}{\mu}{\gamma}
$$
\item{(g)}  By left multiplication, the algebra $B\otimes_A B$ is
a left $B$-module.  If $\Rad(B\otimes_A B)$ is a $B$-submodule
of $B\otimes_A B$ and $\iota\colon B\to (B\otimes_A B)/\Rad(B\otimes_A B)$
is a left $B$-module homomorphism then
$$\iota\Big(
b_\dtrisub{R}{S}{\tau}{\beta}{\pi}
\Big)
=\sum_{\pi\to\gamma}
e_\ddtrisub{R}{S}{\tau}{\beta}{\pi}{\pi}{\gamma}.
$$
\endthm

\pf By (4.11) and (4.25),
$$
\matrix{ B &\mapright{\sim} &\displaystyle{
\bigoplus_{\mu\in \hat A} \overrightarrow A^\mu \otimes L^\mu } \cr
b_\dtrisub{P}{Q}{\mu}{\nu}{\lambda}
&\longmapsto
&\overrightarrow a_\vsub{P}{\mu}\otimes
\ell^\mu_\dtrisub{}{Q}{\mu}{\nu}{\lambda} \cr}
\qquad\hbox{and}\qquad
\matrix{ B &\mapright{\sim} &\displaystyle{
\bigoplus_{\nu\in \hat A} R^\nu\otimes \overleftarrow A^\nu} \cr
b_\dtrisub{P}{Q}{\mu}{\nu}{\lambda}
&\longmapsto
&r^\nu_\dtrisub{P}{}{\mu}{\nu}{\lambda}
\otimes\overleftarrow
a_\vsub{Q}{\nu} \cr
}
\formula$$
as left $A$-modules and as right $A$-modules, respectively.
Identify the left and right hand sides of these isomorphisms.
Then, by (4.17),  the elements of $C= B\otimes_A B$ given by
$$\bar m_\ddtrisub{P}{Q}{\mu}{\nu}{\lambda}{\sigma}{\gamma}
=
r^\gamma_\dtrisub{P}{}{\mu}{\gamma}{\lambda}
\otimes\overleftarrow a_\vsub{T}{\gamma}
\otimes\overrightarrow a_\vsub{T}{\gamma}
\otimes
\ell^\gamma_\dtrisub{}{Q}{\gamma}{\nu}{\sigma}
=
b_\dtrisub{P}{T}{\mu}{\gamma}{\lambda}
\otimes
b_\dtrisub{T}{Q}{\gamma}{\nu}{\sigma}
\formula
$$
do not depend on ${\scriptstyle{T\to}}\gamma\in \hat A^\gamma$ and
form a basis of $B\otimes_A B$.

\smallskip\noindent
(b) By (4.12), the map $\varepsilon\colon B\otimes_\FF B\to A$ is determined
by the values
$$
\varepsilon^\mu_\ddtrisub{T}{Q}{\gamma}{\tau}{\lambda}{\sigma}{\mu}
\in \FF
\qquad\hbox{given by}\qquad
\varepsilon^\mu_\ddtrisub{T}{Q}{\gamma}{\tau}{\lambda}{\sigma}{\mu}
a_\trisub{P}{P}{\mu}
=
\varepsilon\big(
\overrightarrow a_\vsub{P}{\mu}
\otimes
\ell^\mu_\dtrisub{}{T}{\mu}{\gamma}{\lambda}
\otimes
r^\mu_\dtrisub{Q}{}{\tau}{\mu}{\sigma}
\otimes\overleftarrow a_\vsub{P}{\mu}
\big).
\formula$$
Since
$$
\eqalign{
\varepsilon^\mu_\ddtrisub{T}{Q}{\gamma}{\tau}{\lambda}{\sigma}{\mu}
a_\trisub{P}{P}{\mu}
&=
\varepsilon\big(
b_\dtrisub{P}{T}{\mu}{\gamma}{\lambda}
\otimes
b_\dtrisub{Q}{P}{\tau}{\mu}{\sigma}
\big)
= \varepsilon_1\big(
b_\dtrisub{P}{T}{\mu}{\gamma}{\lambda}
\otimes
b_\dtrisub{Q}{P}{\tau}{\mu}{\sigma}
\big)  \cr
&= \delta_\dtrisub{T}{Q}{\gamma}{\tau}{\lambda\,\sigma}
\varepsilon_1\big(b_\dtrisub{P}{P}{\mu}{\mu}{\lambda}\big)
= \delta_\dtrisub{T}{Q}{\gamma}{\tau}{\lambda\,\sigma}
\varepsilon_1\big(b_\dtrisub{P}{P}{\mu}{\mu}{\lambda}
b_\dtrisub{P}{P}{\mu}{\mu}{\lambda} \big)
= \delta_\dtrisub{T}{Q}{\gamma}{\tau}{\lambda\,\sigma}
\varepsilon^\mu_\dtrisub{P}{P}{\mu\mu}{\mu}{\lambda\lambda}
a_\trisub{P}{P}{\mu}
.\cr
}$$
the matrix ${\cal E}^\mu$ given by (4.14) is diagonal
with entries $\varepsilon_\mu^\lambda$ given by (4.15)
and, by (4.17), $\varepsilon_\mu^\lambda$ is independent of
${\scriptstyle{P\to}}\mu\in \hat A^\mu$.
By Theorem 4.18(a),
$$
\bar m_\ddtrisub{P}{Q}{\mu}{\nu}{\lambda}{\sigma}{\gamma}
\bar m_\ddtrisub{R}{S}{\tau}{\xi}{\rho}{\eta}{\pi}
=\delta_{\gamma\pi}
\varepsilon_\ddtrisub{Q}{R}{\nu}{\tau}{\sigma}{\rho}{\gamma}
\bar m_\ddtrisub{P}{S}{\mu}{\xi}{\lambda}{\eta}{\gamma}
=
\delta_{\gamma\pi}
\delta_\dtrisub{Q}{R}{\nu}{\tau}{\sigma\,\rho}
\varepsilon_{\gamma}^\sigma
\bar m_\ddtrisub{P}{S}{\mu}{\xi}{\lambda}{\eta}{\gamma}
$$
in the algebra $C$. The rese of the statements in part (b)
follow from Theorem 4.18(b).

\smallskip\noindent
(c)
Evaluating the equations in (4.31) and using (4.29)
gives
$$\tr_B^\lambda = \tr_B( b_\dtrisub{P}{P}{\mu}{\mu}{\lambda})
= \tr_A( \varepsilon_1(
b_\dtrisub{P}{P}{\mu}{\mu}{\lambda}))
= \varepsilon_\mu^\lambda
\tr_A( a_\trisub{P}{P}{\mu})
= \varepsilon_\mu^\lambda
\tr_A^\mu,
\formula$$

\smallskip\noindent
(d) Since
$$1 =
\sum_{P\to\mu\to\lambda}
b_\dtrisub{P}{P}{\mu}{\mu}{\lambda}
\qquad\hbox{in the algebra $B$,}
$$
it follows from part (b) and (4.16) that
$$
1\otimes 1
= \Big(
\sum_{P\to\mu\to\lambda}
b_\dtrisub{P}{P}{\mu}{\mu}{\lambda}
\Big)
\otimes
\Big(
\sum_{Q\to\nu\to\gamma}
b_\dtrisub{Q}{Q}{\nu}{\nu}{\gamma}
\Big)
=
\sum_{P\to\mu\to\lambda \atop Q\to\nu\to\gamma}
\delta_{PQ}\delta_{\mu\nu}
\Big(
b_\dtrisub{P}{P}{\mu}{\mu}{\lambda}
\otimes
b_\dtrisub{Q}{Q}{\nu}{\nu}{\gamma}
\Big)
=
\sum_\starsub{P}{\mu}{\lambda}{\gamma}
\bar m_\ddtrisub{P}{P}{\mu}{\mu}{\lambda}{\gamma}{\mu}
$$
giving part (d).
\smallskip\noindent
(e) By left multiplication, the algebra $B\otimes_A B$ is a left $B$-module.
If $\varepsilon_\gamma^\lambda\ne 0$ and
$\varepsilon_\gamma^\sigma\ne 0$ then
$$
b_\dtrisub{R}{S}{\tau}{\beta}{\pi}
e_\ddtrisub{P}{Q}{\mu}{\nu}{\lambda}{\sigma}{\gamma}
=
\left({1\over \varepsilon_\gamma^\sigma}\right)
b_\dtrisub{R}{S}{\tau}{\beta}{\pi}
\Big(
b_\dtrisub{P}{T}{\mu}{\gamma}{\lambda}
\otimes
b_\dtrisub{T}{Q}{\gamma}{\nu}{\sigma}
\Big)
=
\left({1\over \varepsilon_\gamma^\sigma}\right)
\delta_\dtrisub{S}{P}{\beta}{\mu}{\pi\,\lambda}
\Big(
b_\dtrisub{R}{T}{\tau}{\gamma}{\lambda}
\otimes
b_\dtrisub{T}{Q}{\gamma}{\nu}{\sigma}
\Big)
=
\delta_\dtrisub{S}{P}{\beta}{\mu}{\pi\,\lambda}
e_\ddtrisub{R}{Q}{\tau}{\nu}{\pi}{\sigma}{\gamma}.
$$
Thus, if
$\iota\colon B\to (B\otimes_A B)/\Rad(B\otimes_A B)$
is a left $B$-module homomorphism then
$$
\hskip.5truein
\iota\Big(
b_\dtrisub{R}{S}{\tau}{\beta}{\pi}
\Big)
=\iota\Big(
b_\dtrisub{R}{S}{\tau}{\beta}{\pi}
\Big)\cdot 1
=
b_\dtrisub{R}{S}{\tau}{\beta}{\pi}
\sum_{P\to\mu\to\lambda\to\gamma}
e_\ddtrisub{P}{P}{\mu}{\mu}{\lambda}{\lambda}{\gamma}
=\sum_{P\to\mu\to\lambda\to\gamma}
\delta_\dtrisub{S}{P}{\beta}{\mu}{\pi\,\lambda}\,
e_\ddtrisub{R}{P}{\tau}{\mu}{\pi}{\lambda}{\gamma}
=\sum_{\pi\to\gamma}
e_\ddtrisub{R}{S}{\tau}{\beta}{\pi}{\pi}{\gamma}.\hskip.5truein\hbox{\qed}
$$

\section 5. Semisimple Algebras

\bigskip
Let $R$ be a integral domain and let $A_R$ be an algebra
over $R$, so that $A_R$ has an $R$-basis $\{b_1,\ldots, b_d\}$,
$$A_R = R\hbox{-span}\{b_1,\ldots, b_d\}
\qquad\hbox{and}\quad
b_ib_j = \sum_{k=1}^d r_{ij}^k b_k,\qquad\hbox{
with $r_{ij}^k\in R$,}$$
making $A_R$ a ring with identity.
Let $\FF$ be the field of fractions of $R$, let
$\bar \FF$ be the algebraic closure of $\FF$ and set
$$A = \bar \FF \otimes_R A_R = \bar \FF\hbox{-span}\{b_1,\ldots, b_d\},$$
with multiplication determined by the multiplication in $A_R$.
Then $A$ is an algebra over $\bar \FF$.

A {\it trace} on $A$ is a linear map $\vec t\colon A\to \bar \FF$ such that
$$\vec t(a_1a_2)=\vec t(a_2a_1), \qquad \hbox{for all $a_1,a_2\in A$.}$$
A trace $\vec t$ on $A$
is {\it nondegenerate} if for each $b\in A$ there is an $a\in A$ such that
$\vec t(ba)\ne 0$.

\lemma  Let $A$ be a finite dimensional algebra over a field $\FF$,
let $\vec t$ be a trace on $A$.  Define a symmetric bilinear form
$\langle,\rangle\colon A\times A\to \FF$ on $A$ by
$\langle a_1,a_2\rangle = \vec t(a_1a_2)$, for all $a_1,a_2\in A$.
Let $B$ be a basis of $A$.
Let $G = \big(\langle b,b'\rangle\big)_{b,b'\in B}$ be the matrix of the 
form
$\langle,\rangle$ with respect to $B$.  The following are equivalent:
\smallskip
\item{(1)} The trace $\vec t$ is nondegenerate.
\smallskip
\item{(2)} $\det G \ne 0$.
\smallskip
\item{(3)} The dual basis $B^*$ to the basis $B$ with respect to the form
$\langle,\rangle$ exists.
\pf
$(2)\Leftrightarrow (1)\colon$ \enspace
The trace $\vec t$ is degenerate if there is an element $a\in A$, $a\ne 0$,
such that $\vec t(ac)=0$ for all $c\in B$.  If $a_b\in \bar \FF$ are such 
that
$$a=\sum_{b\in B} a_b b,
\qquad\quad\hbox{then}\qquad\quad
0= \langle a, c\rangle = \sum_{b\in B} a_b\langle b,c\rangle$$
for all $c\in B$.  So $a$ exists if and only if
the columns of $G$ are linearly dependent, i.e. if and only if
$G$ is not invertible.

$(3)\Leftrightarrow (2)\colon$ \enspace
Let $B^*=\{b^*\}$ be the dual basis to $\{b\}$ with respect to
$\langle,\rangle$ and let $P$ be the change of basis matrix from $B$
to $B^*$.  Then
$$d^* = \sum_{b\in B} P_{db}b,
\qquad\hbox{and}\qquad
\delta_{bc}=\langle b, d^*\rangle = \sum_{d\in B} P_{dc}\langle b,c\rangle
=(GP^t)_{b,c}.$$
So $P^t$, the transpose of $P$, is the inverse of
the matrix $G$.  So the dual basis to $B$ exists
if and only if $G$ is invertible, i.e. if and only if $\det G\ne 0$.
\endpf

\prop
Let $A$ be an algebra and let $\vec t$ be a nondegenerate trace on $A$.
Define a symmetric bilinear form
$\langle,\rangle\colon A\times A\to \bar \FF$
on $A$ by $\langle a_1,a_2\rangle = \vec t(a_1a_2)$,
for all $a_1,a_2\in A$.
Let $B$ be a basis of $A$ and let $B^*$ be the dual basis to $B$
with respect to $\langle\ ,\ \rangle$.
\smallskip\noindent
(a) Let $a\in A$.  Then
$$[a] = \sum_{b\in B} bab^*
\qquad\hbox{is an element of the center $Z(A)$ of $A$}$$
and $[a]$ does not
depend on the choice of the basis $B$.
\smallskip\noindent
(b)  Let $M$ and $N$ be $A$-modules and let $\phi\in \Hom_{\bar \FF}(M,N)$ 
and
define
$$[\phi] = \sum_{b\in B} b\phi b^*.$$
Then $[\phi]\in \Hom_A(M,N)$ and $[\phi]$ does not depend on the choice
of the basis $B$.
\pf
(a) Let $c\in A$.  Then
$$
c[a] = \sum_{b\in B} cbab^*
= \sum_{b\in B} \sum_{d\in B} \langle cb,d^*\rangle dab^*
= \sum_{d\in B} da \sum_{b\in B} \langle d^*c,b\rangle b^*
= \sum_{d\in B} dad^*c
= [a]c,
$$
since $\langle cb,d^*\rangle = \vec t(cbd^*)=\vec t(d^*cb) =
\langle d^*c,b\rangle$.
So $[a]\in Z(A)$.

Let $D$ be another basis of $A$ and let $D^*$ be the dual basis to
$D$ with respect to $\langle, \rangle$.  Let $P=\big(P_{db}\big)$ be the
transition matrix from $D$ to $B$ and let
$P^{-1}$ be the inverse of $P$.  Then
$$d = \sum_{b\in B} P_{db}b
\qquad\hbox{and}\qquad
d^* = \sum_{\tilde b\in B} (P^{-1})_{\tilde bd}\tilde b^*,$$
since
$$
\langle d,\tilde d^*\rangle
=\left\langle\sum_{b\in B} P_{db}b,
\sum_{\tilde b\in B} (P^{-1})_{\tilde b\tilde d}\tilde b^*\right\rangle
=\sum_{b,\tilde b\in B} P_{db}(P^{-1})_{\tilde b\tilde d}
\delta_{b\tilde b}
= \delta_{d\tilde d}.
$$
So
$$
\sum_{d\in D} dad^* = \sum_{d\in D} \sum_{b\in B} P_{db}ba
\sum_{\tilde b\in B} (P^{-1})_{\tilde bd}\tilde b^*
= \sum_{b,\tilde b\in B} ba\tilde b^* \delta_{b\tilde b}
= \sum_{b\in B} bab^*.
$$
So $[a]$ does not depend on the choice of the basis $B$.

\smallskip\noindent
The proof of part (b) is the same as the proof of part (a)
except with $a$ replaced by $\phi$.
\endpf

\bigskip
Let $A$ be an algebra and let $M$ be an $A$-module.  Define
$$\End_A(M)=\{ T\in \End(M) \ |\ Ta=aT \hbox{\ for all $a\in A$}\}.$$

\thm (Schur's Lemma)
Let $A$ be a finite dimensional algebra over an algebraically closed
field $\bar \FF$.
\medskip
\item{(1)}  Let $A^\lambda$ be a simple $A$-module.  Then
$\End_A(A^\lambda) = \bar \FF\cdot \Id_{A^\lambda}.$
\medskip
\item{(2)}
If $A^\lambda$ and $A^\mu$ are nonisomorphic simple $A$-modules then
$\Hom_A(A^\lambda,A^\mu) = 0.$
\pf
Let $T\colon A^\lambda\to A^\mu$ be a nonzero $A$-module homomorphism.
Since $A^\lambda$ is simple, $\ker T = 0$ and so $T$ is injective.
Since $A^\mu$ is simple, $\im T = A^\mu$ and so $T$ is surjective.
So $T$ is an isomorphism.  Thus we may assume that
$T\colon A^\lambda\to A^\lambda$.

Since $\bar \FF$ is algebraically closed $T$ has an eigenvector and
a corresponding eigenvalue $\alpha\in \bar \FF$.  Then $T-\alpha\cdot Id\in
\Hom_A(A^\lambda,A^\lambda)$ and so $T-\alpha\cdot Id$ is either $0$
or an isomorphism.  However, since $\det(T-\alpha\cdot Id)=0$,
$T-\alpha\cdot Id$ is not invertible.  So $T-\alpha\cdot Id=0$.
So $T=\alpha\cdot Id$.  So $\End_A(A^\lambda) = \bar \FF\cdot Id$.
\endpf

\thm  (The Centralizer Theorem)
Let $A$ be a finite dimensional algebra over an
algebraically closed field $\bar \FF$.  Let $M$ be a semisimple
$A$-module and set $Z = \End_A(M)$.  Suppose that
$$M\cong \bigoplus_{\lambda\in \hat M} (A^\lambda)^{\oplus m_\lambda},$$
where $\hat M$ is an index set for the irreducible $A$-modules
$A^\lambda$ which appear in $M$ and the $m_\lambda$ are positive
integers.
\smallskip
\item{(a)}  $\displaystyle{
Z \cong \bigoplus_{\lambda\in \hat M} M_{m_\lambda}(\bar \FF) }$.
\smallskip
\item{(b)}  As an $(A, Z)$-bimodule
$$M\cong \bigoplus_{\lambda\in \hat M} A^\lambda \otimes Z^\lambda,$$
where the $Z^\lambda$, $\lambda\in \hat M$, are the simple $Z$-modules.
\pf
Index the components in the decomposition of $M$ by dummy variables
$\epsilon_i^\lambda$ so that we may write
$$M\cong \bigoplus_{\lambda\in \hat M} \bigoplus_{i=1}^{m_\lambda}
A^\lambda\otimes \epsilon_i^\lambda.$$
For each $\lambda\in \hat M$, $1\le i,j\le m_\lambda$ let
$\phi_{ij}^\lambda\colon A^\lambda\otimes \epsilon_j
\to A^\lambda\otimes \epsilon_i$ be the $A$-module isomorphism given by
$$\phi_{ij}^\lambda(m\otimes \epsilon_j^\lambda)
= m\otimes \epsilon_i^\lambda,
\qquad\hbox{for $m\in A^\lambda$.}$$
By Schur's Lemma,
$$\eqalign{
\End_A(M) = \Hom_A(M,M)
&\cong \Hom_A\left(
\bigoplus_\lambda \bigoplus_j
A^\lambda\otimes\epsilon_j^\lambda,
\bigoplus_\mu \bigoplus_i
A^\mu\otimes\epsilon_i^\mu\right) \cr
&\cong \bigoplus_{\lambda,\mu} \bigoplus_{i,j}
\delta_{\lambda\mu}
\Hom_A(A^\lambda\otimes \epsilon_j^\lambda,A^\mu\otimes \epsilon_i^\mu)
\cong \bigoplus_{\lambda}\bigoplus_{i,j=1}^{m_\lambda}
\bar \FF\phi^\lambda_{ij}. \cr
}$$
Thus each element $z\in \End_A(M)$ can be written as
$$z= \sum_{\lambda\in \hat M} \sum_{i,j=1}^{m_\lambda} z_{ij}^\lambda
\phi_{ij}^\lambda,
\qquad\hbox{for some $z_{ij}^\lambda\in \bar \FF$,}$$
and identified with an element of $~\bigoplus_\lambda M_{m_\lambda}(\bar 
\FF)$.
Since $\phi_{ij}^\lambda\phi_{kl}^\mu=\delta_{\lambda\mu}\delta_{jk}
\phi_{il}^\lambda$ it follows that
$$\End_A(M) \cong \bigoplus_{\lambda\in \hat M} M_{m_\lambda}(\bar \FF).$$

(b)  As a vector space
$Z^\mu = {\rm span}\{ \epsilon_i^\mu\ |\ 1\le i\le m_\mu\}$
is isomorphic to the simple $~\bigoplus_\lambda M_{m_\lambda}(\bar \FF)$
module of column vectors of length $m_\mu$.  The decomposition of
$M$ as $A\otimes Z$ modules follows since
$$(a\otimes \phi_{ij}^\lambda) (m\otimes \epsilon_k^\mu)
= \delta_{\lambda\mu}\delta_{jk} (a\otimes \epsilon_i^\mu),
\qquad\hbox{for all $m\in A^\mu$, $a\in A$,}
\qquad\hbox{\qed}$$

\medskip
If $A$ is an algebra then $A^{\rm op}$ is the algebra $A$ except with
the opposite multiplication, i.e.
$$A^{\rm op} = \{ a^{\rm op} \ |\ a\in A\} \qquad\hbox{with}\qquad
a_1^{\rm op}a_2^{\rm op} = (a_2a_1)^{\rm op},
\quad\hbox{for all $a_1,a_2\in A$.}$$
The left {\it regular representation} of $A$ is the vector space $A$
with $A$ action given by left multiplication.  Here $A$ is serving both as 
an
algebra and as an $A$-module.  It is often useful to distinguish the two
roles of $A$ and use the notation $\vec A$ for the $A$-module, i.e.
$\vec A$ is the vector space
$$\vec A =\{ \vec b\ |\ b\in A\}
\qquad\hbox{with $A$-action}\qquad
a\vec b = \overrightarrow{ab},
\qquad\hbox{for all $a\in A$, $\vec b\in \vec A$.}$$

\prop Let $A$ be an algebra and let $\vec A$ be the regular representation
of $A$.  Then $\End_A(\vec A)\cong A^{\rm op}$.  More precisely,
$$\End_A(\vec A) = \{ \phi_b\ |\ b\in A\}, \qquad\hbox{where
$\phi_b$ is given by}\qquad
\phi_b(\vec a) = \vec{ab}, \quad\hbox{for all $\vec a\in \vec A$.}$$
\pf
Let $\phi\in \End_A(\vec A)$ and let $b\in A$ be such that
$\phi(\vec 1)=\vec b$.  For all $\vec a\in \vec A$,
$$\phi(\vec a) = \phi(a\cdot \vec 1) = a\phi(\vec 1)
=a\vec b = \vec{ab},$$
and so $\phi=\phi_b$.
Then $\End_A(\vec A)\cong A^{\rm op}$ since
$$(\phi_{b_1}\circ \phi_{b_2})(\vec a)=\vec{ab_2b_1}
=\phi_{b_2b_1}(\vec a),$$
for all $b_1,b_2\in A$ and $\vec a\in \vec A$.
\endpf

\thm
Suppose that $A$ is a finite dimensional algebra over an algebraically
closed field $\FF$ such that the regular representation
$\vec A$ of $A$ is completely decomposable.
Then $A$ is isomorphic to a direct sum of matrix algebras, i.e.
$$A\cong \bigoplus_{\lambda\in \hat A} M_{d_\lambda}(\bar \FF),$$
for some set $\hat A$ and some positive integers $d_\lambda$, indexed
by the elements of $\hat A$.
\pf
If $\vec A$ is completely decomposable then, by Theorem 5.4,
$\End_A(\vec A)$ is isomorphic to a direct sum of matrix algebras.
By Proposition 5.5,
$$A^{\rm op} \cong
\bigoplus_{\lambda\in \hat A} M_{d_\lambda}(\bar \FF),$$
for some set $\hat A$ and some positive integers $d_\lambda$, indexed
by the elements of $\hat A$.
The map
$$\matrix{
\left(\bigoplus_{\lambda\in \hat A} M_{d_\lambda}(\bar \FF)\right)^{\rm op}
&\longrightarrow &
\bigoplus_{\lambda\in \hat A} M_{d_\lambda}(\bar \FF) \cr
a &\longmapsto &a^t, \cr
}$$
where $a^t$ is the transpose of the matrix $a$, is an algebra
isomorphism.  So $A$ is isomorphic to a direct sum of matrix algebras.
\endpf

If $A$ is an algebra then the trace $\tr$  of the regular representation
is the trace on $A$ given by
$$\tr(a) = \Tr (\vec A(a)), \qquad\hbox{for $a\in A$,}$$
where $\vec A(a)$ is the linear transformation of $A$ induced by
the action of $a$ on $A$ by left multiplication.

\prop Let
$A = \bigoplus_{\lambda\in \hat A} M_{d_\lambda}(\bar \FF).$
Then the trace $\tr$ of the regular representation of $A$ is nondegenerate.
\pf
As $A$-modules, the regular representation
$$\vec A\cong \bigoplus_{\lambda\in \hat A} (A^\lambda)^{\oplus 
d_\lambda},$$
where $A^\lambda$ is the irreducible $A$-module consisting
of column vectors of length $d_\lambda$.  For $a\in A$ let
$A^\lambda(a)$ be the linear transformation of $A^\lambda$ induced
by the action of $a$.  Then the trace $tr$ of the
regular representation is given by
$$\tr = \sum_{\lambda\in \hat A} d_\lambda\chi^\lambda,
\qquad\hbox{where}\qquad
\matrix{\chi^\lambda_A\colon &A &\to &\bar \FF \cr
&a &\longmapsto &\Tr(A^\lambda(a))\,,\cr
}$$
where $\chi^\lambda_A$ are the irreducible characters of $A$.  Since
the $d_\lambda$ are all nonzero the trace $tr$ is nondegenerate.
\endpf

\thm  (Maschke's theorem)
Let $A$ be a finite dimensional algebra over a field $\FF$
such that the trace $tr$ of the regular
representation of $A$ is nondegenerate.  Then every representation
of $A$ is completely decomposable.
\pf
Let $B$ be a basis of $A$ and let $B^*$ be the dual basis of $A$ with
respect to the form $\langle, \rangle\colon A\times A\to \bar \FF$ defined 
by
$$\langle a_1,a_2\rangle = tr(a_1a_2), \qquad
\hbox{for all $a_1,a_2\in A$.}$$
The dual basis $B^*$ exists because the trace $tr$ is nondegenerate.

Let $M$ be an $A$-module.  If $M$ is irreducible then the result is
vacuously true, so we may assume that $M$ has a proper submodule
$N$.  Let $p\in \End(M)$ be a projection onto $N$, i.e. $pM=N$ and
$p^2=p$.  Let
$$[p] = \sum_{b\in B} bpb^*,
\qquad\hbox{and}\qquad e = \sum_{b\in B} bb^*.$$
For all $a\in A$,
$$
\tr(ea)
= \sum_{b\in B} \tr(bb^*a)
= \sum_{b\in B} \langle ab,b^*\rangle
= \sum_{b\in B} ab\big|_b
= \tr(a),
$$
So $\tr((e-1)a)=0$, for all $a\in A$.  Thus, since $tr$ is nondegenerate,
$e=1$.

Let $m\in M$.  Then $pb^*m\in N$ for all $b\in B$,
and so $[p]m\in N$.  So $[p]M\subseteq N$.  Let $n\in N$.  Then
$pb^*n=b^*n$ for all $b\in B$, and so $[p]n = en = 1\cdot n=n$.
So $[p]M=N$ and $[p]^2 = [p]$, as elements of $\End(M)$.

Note that $[1-p] = [1]-[p] = e-[p] = 1-[p]$.
So
$$M = [p]M\oplus (1-[p])M = N \oplus [1-p]M,$$
and, by Proposition 5.2b, $[1-p]M$ is an $A$-module.
So $[1-p]M$ is an $A$-submodule of $M$ which is complementary
to $M$.  By induction on the dimension of $M$, $N$ and $[1-p]M$ are
completely decomposable, and therefore $M$ is completely decomposable.
\endpf

Together, Theorems 5.6, 5.8 and Proposition 5.7 yield the following theorem.

\thm (Artin-Wedderburn theorem)  Let $A$ be a finite dimensional algebra
over an algebraically closed field $\bar \FF$.
Let $\{b_1,\ldots, b_d\}$ be a basis of $A$ and let $\tr$ be the
trace of the regular representation of $A$.
The following are equivalent:
\smallskip
\itemitem{(1)}  Every representation of $A$ is completely decomposable.
\smallskip
\itemitem{(2)}  The regular representation of $A$ is
completely decomposable.
\smallskip
\itemitem{(3)}  $A\cong \bigoplus_{\lambda\in \hat A} M_{d_\lambda}(\bar 
\FF)$
for some finite index set $\hat A$, and some $d_\lambda\in \ZZ_{>0}$.
\smallskip
\itemitem{(4)}  The trace of the regular representation of $A$ is
nondegenerate.
\smallskip
\itemitem{(5)}
$\det(\tr(b_ib_j))\ne0$.
\endthm

\smallskip\noindent
{\bf Remark.}  Let $R$ be an integral domain, and let $A_R$
be an algebra over $R$ with basis $\{b_1,\ldots, b_d\}$.  Then
$\det(\tr(b_ib_j))$ is an element of $R$ and
$\det(\tr(b_ib_j))\ne 0$ in $\bar \FF$ if and only if
$\det(\tr(b_ib_j))\ne 0$ in $R$.
In particular, if $R = \CC[x]$, then $\det(\tr(b_ib_j))$ is a polynomial.
Since a polynomial has only a finite number of roots,
$\det(\tr(b_ib_j))(n)=0$ for only a finite number of values $n\in \CC$.

\thm (Tits deformation theorem)
Let $R$ be an integral domain, $\FF$, the field of fractions of $R$,
$\bar \FF$ the algebraic closure of $\FF$, and
$\bar R$, the integral closure of $R$ in $\bar \FF$.
Let $A_R$ be an $R$-algebra and let $\{b_1,\ldots, b_d\}$ be a basis
of $A_R$.  For $a\in A_R$ let
$\vec A(a)$ denote the linear transformation of $A_R$
induced by left multiplication by $a$.
Let $t_1,\ldots, t_d$ be indeterminates and let
$$\vec p(t_1,\ldots, t_d;x)
=\det(x\cdot\Id
-(t_1\vec A(b_1)+\cdots t_d\vec A(b_d)))\in R[t_1,\ldots, t_d][x],
$$
so that $\vec p$ is the characteristic polynomial of a ``generic''
element of $A_R$.
\smallskip\noindent
(a)  Let $A_{\bar\FF} = \bar \FF\otimes_R A_R$.  If
$$A_{\bar \FF} \cong \bigoplus_{\lambda\in \hat A} M_{d_\lambda}(\bar 
\FF),$$
then the factorization of $\vec p(t_1,\ldots, t_d,x)$
into irreducibles in $\bar \FF[t_1,\ldots, t_d,x]$ has the form
$$\vec p = \prod_{\lambda\in \hat A} (\vec p^\lambda)^{d_\lambda},
\qquad\hbox{with}\qquad
\vec p^\lambda\in \bar R[t_1,\ldots,t_d,x]
\qquad\hbox{and}\qquad
d_\lambda = {\rm deg}(\vec p^\lambda).$$
If $\chi^\lambda(t_1,\ldots, t_d)\in \bar R[t_1,\ldots, t_d]$ is
given by
$$\vec p^\lambda(t_1,\ldots, t_d,x)
= x^{d_\lambda} - \chi^\lambda(t_1,\ldots, t_d)x^{d_\lambda-1}
+\cdots,
$$
then
$$
\matrix{
\chi^\lambda_{A_{\bar \FF}}\colon &A_{\bar \FF} &\longmapsto &\bar \FF \cr
&\alpha_1b_1+\cdots+\alpha_d b_d &\longmapsto
&\chi^\lambda(\alpha_1,\ldots, \alpha_d)\,, \cr
}\qquad \lambda\in \hat A,$$
are the irreducible characters of $A_{\bar \FF}$.
\smallskip\noindent
(b)  Let $\KK$ be a field and let $\bar \KK$ be the algebraic closure
of $\KK$. Let $\gamma\colon R\to \KK$ be a ring
homomorphism and let $\bar \gamma\colon \bar R\to \bar \KK$ be the
extension of $\gamma$.  Let $\chi^\lambda(t_1,\ldots, t_d)\in
\bar R[t_1,\ldots, t_d]$ be as in (a).
If $A_{\bar \KK} = \bar \KK\otimes_R A_R$ is semisimple then
$$A_{\bar \KK} \cong \bigoplus_{\lambda\in \hat A} M_{d_\lambda}(\bar \KK),
\qquad\hbox{and}\qquad
\matrix{
\chi^\lambda_{A_{\bar \KK}}\colon &A_{\bar \KK} &\longmapsto &\bar \KK \cr
&\alpha_1b_1+\cdots+\alpha_d b_d &\longmapsto
&(\bar\gamma\chi^\lambda)(\alpha_1,\ldots, \alpha_d)\,, \cr
}
$$
for $\lambda\in \hat A$, are the irreducible characters of $A_{\bar \KK}$.
\pf
First note that if $\{b_1',\ldots, b_d'\}$ is another basis
of $A_R$ and the change of basis matrix $P=(P_{ij})$ is given by
$$b_i' = \sum_j P_{ij}b_j
\qquad\hbox{then the transformation}\qquad
t_i' = \sum_j P_{ij} t_j,$$
defines an isomorphism of polynomial rings $R[t_1,\ldots,t_d]\cong
R[t_1',\ldots,t_d']$.  Thus it follows that if the statements
are true for one basis of $A_R$ (or $A_{\bar \FF}$)
then they are true for every basis of $A_R$ (resp. $A_{\bar \FF}$).

\smallskip\noindent
(a) Using the decomposition of $A_{\bar \FF}$ let
$\{e_{ij}^\mu, \mu\in \hat A, 1\le i,j\le d_\lambda\}$ be a
basis of matrix units in $A_{\bar \FF}$ and let
$t_{ij}^\mu$ be corresponding variables.
Then the decomposition of $A_{\bar \FF}$ induces a factorization
$$\vec p(t_{ij}^\mu,x)
 = \prod_{\lambda\in \hat A}
(\vec p^\lambda)^{d_\lambda},
\qquad\hbox{where}\qquad
\vec p^\lambda(t_{ij}^\mu;x)
= \det(x-\sum_{\mu,i,j} t_{ij}^\mu A^\lambda(e_{ij})).
\formula$$
The polynomial $\vec p^\lambda(t_{ij}^\mu;x)$ is irreducible since
specializing the variables gives
$$\vec p^\lambda(t_{j+1,j}^\lambda=1, t_{1,n}^\lambda=t,
t_{i,j}^\mu=0 \hbox{\ otherwise}; x)
=x^{d_\lambda}-t,
\formula$$
which is irreducible in $\bar R[t;x]$.
This provides the factorization of $\vec p$ and
establishes that $\deg(\vec p^\lambda) = d_\lambda$.
By (5.11)
$$
\vec p^\lambda(t_{ij}^\mu;x)
= x^{d_\lambda}
-\Tr(A^\lambda(\sum_{\mu,i,j} t_{ij}^\mu e_{ij}^\mu))x^{d_\lambda-1}
+\cdots.$$
which establishes the last statement.

Any root of
$\vec p(t_1,\ldots, t_d,x)$ is an element of
$\overline{R[t_1,\ldots, t_d]}
=\bar R[t_1,\ldots, t_d]$.
So any root of
$\vec p^\lambda(t_1,\ldots, t_d,x)$ is an element of
$\bar R[t_1,\ldots, t_d]$ and therefore the coefficients of
$\vec p^\lambda(t_1,\ldots, t_d,x)$ (symmetric functions in
the roots of $\vec p^\lambda$) are elements of
$\bar R[t_1,\ldots, t_d]$.

\smallskip\noindent
(b)  Taking the image of the equation (5.11) give a factorization of
$\gamma(\vec p)$,
$$\gamma(\vec p)
= \prod_{\lambda\in \hat A} \gamma(\vec p^\lambda)^{d_\lambda},
\qquad\hbox{in $\bar \KK[t_1,\ldots, t_d,x]$}.$$
For the same reason as in (5.12) the factors
$\gamma(\vec p^\lambda)$ are irreducible polynomials in
$\bar \KK[t_1,\ldots, t_d,x]$.

On the other hand, as in the proof of (a), the decomposition of
$A_{\bar \KK}$ induces a factorization of
$\gamma(\vec p)$ into irreducibles in
$\bar \KK[t_1,\ldots, t_d,x]$.
These two factorizations must coincide, whence the result.
\endpf

Applying the Tits deformation theorem to the case where
$R = \CC[x]$ (so that $\FF = \CC(x)$)
gives the following theorem.  The statement in (a) is a consequence
of Theorem 5.6 and the remark which follows Theorem 5.9.

\thm
Let $\CC A(n)$ be a family of algebras defined by generators
and relations such that the coefficients of the relations are
polynomials in $n$.  Assume that there is an $\alpha\in \CC$ such that
$\CC A(\alpha)$ is semisimple.  Let $\hat A$ be an index set
for the irreducible $\CC A(\alpha)$-modules $A^\lambda(\alpha)$. Then
\item{(a)} $\CC A(n)$ is semisimple for all but a finite number
of $n\in \CC$.
\item{(b)}  If $n\in \CC$ is such that $\CC A(n)$ is semisimple
then $\hat A$ is an index set for the simple $\CC A(n)$-modules
$A^\lambda(n)$ and $\dim(A^\lambda(n))=\dim(A^\lambda(\alpha))$
for each $\lambda\in \hat A$.
\item{(c)}  Let $x$ be an indeterminate and let
$\{b_1,\ldots, b_d\}$ be a basis of $\CC[x] A(x)$.
Then there are polynomials
$\chi^\lambda(t_1,\ldots, t_d)\in \CC[t_1,\ldots, t_d,x]$,
$\lambda\in \hat A$,
such that for every $n\in \CC$ such that $\CC A(n)$ is semisimple,
$$\matrix{
\chi^\lambda_{A(n)}\colon &\CC A(n) &\longrightarrow &\CC \cr
&\alpha_1b_1+\cdots +\alpha_d b_d
&\longmapsto &\chi^\lambda(\alpha_1,\ldots,\alpha_d,n)\,, \cr}
\qquad \lambda\in \hat A,$$
are the irreducible characters of $\CC A(n)$.
\endthm

\section 6. References

\bigskip\noindent

\medskip
\item{[Bou1]} {\smallcaps N.\ Bourbaki},
{\it Algebra}, Chapters 1--3, Elements of Mathematics,
Springer-Verlag, Berlin, 1990.

\medskip
\item{[Bou2]} {\smallcaps N.\ Bourbaki},
{\it Groupes et Alg\`ebres de Lie,
Chapitre IV, V, VI\/}, El\'ements de Math\'ematique, Hermann, Paris (1968).

\medskip\noindent
\item{[Bra]}
{\smallcaps R.\ Brauer},
{\rm  On algebras which are connected with the semisimple continuous 
groups},
{\it Ann.\ Math.} {\bf 38}(1937), 854--872.

\medskip\noindent
\item{[Bro1]}
{\smallcaps W.P.\ Brown},
{\rm  An algebra related to the orthogonal group},
{\it Michigan Math J.\/} {\bf 3}(1955), 1--22.

\medskip\noindent
\item{[Bro2]}
{\smallcaps W.P.\ Brown},
{\rm  Generalized matrix algebras},
{\it Canadian J.\ Math.} {\bf 7}(1955), 188--190.

\medskip\noindent
\item{[CPS]}
{\smallcaps E.\ Cline, B.\ Parshall and L.\ Scott},
{\it Finite dimensional algebras and highest weight categories},
J.\ reine angew.\ Math.\ {\bf 391} (1988), 85-99.


\medskip\noindent
\item{[CR]}
{\smallcaps C.\ Curtis and I.\ Reiner},
{\it Methods of representation theory---With applications to 
finite groups and orders, Vols.~I and II,} Pure and Applied Mathematics. 
Wiley \& Sons, Inc., New York, 1987.

\medskip\noindent
\item{[DR]} 
{\smallcaps V.\ Dlab and C.\ Ringel},
{\rm A construction for quasi-hereditary algebras},
{\it Compositio Math.\ } {\bf 70} (1989), 155-175.

\medskip\noindent
\item{[DW]}
{\smallcaps W.\ Doran and D.\ Wales, David},
{\rm The partition algebra revisited.}
{\it J. Algebra} {\bf 231}(2000),  265--330.

\medskip\noindent
\item{[FH]} {\smallcaps J.\ Farina and T.\ Halverson},
{\rm Character orthogonality for the partition algebra and fixed
points of permutations},
{\it Adv.\ Applied Math.,} to appear.

\medskip\noindent
\item{[FL]} {\smallcaps D.\ Fitzgerald and J.\ Leech},
Dual symmetric inverse monoids and representation theory,
{\it J.\ Australian Math.\ Soc.\ Ser.\ A} {\bf 64}(1998), 345--367

\medskip\noindent
\item{[GHJ]}
{\smallcaps F.\ Goodman, P.\ de  la Harpe, and V.\ Jones, }
{\it  Coxeter graphs and towers of algebras},
Mathematical Sciences Research Institute Publications
{\bf 14} Springer-Verlag, New York, 1989.

\medskip\noindent
\item{[GW]}
{\smallcaps R.\ Goodman and N.\ Wallach,} {\it  Representations and 
invariants of the
classical groups}. Encyclopedia of Mathematics and its Applications {\bf 
68},
Cambridge University Press, Cambridge, 1998.

\medskip\noindent
\item{[Gr]} {\smallcaps C.\ Grood}, {\rm The rook partition algebra},
in preparation.

\medskip\noindent
\item{[Ha]} {\smallcaps T.\ Halverson},
{\rm Characters of the partition algebras},
{\it J. Algebra} {\bf 238}(2001),  502--533.

\medskip\noindent
\item{[Jo]} {\smallcaps V.\ F.\ R.\ Jones},
{\rm The Potts model and the symmetric group},
in {\it Subfactors: Proceedings of the Taniguchi Symposium
on Operator Algebras (Kyuzeso, 1993),}
World Sci.\ Publishing, River Edge, NJ, 1994, 259--267.

\medskip\noindent
\item{[Mac]}
{\smallcaps I.\ Macdonald},
{\it Symmetric Functions and Hall Polynomials},
2nd ed., Oxford Univ.\ Press, New York, 1995.

\medskip\noindent
\item{[Ma1]} {\smallcaps P.\ Martin},
Potts models and related problems in statistical mechanics,
{\it Series on Advances in Statistical Mechanics}, 5.
World Scientific Publishing Co., Inc., Teaneck, NJ, 1991.

\medskip\noindent
\item{[Ma2]}
{\smallcaps P.\ Martin},
{\rm Temperley-Lieb algebras for nonplanar statistical
mechanics---the partition algebra construction},
{\it J.\ Knot Theory Ramifications} {\bf 3}(1994), 51--82.

\medskip\noindent
\item{[Ma3]}
{\smallcaps P.\ Martin},
{\rm The structure of the partition algebras},
{\it J.\ Algebra} {\bf 183}(1996), 319--358.

\medskip\noindent
\item{[Ma4]}
{\smallcaps P.\ Martin},
{\rm The partition algebra and the Potts model transfer matrix spectrum in 
high
dimensions}, {\it J.\ Phys.\ A:Math.\ Gen.\ } {\bf 33} (2000), 3669--3695.

\medskip\noindent
\item{[MR]}
{\smallcaps P.\ Martin and G.\ Rollet},
{\rm The Potts model representation and a Robinson-Schensted
correspondence for the partition algebra},
{\it Compositio Math.} {\bf 112}(1998), 237--254.

\medskip\noindent
\item{[MS]}
{\smallcaps P.\ Martin and H.\ Saleur},
{\rm Algebras in higher-dimensional statistical mechanics---the
exceptional partition (mean field) algebras},
{\it Lett.\  Math.\ Phys.} {\bf 30}(1994), 179--185.

\medskip\noindent
\item{[MW1]}
{\smallcaps P.\ Martin and D.\ Woodcock},
{\rm On central idempotents in the partition algebra.}
{\it J. Algebra} {\bf 217}(1999), no. 1, 156--169.

\medskip\noindent
\item{[MW2]}
{\smallcaps P.\ Martin and D.\ Woodcock},
{\rm The partition algebras and a new deformation of the Schur algebras},
{\it J. Algebra} {\bf 203}(1998), 91--124.

\medskip\noindent
\item{[MV]}
{\smallcaps R.\ Mirollo and K.\ Vilonen},
{\rm Berstein-Gelfand-Gelfand reciprocity on perverse sheaves},
{\it Ann.\ Scient.\ \'Ec.\ Norm.\ Sup.\ } $4^{\rm e}$ s\'erie {\bf 20}
(1987), 311-324.

\medskip\noindent
\item{[Ow]}
{\smallcaps W.\ Owens},
{\ The Partition Algebra}, Honors Thesis, Macalester College,
 May 2002.

\medskip\noindent
\item{[Sta]}
{\smallcaps R.\ Stanley},
{\it Enumerative Combinatorics}, Vol.\ 1,
Cambridge Studies in Adv. Math. {\bf 49} 1997.

\medskip\noindent
\item{[Wz]}
{\smallcaps H.\ Wenzl},
{\rm On the structure of Brauer's centralizer algebras},
{\it Ann.\ Math.} {\bf 128}(1988), 173-193.

\medskip\noindent
\item{[Xi]}
{\smallcaps C.\ Xi},
{\rm Partition algebras are cellular},
{\it Compositio Math.} {\bf 119}(1999), no. 1, 99--109.

\vfill\eject
\end